\documentclass{svmult}
\usepackage{bbm}
\usepackage{amsfonts}
\usepackage{mathrsfs}
\usepackage{amssymb}
\newenvironment{enumerate-roman}{\begin{enumerate}}{\end{enumerate}}
\usepackage{makeidx}         
\usepackage{graphicx}        
\usepackage{multicol}        
\usepackage[bottom]{footmisc}

\makeindex             

\newtheorem{thm}{Theorem}[section]
\newtheorem{lem}[thm]{Lemma}
\newtheorem{prop}[thm]{Proposition}

\newtheorem{rmk}[thm]{Remark}

\newtheorem{defi}[thm]{Definition}

\begin{document}
\baselineskip13pt

\title*{Stationary Solutions of SPDEs and Infinite Horizon BDSDEs}
\titlerunning{Stationary Solution of SPDEs}
\author{Qi Zhang\inst{1,2},
Huaizhong Zhao\inst{1}}
\authorrunning{Q. Zhang and H.Z. Zhao}
\institute{ Department of Mathematical Sciences, Loughborough
University, Loughborough, LE11 3TU, UK.
\texttt{Q.Zhang3@lboro.ac.uk}, \texttt{H.Zhao@lboro.ac.uk} \and
School of Mathematics and System Sciences, Shandong University,
Jinan, 250100, China.}
\maketitle
\newcounter{bean}
\begin{abstract}

  In this paper we study the existence of
  stationary solutions for stochastic partial differential
  equations. We establish a new connection between $L_{\rho}^2({\mathbb{R}^{d}};{\mathbb{R}^{1}})\otimes L_{\rho}^2({\mathbb{R}^{d}};{\mathbb{R}^{d}})$ valued solutions of backward doubly
    stochastic differential equations (BDSDEs) on infinite horizon and
    the stationary solutions of the SPDEs. Moreover, we prove the existence and
    uniqueness of the solutions of BDSDEs on both finite and infinite horizons, so obtain
   the solutions of initial value problems and the stationary solutions (independent of any
   initial value) of SPDEs. The connection of the weak solutions of SPDEs and
    BDSDEs has independent interests in the areas of both SPDEs and BSDEs.

    \end{abstract}
 \textbf{Keywords:} backward doubly stochastic differential
equations, weak solutions, stochastic partial differential
equations, stationary solution, random dynamical systems. \vskip5pt

\noindent {AMS 2000 subject classifications}: 60H15, 60H10, 37H10.
\vskip5pt

\renewcommand{\theequation}{\arabic{section}.\arabic{equation}}

\section{Introduction}

Let $u: [0,\infty)\times U\times \Omega \to U$ be a measurable
random dynamical system on a measurable space $(U,\mathcal{B})$ over
a metric dynamical system ($\Omega$, $\mathscr{F}$, $P$,
$(\theta_t)_{t\geq0})$, then a stationary solution is a
$\mathscr{F}$ measurable random variable $Y:\Omega \to U$ such that
(Arnold \cite{ar})
\begin{eqnarray}\label{zhao001}
u(t,Y(\omega),\omega)=Y(\theta_t\omega)\ \ \ {\rm for}\ {\rm all}\
t\geq 0\ \rm{a.s.}.
\end{eqnarray}
This ``one-force, one-solution" setting is a natural extension of
equilibria or fixed points in deterministic systems to stochastic
counterparts. The simplest nontrivial example is the
Ornstein-Uhlenbeck process defined by the stochastic differential
equation $du(t)=-u(t)dt+dB_t$. It defines a random dynamical system
$u(t,u_0)=u_0{\rm e}^{-t}+\int_0^t {\rm e}^{-(t-s)}dB_s$ and its
stationary point is given by $Y(\omega)=\int_{-\infty}^0 {\rm
e}^{s}dB_s$. Moreover, for any $u_0$, $u(t,u_0,\theta_{-t}\omega)\to
Y(\omega)$ as $t\to \infty$, where $\theta_t$ is the shift operator
of the Brownian path: $(\theta_tB)(s)=B(t+s)-B(s)$ for any $s\in
(-\infty,+\infty)$. A pathwise stationary solution describes the
pathwise invariance of the stationary solution over time along the
measurable and $P$-preserving transformation $\theta_t$:
$\Omega\longrightarrow\Omega$, and the pathwise limit of the
solutions of random dynamical systems.
Needless to say, it is one of the fundamental questions of basic
importance (\cite{ar}, \cite{kloeden}, \cite{kh-ma-si},
\cite{mo-zh-zh}, \cite{si1}, \cite{si2}). For random dynamical systems generated 
by stochastic partial differential equations (SPDEs),
such random fixed points
consist of infinitely many random moving invariant
surfaces on the configuration space due to the random external force
pumped to the system constantly. They are more realistic models than
many deterministic models as it demonstrates some complicated
phenomena such as turbulence. Their existence and stability are of
great interests in both mathematics and physics. However, in
contrast to the deterministic dynamical systems, also due to the
fact that the external random force exists at all time, the
existence of stationary solutions of stochastic dynamical systems
generated e.g. by stochastic differential equations (SDEs) or SPDEs,
is a difficult and subtle problem. We would like to point out that
there have been extensive works on stability and invariant
manifolds of random dynamical systems, and researchers usually assume
there is an invariant set (or a single point: a stationary solution
or a fixed point, often assumed to be $0$), then prove invariant
manifolds and stability results at a point of the invariant set
(Arnold \cite{ar} and references therein, Ruelle \cite{ruelle},
Duan, Lu and Schaumulfuss \cite{du-lu-sc1}, Li and Lu \cite{li-lu},
Mohammed, Zhang and Zhao \cite{mo-zh-zh} to name but a few). But the
invariant manifolds theory gives neither the existence results of the
invariant set and the stationary solution nor a way to find them. In
particular, for  the existence of stationary solutions for SPDEs,
results are only known in very few cases (\cite{kloeden},
\cite{kh-ma-si}, \cite{mo-zh-zh}, \cite{si1}, \cite{si2}). In
\cite{si1}, \cite{si2}, the stationary strong solution of the
stochastic Burgers' equations with periodic or random forcing ($C^3$
in the space variable) was established by Sinai using the Hopf-Cole
transformation. In \cite{mo-zh-zh}, the stationary solution of the
stochastic evolution equations was identified as a solution of the
corresponding integral equation up to time $+\infty$ and the
existence was established for certain SPDEs by Mohammed, Zhang and
Zhao. But the existence of solutions of such a stochastic integral
equations in general is far from clear.

The main purpose of this paper is to find the pathwise stationary solution of
the following SPDE
\begin{eqnarray}\label{zhao21}
dv(t,x)&=&[\mathscr{L}v(t,x)+f\big(x,v(t,x),\sigma^*(x)Dv(t,x)\big)]dt\nonumber\\
&&+g\big(x,v(t,x), \sigma^*(x)Dv(t,x)\big)d{B}_t,
\end{eqnarray}
without assumption that there is an invariant set.
Here ${B}$ is a two-sided cylindrical Brownian motion on a separable
Hilbert space $U_0$; $\mathscr{L}$ is the infinitesimal generator of
a diffusion process $X_{s}^{t,x}$ (solution of Eq.(\ref{qi17}))
given by
\begin{eqnarray}\label{buchong1}
\mathscr{L}={1\over2}\sum_{i,j=1}^na_{ij}(x){{\partial^2}\over{\partial
x_i\partial x_j}}+\sum_{i=1}^nb_i(x){\partial\over{\partial x_i}}
\end{eqnarray}
with $\big(a_{ij}(x)\big)=\sigma\sigma^*(x)$. Eq. (\ref{zhao21}) is
very general, especially the nonlinear functions $f$ and $g$ can
include $\nabla u$ and the second order differential operator
${\mathscr {L}}$ is allowed to be degenerate, while in most literature, $g$ is not
allowed to depend on $\nabla u$ or $g$ only depends on $\nabla u$
linearly (Da Prato and Zabzcyk \cite{pr-za1}, Krylov \cite{krylov},
Pardoux \cite{pa2}).
 As an intermediate step, the result of
existence and uniqueness of the weak solutions of (\ref{zhao21}),
obtained by solving the corresponding backward doubly stochastic 
differential equations (BDSDEs), appears also new. The
existence and uniqueness of such equations when $g$ is independent
of $\nabla u$ or linearly dependent of $\nabla u$ were studied by Da
Prato and Zabzcyk \cite{pr-za1}, Krylov \cite{krylov}. But we don't
claim here our results on the existence and uniqueness for the types
of SPDEs studied in \cite{pr-za1} and \cite{krylov} have superseded
their previous results.

Note that from the pathwise stationary solution obtained in this paper,
 we can construct an invariant
measure for the skew product of the metric dynamical system and the
random dynamical system. In this connection, we mention that in
 recent years, substantial results on the
existence and uniqueness of invariant measures for SPDEs and weak
convergence of the law of the solutions as time tends to infinity
have been proved for many important SPDEs (\cite{br-ga},
\cite{br-li}, \cite{pr-za1}, \cite{fl}, \cite{ha-ma} to name but a
few). The invariant measure describes the invariance of a certain
solution in law when time changes, therefore it is a stationary
measure of the Markov transition probability. It is well known that
an invariant measure gives a stationary solution when it is a random
Dirac measure. Although an invariant measure of a random dynamical
system on ${\mathbb{R}^{1}}$ gives a stationary solution, in
general, this is not true unless one considers an extended
probability space. However, considering the extended probability
space, one essentially regards the random dynamical system as noise as well,
so the dynamics is different. See \cite{li-zh} for some examples of
SDEs on ${\mathbb{R}^{1}}$ and a perfect cocycle on
${\mathbb{S}^{1}}$ having an invariant measure, but not a stationary
solution. In fact, the stationary solution we study in this paper
gives the support of the corresponding invariant measure, so reveals
more detailed information than an invariant measure.

In this paper, BDSDEs will be used as our tool to study stationary
solutions of SPDEs. We will prove that the solutions of the
corresponding infinite horizon BDSDEs give the desired stationary
solutions of the SPDEs (\ref{zhao21}). Backward stochastic
differential equations (BSDEs) have been studied extensively in the
last 16 years since the pioneering work of Pardoux and Peng
\cite{pa-pe1}. The connection between BSDEs and quasilinear
parabolic partial differential equations (PDEs) was discovered by
Pardoux and Peng in \cite{pa-pe2} and Peng in \cite{pe}. The study
of the connection of weak solutions of PDEs and BSDEs began in
Barles and Lesigne \cite{ba-le}. The BDSDEs and their connections
with the SPDEs were studied by Pardoux and Peng in \cite{pa-pe3} for
the strong solutions, and by Bally and Matoussi in \cite{ba-ma} for
the weak solutions. On the other hand, the infinite horizon BSDE was
first studied by Peng in \cite{pe} and it was shown that the
corresponding PDE is a Poisson equation (elliptic equation). This
was studied systematically by Pardoux in \cite{pa}. Notice that the
solutions of the Poisson equations can be regarded as the stationary
solutions of the parabolic PDEs. Deepening this idea, it would not
be unreasonable to conjecture that the solutions of infinite horizon
BDSDEs (if exists) be the stationary solutions of the corresponding
SPDEs. Of course, we cannot write them as solutions of Poisson
equations or stochastic Poisson equations like in the deterministic
cases. However, it is very natural to describe the stationary
solutions of SPDEs by the solutions of infinite horizon BDSDEs. In
this sense, BDSDEs (or BSDEs) can be regarded as more general SPDEs
(or PDEs).

As far as we know, the connection of the pathwise stationary
solutions of the SPDEs and infinite horizon BDSDEs we study in this
paper is new (section 2). We believe this new method can be used to
many SPDEs such as those with quadratic or polynomial growth
nonlinear terms. We don't intend to include all these results in the
present paper, but only study Lipschitz continuous nonlinear term to
initiate this intrinsic method to the study of this basic problem in
dynamics of SPDEs. We would like to point out that our BDSDE method
depends on neither the continuity of the random dynamical system
(continuity means $u(t,\cdot,\omega):U\to U$ is a.s. continuous) nor
on the method of the random attractors. The continuity problem for
the SPDE (\ref{zhao21}) with the nonlinear noise considered in this
paper still remains open mainly due to the failure of Kolmogorov's
continuity theorem in infinite dimensional setting as pointed out by
some researchers (e.g. \cite{du-lu-sc1}, \cite{mo-zh-zh}).

One of the necessary intermediate steps is to study the BDSDEs on
finite horizon and establish their connections with the weak
solutions of SPDEs (Sections 3 and 4). Our method to study the
$L_{\rho}^2({\mathbb{R}^{d}};{\mathbb{R}^{1}})\otimes
L_{\rho}^2({\mathbb{R}^{d}};{\mathbb{R}^{d}})$ valued solutions of
BDSDEs on finite horizon was inspired by Bally and Matoussi's
approach on the existence and uniqueness of solutions of BDSDEs with
finite dimensional Brownian motions (\cite{ba-ma}). But our results are stronger
and our conditions are weaker. We will solve the BDSDEs driven by
the cylindrical Brownian motion and nonlinear terms satisfying
Lipschtz conditions in the space
$L_{\rho}^2({\mathbb{R}^{d}};{\mathbb{R}^{1}})\otimes
L_{\rho}^2({\mathbb{R}^{d}};{\mathbb{R}^{d}})$. We obtain a unique
solution $(Y_.^{t,\cdot},Z_.^{t,\cdot})\in
S^{2,0}([t,T];L_{\rho}^2({\mathbb{R}^{d}};{\mathbb{R}^{1}}))\otimes
M^{2,0}([t,T];L_{\rho}^2({\mathbb{R}^{d}};{\mathbb{R}^{d}}))$. The
result $Y_.^{t,\cdot}\in
S^{2,0}([t,T];L_{\rho}^2({\mathbb{R}^{d}};{\mathbb{R}^{1}})$, which
plays an important role in solving the nonlinear BDSDEs and proving
the connection with the weak solutions of SPDEs (also BSDEs and
PDEs), was not obtained in \cite{ba-ma}.
The generalized equivalence of norm principle (Section 2), which is a simple extension of the
equivalence of norm principle obtained by Kunita (\cite{ku2}), Barles and Lesigne (\cite{ba-le}),
Bally and Matoussi (\cite{ba-ma}) to random functions, also plays an important role in the proofs
of our results. We believe
our results for finite time BDSDEs are new even for BSDEs.

In section 5, we will solve the BDSDEs on infinite
horizon and in section 6, we study continuity of the solution in
order to ensure that it gives the perfect stationary solutions of
the SPDEs.

\section{The stationarity of the solutions of infinite horizon BDSDEs and stationary solutions of SPDEs}
\setcounter{equation}{0}

On a probability space $(\Omega, {\mathscr{F}}, P)$, let
$(\hat{B}_t)_{t\geq0}$ and $(W_t)_{t\geq0}$ be two mutually
independent $Q$-Wiener process valued on $U$ and a standard Brownian
motion valued on ${\mathbb{R}^{d}}$ respectively. Here $U$ is a
separable Hilbert space with countable base
$\{e_i\}_{i=1}^{\infty}$; $Q\in L(U)$ is a symmetric nonnegative
trace class operator such that $Qe_i=\lambda _ie_i$ and $\sum
\limits _{i=1}^{\infty}\lambda _i<\infty$. It is well known that
$\hat B$ has the following expansion (\cite{pr-za1}): for each $t$
\begin{eqnarray}\label{bmexpansion}
\hat B_t=\sum\limits _{j=1}^{\infty} \sqrt {\lambda _j}\hat \beta
_j(t)e_j,
\end{eqnarray}
where
\begin{eqnarray*}
\hat \beta _j(t)={1\over \sqrt{\lambda _j}}<\hat B_t,e_j>_U, \ \
j=1,2,\cdots
\end{eqnarray*}
are mutually independent real-valued Brownian motion on $(\Omega,
{\mathscr{F}}, P)$ and the series (\ref{bmexpansion}) is convergent
in $L^2(\Omega, {\mathscr{F}}, P)$. Let ${\cal N}$ denote the class
of $P$-null sets of ${\mathscr{F}}$. We define
\begin{eqnarray*}
&&\mathscr{F}_{t,T}\triangleq{\mathscr{F}_{t,T}^{\hat{B}}}\otimes
\mathscr{F}_t^W\bigvee{\cal N},\ \ \ {\rm for}\ 0\leq t\leq T;\\
&&\mathscr{F}_t\triangleq{\mathscr{F}_{t,\infty}^{\hat{B}}}\otimes
\mathscr{F}_t^W\bigvee{\cal N},\ \ \ \ \ {\rm for}\ t\geq0.
\end{eqnarray*}
Here for any process $(\eta_t)_{t\geq0}$,
$\mathscr{F}_{s,t}^\eta=\sigma\{\eta_r-\eta_s$; ${0\leq s\leq r\leq
t}\}$, ${\mathscr{F}}_t^\eta={\mathscr{F}}_{0,t}^\eta$,
$\mathscr{F}_{t,\infty}^{\eta}=\bigvee_{T\geq0}{\mathscr{F}_{t,T}^\eta}$.
\begin{defi}\label{qi00}
Let $\mathbb{S}$ be a Hilbert space with norm $\|\cdot\|_\mathbb{S}$
and Borel $\sigma$-field $\mathscr{S}$. For $K\in\mathbb{R}^+$, we
denote by $M^{2,-K}([0,\infty);\mathbb{S})$ the set of
$\mathscr{B}_{\mathbb{R}^+}\otimes\mathscr{F}/\mathscr{S}$
measurable random processes $\{\phi(s)\}_{s\geq0}$ with values on
$\mathbb{S}$ satisfying
\begin{enumerate-roman}
\item $\phi(s):\Omega\rightarrow\mathbb{S}$ is $\mathscr{F}_s$ measurable for $s\geq 0$;
\item $E[\int_{0}^{\infty}{\rm e}^{-Ks}\|\phi(s)\|_\mathbb{S}^2ds]<\infty$.
\end{enumerate-roman}
Also we denote by $S^{2,-K}([0,\infty);\mathbb{S})$ the set of
$\mathscr{B}_{\mathbb{R}^+}\otimes\mathscr{F}/\mathscr{S}$
measurable random processes $\{\psi(s)\}_{s\geq0}$ with values on
$\mathbb{S}$ satisfying
\begin{enumerate-roman}
\item $\psi(s):\Omega\rightarrow\mathbb{S}$ is $\mathscr{F}_s$ measurable for $s\geq0$ and $\psi(\cdot,\omega)$ is continuous $P$-a.s.;
\item $E[\sup_{s\geq0}{\rm
e}^{-Ks}\|\psi(s)\|_\mathbb{S}^2]<\infty$.
\end{enumerate-roman}
\end{defi}

Similarly, for $0\leq t\leq T<\infty$, we define
$M^{2,0}([t,T];\mathbb{S})$ and $S^{2,0}([t,T];\mathbb{S})$ on
finite time interval.
\begin{defi}\label{zhao005}
Let $\mathbb{S}$ be a Hilbert space with norm $\|\cdot\|_\mathbb{S}$
and Borel $\sigma$-field $\mathscr{S}$. We denote by
$M^{2,0}([t,T];\mathbb{S})$ the set of
$\mathscr{B}_{[t,T]}\otimes\mathscr{F}/\mathscr{S}$ measurable
random processes $\{\phi(s)\}_{t\leq s\leq T}$ with values on
$\mathbb{S}$ satisfying
\begin{enumerate-roman}
\item $\phi(s):\Omega\rightarrow\mathbb{S}$ is $\mathscr{F}_{s,T}\bigvee{\mathscr{F}_{T,\infty}^{\hat{B}}}$ measurable for $t\leq s\leq T$;
\item $E[\int_{t}^{T}\|\phi(s)\|_\mathbb{S}^2ds]<\infty$.
\end{enumerate-roman}
Also we denote by $S^{2,0}([t,T];\mathbb{S})$ the set of
$\mathscr{B}_{[t,T]}\otimes\mathscr{F}/\mathscr{S}$ measurable
random processes $\{\psi(s)\}_{t\leq s\leq T}$ with values on
$\mathbb{S}$ satisfying
\begin{enumerate-roman}
\item $\psi(s):\Omega\rightarrow\mathbb{S}$ is
$\mathscr{F}_{s,T}\bigvee{\mathscr{F}_{T,\infty}^{\hat{B}}}$
measurable for $t\leq s\leq T$ and $\psi(\cdot,\omega)$ is
continuous $P$-a.s.;
\item $E[\sup_{t\leq s\leq T}\|\psi(s)\|_\mathbb{S}^2]<\infty$.
\end{enumerate-roman}
\end{defi}

For a positive $K$, we consider the following infinite horizon BDSDE
with the infinite dimensional Brownian motion $\hat B$ as noise and
$Y_t$ taking values on a separable Hilbert space $H$, $Z_t$ taking
values on $\mathcal{L}^2_{\mathbb{R}^{d}}(H)$ (the space of all
Hilbert-Schmidt operators from ${\mathbb{R}^{d}}$ to $H$ with the
Hilbert-Schmidt norm):
\begin{eqnarray}\label{qi1}
{\rm e}^{-Kt}Y_{t}&=&\int_{t}^{\infty}{\rm e}^{-Kr}f(r,Y_{r},Z_{r})dr+\int_{t}^{\infty}K{\rm e}^{-Kr}Y_{r}dr\nonumber\\
&&-\int_{t}^{\infty}{\rm e}^{-Kr}
g(r,Y_{r},Z_{r})d^\dagger{\hat{B}}_r-\int_{t}^{\infty}{\rm
e}^{-Kr}Z_{r}dW_r,\ \ \ t\geq0.
\end{eqnarray}
Assume $f:[0,\infty)\times\Omega\times H\times
\mathcal{L}^2_{\mathbb{R}^{d}}(H){\longrightarrow H}$,
$g:[0,\infty)\times\Omega\times H\times
\mathcal{L}^2_{\mathbb{R}^{d}}(H){\longrightarrow\mathcal{L}^2_{U_0}(H)}$
are
$\mathscr{B}_{\mathbb{R}^{+}}\otimes\mathscr{F}\otimes\mathscr{B}_H\otimes\mathscr{B}_{\mathcal{L}^2_{\mathbb{R}^{d}}(H)}$
measurable such that for any $(t,Y,Z)\in[0,\infty)\times H\times
\mathcal{L}^2_{\mathbb{R}^{d}}(H)$, $f(t,Y,Z)$, $g(t,Y,Z)$ are
$\mathscr{F}_t$ measurable, where $U_0=Q^{1\over 2}(U)\subset U$ is
a separable Hilbert space with the norm $<u,v>_{U_0}=<Q^{-{1\over
2}}u, Q^{-{1\over 2}}v>_U$ and the complete orthonormal base
$\{\sqrt{\lambda_i}e_i\}_{i=1}^{\infty}$,
 $\mathcal{L}^2_{U_0}(H)$ is the space of
all Hilbert-Schmidt operators from $U_0$ to $H$ with the
Hilbert-Schmidt norm. It is noted that the $Q$-Wiener process $(\hat
B_t)_{t\geq 0}$ is a cylinderical Wiener process on $U_0$, and both
$\mathcal{L}^2_{U_0}(H)$ and $\mathcal{L}^2_{\mathbb{R}^{d}}(H)$ are
Hilbert spaces.

Note that the integral w.r.t. $\hat{B}$ is a "backward It$\hat {\rm
o}$'s integral" and the integral w.r.t. $W$ is a standard forward
It$\hat {\rm o}$'s integral. The forward integrals in Hilbert space
with respect to $Q$-Wiener processes were defined in Da Prato and
Zabczyk \cite{pr-za1}. To see the backward one, let
$\{h(s)\}_{s\geq0}$ be a stochastic process with values on
${\mathcal{L}^2_{U_0}(H)}$ such that $h(s)$ is $\mathscr{F}_s$
measurable for any $s\geq0$ and locally square integrable, i.e. for
any $0\leq a\leq b<\infty$,
$\int_{a}^{b}\|h(s)\|_{\mathcal{L}^2_{U_0}(H)}^2ds<\infty$ a.s..
Since ${\mathscr{F}}_s$ is a backward filtration with respect to
$\hat B$, so from the one-dimensional backward It$\hat {\rm o}$'s
integral and relation with forward integral, for $0\leq T\leq
T^{\prime}$, we have
\begin{eqnarray*}
\int _t^T\sqrt{\lambda_j}<h(s)e_j,f_k>d^\dagger\hat \beta_j(s)=-\int
_{T^{\prime}-T}^{T^{\prime}-t}\sqrt{\lambda_j}<h(T^{\prime}-s)e_j,f_k>
d \beta_j(s), \ \ j,k=1,2,\cdots
\end{eqnarray*}
where $\beta _j(s)=\hat \beta _j(T^{\prime}-s)-\hat \beta
_j(T^{\prime})$, $j=1,2,\cdots$, and so $B_s=\hat
B_{T^{\prime}-s}-\hat B_{T^{\prime}}$. Here $\{f_k\}$ is the
complete orthonormal basis in $H$. From approximation theorem of the
stochastic integral in Hilbert space (\cite{pr-za1}), we have
\begin{eqnarray*}
\int _{T^{\prime}-T}^{T^{\prime}-t}h(T^{\prime}-s)d
B_s=\sum_{j,k=1}^{\infty}\int
_{T^{\prime}-T}^{T^{\prime}-t}\sqrt{\lambda_j}<h(T^{\prime}-s)e_j,f_k>
d\beta_j(s)f_k.
\end{eqnarray*}
Similarly we also have
\begin{eqnarray*}
\int _t^Th(s)d^\dagger\hat B_s=\sum_{j,k=1}^{\infty} \int
_t^T\sqrt{\lambda_j}<h(s)e_j,f_k>d^\dagger\hat \beta_j(s)f_k.
\end{eqnarray*}
It turns out that
\begin{eqnarray}\label{bmrelation}
\int_{t}^{T}h(s)d^\dagger{\hat{B}}_s=-\int_{T'-T}^{T'-t}
h({T'-s})d{B}_s\ \ \ \ \rm{a.s.}.
\end{eqnarray}
Later we will consider another Hilbert space
${\mathcal{L}^p_{U_0}(H)}$ $(p>2)$, a subspace of
${\mathcal{L}^2_{U_0}(H)}$, including all
$h\in{\mathcal{L}^2_{U_0}(H)}$ which satisfy
$$\|h\|^p_{\mathcal{L}^p_{U_0}(H)}\triangleq\sum_{j,k=1}^{\infty}\lambda_j^{p\over2}|\langle he_j,f_k\rangle|^p<\infty.$$
\begin{defi}\label{qi021}
Let $H_0$ be a dense subset of $H$. If $(Y,Z)\in S^{2,-K}\bigcap
M^{2,-K}([0,\infty);H)\bigotimes
M^{2,-K}\\([0,\infty);\mathcal{L}^2_{\mathbb{R}^{d}}(H))$,
and for any $\varphi\in H_0$,
\begin{eqnarray}\label{qi3}
&&\langle{\rm e}^{-Kt}Y_t,\varphi\rangle=\langle\int_{t}^{\infty}{\rm e}^{-Kr}f(r,Y_r,Z_r)dr,\varphi\rangle+\langle\int_{t}^{\infty}K{\rm e}^{-Kr}Y_{r}dr,\varphi\rangle\nonumber\\
&&\ \ \ \ \ \ \ \ \ \ \ \ \ \ -\langle\int_{t}^{\infty}{\rm
e}^{-Kr}g(r,Y_r,Z_r)d^\dagger{\hat{B}}_r,\varphi\rangle-\langle\int_{t}^{\infty}{\rm
e}^{-Kr}Z_r dW_r,\varphi\rangle,\ t\geq0\ P-{\rm a.s.},
\end{eqnarray}
or equivalently
\begin{eqnarray}
\left\{\begin{array}{l}\label{qi4}
\langle Y_t,\varphi\rangle=\langle Y_T,\varphi\rangle+\langle\int_{t}^{T}f(r,Y_r,Z_r)dr,\varphi\rangle-\langle\int_{t}^{T}g(r,Y_r,Z_r)d^\dagger{\hat{B}}_r,\varphi\rangle-\langle\int_{t}^{T}Z_rdW_r,\varphi\rangle\\
\lim_{T\rightarrow\infty}\langle{\rm e}^{-{K}T}Y_T,\varphi\rangle=0\
\ \ \rm{a.s.},
\end{array}\right.
\end{eqnarray}
then we call $(Y,Z)$ a solution of Eq.(\ref{qi1}) in $H$.
\end{defi}
\begin{rmk}\label{zhang315}
\begin{enumerate-roman}
\item Applying It$\hat {\rm o}$'s formula in $H$ (see
\cite{pr-za1}), we have the equivalent form of Eq.(\ref{qi1})
\begin{eqnarray}
\left\{\begin{array}{l}\label{qi2}
Y_t=Y_T+\int_{t}^{T}f(r,Y_r,Z_r)dr-\int_{t}^{T}g(r,Y_r,Z_r)d^\dagger{\hat{B}}_r-\int_{t}^{T}Z_rdW_r\\
\lim_{T\rightarrow\infty}{\rm e}^{-{K}T} Y_T=0\ \ \ \rm{a.s.};
\end{array}\right.
\end{eqnarray}
\item one can easily verify that the
above definition doesn't depend on the choice of $H_0$ due to the
continuity of the inner product;
\item the uniqueness of $Y$ in $S^{2,-K}([0,\infty);H)$ implies if
$(Y^\prime,Z^\prime)$ is another solution, then $Y_s=Y^\prime_s$ for
all $s\geq0$ a.s.. The uniqueness of $Z$ implies $Z_s=Z^\prime_s$
for a.e. $s\in[0,\infty)$ a.s.. But we can modify the $Z$ at the
measure zero exceptional set of $s$ such that $Z_s=Z^\prime_s$ for
all $s\geq0$ a.s..
\end{enumerate-roman}
\end{rmk}

The first main purpose of this section is to study the stationary
property of the solution of BDSDE (\ref{qi1}) on $H$ if the solution
exists and is unique. In order to show the main idea, we first
assume that there exists a unique solution of Eq.(\ref{qi1}). The
study of the existence and uniqueness of Eq.(\ref{qi1}) will be
deferred to later sections (sections 3-5).

We now construct the measurable metric dynamical system through
defining a measurable and measure-preserving shift. Let
$\hat{\theta}_t:\Omega\longrightarrow\Omega$, $t\geq0$, be a
measurable mapping on $(\Omega, {\mathscr{F}}, P)$, defined by
$\hat{\theta}_{t}\circ \hat{B}_s=\hat{B}_{s+t}-\hat{B}_t$, \
$\hat{\theta}_{t}\circ W_s=W_{s+t}-W_t$. Then for any $s,t\geq0$,
\begin{description}
\item[$(\textrm{i})$]$P\cdot\hat{\theta}_{t}^{-1}=P$;
\item[$(\textrm{i}\textrm{i})$]$\hat{\theta}_{0}=I$, where $I$ is the identity transformation on $\Omega$;
\item[$(\textrm{i}\textrm{i}\textrm{i})$]$\hat{\theta}_{s}\circ\hat{\theta}_{t}=\hat{\theta}_{s+t}$.
\end{description}

Also for an arbitrary $\mathscr{F}$ measurable
$\phi:\Omega\longrightarrow H$, set
\begin{eqnarray*}
\hat{\theta}\circ\phi(\omega)=\phi\big(\hat{\theta}(\omega)\big).
\end{eqnarray*}

We give the following bounded and stationary conditions for $f$, $g$
w.r.t. $\hat{\theta}_{\cdot}$:
\begin{description}
\item[(A.1).] There exist a constant $M_1\geq0$, and functions $\tilde{f}(\cdot)\in M^{2,-K}([0,\infty);R^+)$, $\tilde{g}(\cdot)\in
M^{2,-K}([0,\infty);R^+)$ s.t. for any $s\geq0$, $Y\in H$ and $Z\in
\mathcal{L}^2_{\mathbb{R}^{d}}(H)$,
$$\|f(s,Y,Z)\|_H^2\leq\tilde{f}^2(s)+M_1\|Y\|^2_H+M_1\|Z\|^2_{\mathcal{L}^2_{\mathbb{R}^{d}}(H)},$$$$\|g(s,Y,Z)\|^2_{\mathcal{L}^2_{U_0}(H)}\leq\tilde{g}^2(s)+M_1\|Y\|^2_H+M_1\|Z\|^2_{\mathcal{L}^2_{\mathbb{R}^{d}}(H)};$$
\item[(A.2).] For any $r,s\geq0$, $Y\in H$ and $Z\in \mathcal{L}^2_{\mathbb{R}^{d}}(H)$, $\hat{\theta}_r\circ f(s,Y,Z)=f(s+r,Y,Z)$, $\hat{\theta}_r\circ g(s,Y,Z)=g(s+r,Y,Z)$.
\end{description}

We start from the following general result about the stationarity of
the solution of infinite horizon BDSDE.
\begin{prop}\label{qi031}
Assume Eq.(\ref{qi1}) has a unique solution $(Y,Z)$, then under
Conditions (A.1) and (A.2), $(Y_t, Z_t)_{t\geq0}$ is a "perfect"
stationary solution, i.e.
\begin{eqnarray*}
\hat{\theta}_r\circ Y_t=Y_{t+r}, \ \ \hat{\theta}_r\circ
Z_t=Z_{t+r}\ \ {\rm for}\ {\rm all}\ r,\ t\geq0\ {\rm a.s.}.
\end{eqnarray*}
\end{prop}
{\em Proof}. Let ${B}_s=\hat{B}_{T'-s}-\hat{B}_{T'}$ for arbitrary
$T'>0$ and $-\infty<s\leq T'$. Then ${B}_s$ is a Brownian motion
with ${B}_0=0$. For any $r\geq0$, applying $\hat{\theta}_r$ on
${B}_s$, we have
\begin{eqnarray*}
\hat{\theta}_r\circ{B}_s&=&\hat{\theta}_r\circ(\hat{B}_{T'-s}-\hat{B}_{T'})=\hat{B}_{T'-s+r}-\hat{B}_{T'+r}\\
&=&(\hat{B}_{T'-s+r}-\hat{B}_{T'})-(\hat{B}_{T'+r}-\hat{B}_T')={B}_{s-r}-{B}_{-r}.
\end{eqnarray*}
So for $0\leq t\leq T\leq T'$ and a locally square integrable
process $\{h(s)\}_{s\geq0}$, by (\ref{bmrelation})
\begin{eqnarray*}
\hat{\theta}_r\circ\int_{t}^{T}h(s)d^\dagger{\hat{B}}_s&=&-\hat{\theta}_r\circ\int_{T'-T}^{T'-t}h(T'-s)d{B}_s\nonumber\\
&=&-\int_{T'-T}^{T'-t}\hat{\theta}_r\circ h(T'-s)d{B}_{s-r}\nonumber\\
&=&-\int_{T'-T-r}^{T'-t-r}\hat{\theta}_r\circ h(T'-s-r)d{B}_{s}\nonumber\\
&=&\int_{t+r}^{T+r}\hat{\theta}_r\circ h(s-r)d^\dagger{\hat{B}}_s.
\end{eqnarray*}
As $T'$ can be chosen arbitrarily, so we can get for arbitrary
$T\geq0$, $0\leq t\leq T$, $r\geq0$,
\begin{eqnarray}\label{qi5}
\hat{\theta}_r\circ\int_{t}^{T}h(s)d^\dagger
\hat{B}_s=\int_{t+r}^{T+r}\hat{\theta}_r\circ
h(s-r)d^\dagger{\hat{B}}_s.
\end{eqnarray}
It is easy to see that $g(\cdot,Y_{\cdot},Z_{\cdot})$ is locally
square integrable from Condition (A.1), hence by Condition (A.2) and
(\ref{qi5})
\begin{eqnarray}\label{qi6}
\hat{\theta}_r\circ\int_{t}^{T}g(s,Y_s,Z_s)d^\dagger
\hat{B}_s=\int_{t+r}^{T+r}g(s,\hat{\theta}_r\circ
Y_{s-r},\hat{\theta}_r\circ Z_{s-r})d^\dagger{\hat{B}}_s.
\end{eqnarray}

We consider the equivalent form Eq.(\ref{qi2}) instead of
Eq.(\ref{qi1}). Applying the operator $\hat{\theta}_r$ on both sides
of Eq.(\ref{qi2}) and by (\ref{qi6}), we know that
$\hat{\theta}_r\circ Y_t$ satisfies the following equation
\begin{eqnarray}
\left\{\begin{array}{l}\label{qi7}
\hat{\theta}_r\circ Y_t=\hat{\theta}_r\circ Y_T+\int_{t+r}^{T+r}f(s,\hat{\theta}_r\circ Y_{s-r}, \hat{\theta}_r\circ Z_{s-r})ds\\
\ \ \ \ \ \ \ \ \ \ \ \ -\int_{t+r}^{T+r}g(s,\hat{\theta}_r\circ
Y_{s-r}, \hat{\theta}_r\circ Z_{s-r})d^\dagger{\hat{B}}_s
 -\int_{t+r}^{T+r}\hat{\theta}_r\circ Z_{s-r}dW_s\\
\lim_{T\rightarrow\infty}{\rm e}^{-K(T+r)}(\hat{\theta}_r\circ
Y_T)=0\ \ \ \rm{a.s.}.
\end{array}\right.
\end{eqnarray}
On the other hand, from Eq.(\ref{qi2}), it follows that
\begin{eqnarray}\label{qi8}
\left\{\begin{array}{l} Y_{t+r}=Y_{T+r}+\int_{t+r}^{T+r}f(s,Y_s, Z_s)ds-\int_{t+r}^{T+r}g(s,Y_s, Z_s)d^\dagger{\hat{B}}_s -\int_{t+r}^{T+r}Z_sdW_s\\
\lim_{T\rightarrow\infty}{\rm e}^{-K(T+r)}Y_{T+r}=0 \ \ \ \rm{a.s.}.
\end{array}\right.
\end{eqnarray}
Let $\hat{Y}_{\cdot}=\hat{\theta}_r\circ Y_{\cdot-r}$,
$\hat{Z}_{\cdot}=\hat{\theta}_r\circ Z_{\cdot-r}$. By the uniqueness
of solution of Eq.(\ref{qi2}) and Remark \ref{zhang315} (iii), it
follows from comparing (\ref{qi7}) with (\ref{qi8}) that for any
$r\geq0$,
\begin{eqnarray*}
\hat{\theta}_r\circ Y_t=\hat{Y}_{t+r}=Y_{t+r}, \ \
\hat{\theta}_r\circ Z_t=\hat{Z}_{t+r}=Z_{t+r}\ \ {\rm for}\ {\rm
all}\ t\geq0\ \rm{a.s.}.
\end{eqnarray*}
Then by perfection procedure (\cite{ar}, \cite{ar-sc}), we can prove
above identities are true for all $t$, $r\geq0$ a.s.. We proved the desired result. $\hfill\diamond$\\

An important application of the BDSDEs is to connect its solution
with the solution of the corresponding SPDEs. If some kind of
relationship is established, we can transfer stationary solutions
from the infinite horizon BDSDEs to SPDEs. In this way, we are in
access to stationary solutions of the SPDEs due to the stationary
property of solutions of infinite horizon BDSDEs. For this, a specific
Hilbert space  $H=L_{\rho}^2({\mathbb{R}^{d}};{\mathbb{R}^{1}})$
defined below is considered.
The main aim of rest of this section is to construct the
stationary solution of the SPDEs. Some proofs are given in this
sections. But many detailed proofs are postponed to later sections.

In the following we consider the case
$H=L_{\rho}^2({\mathbb{R}^{d}};{\mathbb{R}^{1}})$ with the inner
product $\langle
u_1,u_2\rangle=\int_{\mathbb{R}^d}u_1(x)u_2(x)\rho^{-1}(x)dx$, a
$\rho$-weighted $L^2$ space. Here $\rho(x)=(1+|x|)^q$, $q>3$, is a
weight function. It is easy to see that
$\rho(x):\mathbb{R}^d\longrightarrow\mathbb{R}^1$ is a continuous
positive function satisfying
$\int_{\mathbb{R}^{d}}|x|^p\rho^{-1}(x)dx<\infty$ for any
$p\in(2,q-1)$. Note that we can consider more general $\rho$ which
satisfies the above condition and conditions in \cite{ba-ma} and all the results
of this paper still hold.

We can write down the solution spaces following Definition
\ref{qi00}:
$M^{2,-K}([0,\infty);L_{\rho}^2({\mathbb{R}^{d}};{\mathbb{R}^{1}}))$,
$M^{2,-K}([0,\infty);L_{\rho}^2({\mathbb{R}^{d}};{\mathbb{R}^{d}}))$
and
$S^{2,-K}([0,\infty);L_{\rho}^2({\mathbb{R}^{d}};{\mathbb{R}^{1}}))$.
Similar to the definition for
$M^{2,-K}\\([0,\infty);L_{\rho}^2({\mathbb{R}^{d}};{\mathbb{R}^{d}}))$,
we can also define
$M^{p,-K}([0,\infty);L_{\rho}^p({\mathbb{R}^{d}};{\mathbb{R}^{d}}))$.

For $k\geq0$, we denote by $C_{l,b}^k(\mathbb{R}^p,\mathbb{R}^q)$
the set of $C^k$-functions whose partial derivatives of order less
than or equal to $k$ are bounded and by
$H^k_\rho(\mathbb{R}^d;\mathbb{R}^{1})$ the $\rho$-weighted Sobolev
space (See e.g. \cite{ba-ma}). In order to connect BDSDEs with
SPDEs, the form of BDSDEs should be a kind of FBDSDEs (forward and
backward doubly SDEs). So we first give the following forward SDE.

For $s\geq t$, let $X_{s}^{t,x}$ be a diffusion process given by the
solution of
\begin{eqnarray}\label{qi17}
X_{s}^{t,x}=x+\int_{t}^{s}b(X_{u}^{t,x})du+\int_{t}^{s}\sigma(X_{u}^{t,x})dW_u,
\end{eqnarray}
where $b\in C_{l,b}^2(\mathbb{R}^{d};\mathbb{R}^{d})$, $\sigma\in
C_{l,b}^3(\mathbb{R}^{d};\mathbb{R}^{d}\times\mathbb{R}^{d})$, and
for $0\leq s<t$, we regulate $X_{s}^{t,x}=x$.

For any $r\geq0$, $s\geq t$, $x\in\mathbb{R}^d$, apply $\theta_r$ on
SDE (\ref{qi17}), then
\begin{eqnarray*}
\hat{\theta}_r\circ
X_{s}^{t,x}=x+\int_{t+r}^{s+r}b(\hat{\theta}_r\circ
X_{u-r}^{t,x})du+\int_{t+r}^{s+r}\sigma(\hat{\theta}_r\circ
X_{u-r}^{t,x})dW_u.
\end{eqnarray*}
So by the uniqueness of the solution and a perfection procedure
(c.f. \cite{ar}), we have
\begin{eqnarray}\label{qi18}
\hat{\theta}_r\circ X_{s}^{t,x}=X_{s+r}^{t+r,x},\ \ {\rm for}\ {\rm
all}\ r,s,t,x\ \ \ {\rm a.s.}.
\end{eqnarray}
Moreover, it is well-known that the solution defines a stochastic
flow of diffeomorphism
$X_s^{t,\cdot}:\mathbb{R}^d\rightarrow\mathbb{R}^d$ and denote by
$\hat{X}_s^{t,\cdot}$ the inverse flow (See e.g. Kunita \cite{ku2}).
Denote by $J(\hat{X}_s^{t,x})$ the determinant of the Jacobi matrix
of $\hat{X}_s^{t,x}$. For $\varphi\in
H^k_\rho(\mathbb{R}^d;\mathbb{R}^{1})$, we define a process
$\varphi_t:\Omega\times[0,T]\times\mathbb{R}^d\rightarrow\mathbb{R}^1$
by $\varphi_t(s,x)=\varphi(\hat{X}_s^{t,x})J(\hat{X}_s^{t,x})$. It
is proved in \cite{ba-ma} that $\varphi_t(s,\cdot)\in
H^k_\rho(\mathbb{R}^d;\mathbb{R}^{1})$ and for $u\in
{H^k_\rho}^*(\mathbb{R}^d;\mathbb{R}^{1})$,
$\int_{\mathbb{R}^d}u(x)\varphi(x)dx\triangleq\sum_{0\leq|\alpha|\leq
k}\int_{\mathbb{R}^d}u_\alpha(x)D^\alpha\varphi(x)dx\leq\sum_{0\leq|\alpha|\leq
k}\sqrt{{\int_{\mathbb{R}^d}|u_\alpha(x)|^2\rho^{-1}(x)dx}{\int_{\mathbb{R}^d}|D^\alpha\varphi(x)|^2\rho(x)dx}}<\infty$
and
$\int_{\mathbb{R}^d}u(y)\varphi_t(s,y)dy=\int_{\mathbb{R}^d}u({X}_s^{t,x})\\
\cdot\varphi(x)dx$.

The following lemma plays an important role in the analysis in this
article. It is an extension of equivalence of norm principle given
in \cite {ku1}, \cite{ba-le}, \cite{ba-ma} to the cases when
$\varphi$ and $\Psi$ are random.
\begin{lem}\label{qi045}(generalized equivalence
of norm principle) Let $\rho$ be the weight function defined at the
beginning of this section and $X$ be a diffusion process defined
above. If $s\in[t,T]$,
$\varphi:\Omega\times\mathbb{R}^d\rightarrow\mathbb{R}^1$ is
independent of $\mathscr{F}^W_{t,s}$ and $\varphi\rho^{-1}\in
L^1(\Omega\otimes\mathbb{R}^{d})$, then there exist two constants
$c>0$ and $C>0$ such that
\begin{eqnarray*}
cE[\int_{\mathbb{R}^{d}}|\varphi(x)|\rho^{-1}(x)dx]\leq
E[\int_{\mathbb{R}^{d}}|\varphi(X_{s}^{t,x})|\rho^{-1}(x)dx]\leq
CE[\int_{\mathbb{R}^{d}}|\varphi(x)|\rho^{-1}(x)dx].
\end{eqnarray*}
Moreover if
$\Psi:\Omega\times[t,T]\times\mathbb{R}^d\rightarrow\mathbb{R}^1$,
$\Psi(s,\cdot)$ is independent of $\mathscr{F}^W_{t,s}$ and
$\Psi\rho^{-1}\in L^1(\Omega\otimes[t,T]\otimes\mathbb{R}^{d})$,
then
\begin{eqnarray*}
&&cE[\int_{t}^{T}\int_{\mathbb{R}^{d}}|\Psi(s,x)|\rho^{-1}(x)dxds]\leq E[\int_{t}^{T}\int_{\mathbb{R}^{d}}|\Psi(s,X_{s}^{t,x})|\rho^{-1}(x)dxds]\\
&\leq&CE[\int_{t}^{T}\int_{\mathbb{R}^{d}}|\Psi(s,x)|\rho^{-1}(x)dxds].
\end{eqnarray*}
\end{lem}
{\em Proof}. Using the conditional expectation w.r.t.
$\mathscr{F}^W_{t,s}$ and noting that
${\rho^{-1}(\hat{X}_{s}^{t,y})J(\hat{X}_s^{t,y})\over{\rho^{-1}(y)}}$
is $\mathscr{F}^W_{t,s}$ measurable and $|\varphi(y)|\rho^{-1}(y)$
is independent of $\mathscr{F}^W_{t,s}$, we have
\begin{eqnarray*}
&&E[\int_{\mathbb{R}^{d}}|\varphi(X_{s}^{t,x})|\rho^{-1}(x)dx]\\
&=&\int_{\mathbb{R}^{d}}E[\ E[\ |\varphi(y)|\rho^{-1}(y){\rho^{-1}(\hat{X}_{s}^{t,y})J(\hat{X}_s^{t,y})\over{\rho^{-1}(y)}}\ |\mathscr{F}^W_{t,s}]\ ]dy\\
&=&\int_{\mathbb{R}^{d}}E[\
|\varphi(y)|\rho^{-1}(y)]E[{\rho^{-1}(\hat{X}_{s}^{t,y})J(\hat{X}_s^{t,y})\over{\rho^{-1}(y)}}]dy.
\end{eqnarray*}
By Lemma 5.1 in \cite{ba-ma}, $c\leq
E[{\rho^{-1}(\hat{X}_{s}^{t,y})J(\hat{X}_s^{t,y})\over{\rho^{-1}(y)}}]\leq
C$ for any $y\in\mathbb{R}^{d}$, $s\in[t,T]$, the claim follows.
$\hfill\diamond$
\\

By Lemma \ref{qi045}, it is easy to deduce that
$X_{\cdot}^{t,\cdot}\in
M^{p,-K}([0,\infty);L_{\rho}^p({\mathbb{R}^{d}};{\mathbb{R}^{d}}))$
for $K\in\mathbb{R}^{+}$.

Now we consider the following BDSDE with infinite dimensional noise
on infinite horizon
\begin{eqnarray}\label{qi13}
{\rm e}^{-Ks}Y_{s}^{t,x}&=&\int_{s}^{\infty}{\rm e}^{-Kr}f(X_{r}^{t,x},Y_{r}^{t,x},Z_{r}^{t,x})dr+\int_{s}^{\infty}K{\rm e}^{-Kr}Y_{r}^{t,x}dr\nonumber\\
&&-\int_{s}^{\infty}{\rm e}^{-Kr}
g(X_{r}^{t,x},Y_{r}^{t,x},Z_{r}^{t,x})d^\dagger{\hat{B}}_r-\int_{s}^{\infty}{\rm
e}^{-Kr}\langle Z_{r}^{t,x},dW_r\rangle.
\end{eqnarray}
Here
$\hat{B}_r=\sum_{j=1}^{\infty}\sqrt{\lambda_j}\hat{\beta}_j(r)e_j$,
$\{\hat{\beta}_j(r)\}_{j=1,2,\cdots}$ are mutually independent
one-dimensional Brownian motions. Note that we will solve
Eq.(\ref{qi13}) for $Y_{r}^{t,\cdot}\in
L_{\rho}^2({\mathbb{R}^{d}};{\mathbb{R}^{1}})$ and
$Z_{r}^{t,\cdot}\in
\mathcal{L}^2_{\mathbb{R}^{d}}(L_{\rho}^2({\mathbb{R}^{d}};{\mathbb{R}^{1}}))=L_{\rho}^2({\mathbb{R}^{d}};{\mathbb{R}^{d}})$.

Set $g_j\triangleq
g\sqrt{\lambda_j}e_j:\mathbb{R}^{d}\times\mathbb{R}^1\times\mathbb{R}^{d}{\longrightarrow{\mathbb{R}^1}}$,
then Eq.(\ref{qi13}) is equivalent to
\begin{eqnarray*}\label{qi14}
{\rm e}^{-Ks}Y_{s}^{t,x}&=&\int_{s}^{\infty}{\rm e}^{-Kr}f(X_{r}^{t,x},Y_{r}^{t,x},Z_{r}^{t,x})dr+\int_{s}^{\infty}K{\rm e}^{-Kr}Y_{r}^{t,x}dr\nonumber\\
&&-\sum_{j=1}^{\infty}\int_{s}^{\infty}{\rm
e}^{-Kr}g_j(X_{r}^{t,x},Y_{r}^{t,x},Z_{r}^{t,x})d^\dagger{\hat{\beta}}_j(r)-\int_{s}^{\infty}{\rm
e}^{-Kr}\langle Z_{r}^{t,x},dW_r\rangle.
\end{eqnarray*}
Referring to Definition \ref{qi021} and noting that
$C_c^{0}(\mathbb{R}^d;\mathbb{R}^1)$ is dense in
$L_{\rho}^2({\mathbb{R}^{d}};{\mathbb{R}^{1}})$ under the norm
$(\int_{\mathbb{R}^d}|\cdot|^2\rho^{-1}(x)dx)^{1\over2}$, we can
define the solution in
$L_{\rho}^2({\mathbb{R}^{d}};{\mathbb{R}^{1}})$ as follows:
\begin{defi}\label{qi041}
A pair of processes $(Y_{\cdot}^{t,\cdot},Z_{\cdot}^{t,\cdot})\in
S^{2,-K}\bigcap
M^{2,-K}([0,\infty);L_{\rho}^2({\mathbb{R}^{d}};{\mathbb{R}^{1}}))\bigotimes
M^{2,-K}\\([0,\infty);L_{\rho}^2({\mathbb{R}^{d}};{\mathbb{R}^{d}}))$
is called a solution of Eq.(\ref{qi13}) if for an arbitrary
$\varphi\in C_c^{0}(\mathbb{R}^d;\mathbb{R}^1)$,
\begin{eqnarray}\label{qi15}
\int_{\mathbb{R}^{d}}{\rm
e}^{-Ks}Y_s^{t,x}\varphi(x)dx&=&\int_{s}^{\infty}\int_{\mathbb{R}^{d}}{\rm e}^{-Kr}f(X_{r}^{t,x},Y_r^{t,x},Z_r^{t,x})\varphi(x)dxdr+\int_{s}^{\infty}\int_{\mathbb{R}^{d}}K{\rm e}^{-Kr}Y_{r}^{t,x}\varphi(x)dxdr\nonumber\\
&&-\sum_{j=1}^{\infty}\int_{s}^{\infty}\int_{\mathbb{R}^{d}}{\rm e}^{-Kr}g_j(X_{r}^{t,x},Y_r^{t,x},Z_r^{t,x})\varphi(x)dxd^\dagger{\hat{\beta}}_j(r)\nonumber\\
&&-\int_{s}^{\infty}\langle \int_{\mathbb{R}^{d}}{\rm
e}^{-Kr}Z_r^{t,x}\varphi(x)dx,dW_r\rangle\ \ \ P-{\rm a.s.}.
\end{eqnarray}
\end{defi}
Note that in (\ref{qi15}) we leave out the weight function $\rho$ in
the inner product due to the arbitrariness of $\varphi$.

If Eq.(\ref{qi13}) has a unique solution, then for an arbitrary $T$,
$Y_{T}^{t,x}$ satisfies
\begin{eqnarray}\label{qia}
Y_{s}^{t,x}&=&Y_{T}^{t,x}+\int_{s}^{T}f(X_{r}^{t,x},Y_{r}^{t,x},Z_{r}^{t,x})dr-\int_{s}^{T}g(X_{r}^{t,x},Y_{r}^{t,x},Z_{r}^{t,x})d^\dagger{\hat{B}}_r-\int_{s}^{T}\langle
Z_{r}^{t,x},dW_r\rangle.\nonumber\\
\end{eqnarray}
In section 4, we will deduce the following SPDE associated with
BDSDE (\ref{qia})
\begin{eqnarray}\label{zhang685}
u(t,x)&=&u(T,x)+\int_{t}^{T}[\mathscr{L}u(s,x)+f\big(x,u(s,x),(\sigma^*\nabla u)(s,x)\big)]ds\nonumber\\
&&-\int_{t}^{T}g\big(x,u(s,x),(\sigma^*\nabla u)(s,x)\big)d^\dagger
\hat{B}_s.
\end{eqnarray}
Here $\mathscr{L}$ is given by (\ref{buchong1}), $u(T,x)=Y_T^{T,x}$.
But we can normally study general $u(T,x)$ unless we consider the
stationary solution.

Now following Definition \ref{zhao005} we write down the solution
spaces needed in our paper:
$M^{2,0}([t,T];L_{\rho}^2({\mathbb{R}^{d}};{\mathbb{R}^{1}}))$,
$M^{2,0}([t,T];L_{\rho}^2({\mathbb{R}^{d}};{\mathbb{R}^{d}}))$ and
$S^{2,0}([t,T];L_{\rho}^2({\mathbb{R}^{d}};{\mathbb{R}^{1}}))$.
\begin{defi}\label{qi042}
A process $u$ is called a weak solution {\rm(}solution in
$L_{\rho}^2({\mathbb{R}^{d}};{\mathbb{R}^{1}})${\rm)} of
Eq.(\ref{zhang685}) if $(u,\sigma^*\nabla u)\in
M^{2,0}([0,T];L_{\rho}^2({\mathbb{R}^{d}};{\mathbb{R}^{1}}))\bigotimes
M^{2,0}([0,T];L_{\rho}^2({\mathbb{R}^{d}};{\mathbb{R}^{d}}))$ and
for an arbitrary $\Psi\in
C_c^{1,\infty}([0,T]\times\mathbb{R}^d;\mathbb{R}^1)$,
\begin{eqnarray}\label{qi16}
&&\int_{t}^{T}\int_{\mathbb{R}^{d}}u(s,x)\partial_s\Psi(s,x)dxds+\int_{\mathbb{R}^{d}}u(t,x)\Psi(t,x)dx-\int_{\mathbb{R}^{d}}u(T,x)\Psi(T,x)dx\nonumber\\
&&-{1\over2}\int_{t}^{T}\int_{\mathbb{R}^{d}}(\sigma^*\nabla u)(s,x)(\sigma^*\nabla\Psi)(s,x)dxds-\int_{t}^{T}\int_{\mathbb{R}^{d}}u(s,x)div\big((b-\tilde{A})\Psi\big)(s,x)dxds\nonumber\\
&=&\int_{t}^{T}\int_{\mathbb{R}^{d}}f\big(x,u(s,x),(\sigma^*\nabla u)(s,x)\big)\Psi(s,x)dxds\nonumber\\
&&-\sum_{j=1}^{\infty}\int_{t}^{T}\int_{\mathbb{R}^{d}}g_j\big(x,u(s,x),(\sigma^*\nabla
u)(s,x)\big)\Psi(s,x)dxd^\dagger{\hat{\beta}}_j(s)\ \ \ P-{\rm
a.s.}.
\end{eqnarray}
Here $\tilde{A}_j\triangleq{1\over2}\sum_{i=1}^d{\partial
a_{ij}(x)\over\partial x_i}$, and
$\tilde{A}=(\tilde{A}_1,\tilde{A}_2,\cdot\cdot\cdot,\tilde{A}_d)^*$.
\end{defi}

This definition can be easily understood if we note the following
integration by parts formula: for $\varphi_1,\varphi_2\in
C^2(\mathbb{R}^{d})$,
\begin{eqnarray*}
-\int_{\mathbb{R}^{d}}\mathscr{L}\varphi_1(x)\varphi_2(x)dx={1\over2}\int_{\mathbb{R}^{d}}(\sigma^*\nabla
\varphi_1)(x)(\sigma^*\nabla\varphi_2)(x)dx+\int_{\mathbb{R}^{d}}\varphi_1(x)div\big((b-\tilde{A})\varphi_2\big)(x)dx.
\end{eqnarray*}

The main purpose of this section is to find the stationary solution
of SPDE (\ref{zhao21}) via the solution of BDSDE (\ref{qi13}). We
consider the following conditions:
\begin{description}
\item[{\rm(A.1)$'$.}] Functions $f:\mathbb{R}^{d}\times\mathbb{R}^1\times\mathbb{R}^{d}{\longrightarrow{\mathbb{R}^1}}$
and
$g:\mathbb{R}^{d}\times\mathbb{R}^1\times\mathbb{R}^{d}\longrightarrow
{\mathcal{L}^2_{U_0}(\mathbb{R}^{1})}$ are
$\mathscr{B}_{\mathbb{R}^{d}}\otimes\mathscr{B}_{\mathbb{R}^{1}}\otimes\mathscr{B}_{\mathbb{R}^{d}}$
measurable, and there exist constants
$M_2,M_{2j},C,C_j,\alpha_j\geq0$ with $\sum_{j=1}^\infty
M_{2j}<\infty$, $\sum_{j=1}^\infty C_j<\infty$ and
$\sum_{j=1}^\infty\alpha_j<{1\over2}$ s.t. for any $Y_1, Y_2\in
L_{\rho}^2({\mathbb{R}^{d}};{\mathbb{R}^{1}})$, $X_1, X_2, Z_1,
Z_2\in L_{\rho}^2({\mathbb{R}^{d}};{\mathbb{R}^{d}})$,
measurable $U:{\mathbb{R}^{d}}\rightarrow [0,1]$,
\begin{eqnarray*}
&&\int_{\mathbb{R}^d}U(x)|f(X_1(x),Y_1(x),Z_1(x))-f(X_2(x),Y_2(x),Z_2(x))|^2\rho^{-1}(x)dx\\
&\leq&\int_{\mathbb{R}^d}U(x)\big(M_2|X_1(x)-X_2(x)|^2+C|Y_1(x)-Y_2(x)|^2+C|Z_1(x)-Z_2(x)|^2\big)\rho^{-1}(x)dx,\\
&&\int_{\mathbb{R}^d}U(x)|g_j(X_1(x),Y_1(x),Z_1(x))-g_j(X_2(x),Y_2(x),Z_2(x))|^2\rho^{-1}(x)dx\\
&\leq&\int_{\mathbb{R}^d}U(x)\big(M_{2j}|X_1(x)-X_2(x)|^2+C_j|Y_1(x)-Y_2(x)|^2+\alpha_j|Z_1(x)-Z_2(x)|^2\big)\rho^{-1}(x)dx;
\end{eqnarray*}
\item[{\rm(A.2)$'$.}] For $p\in(2,q-1)$, $$\int_{\mathbb{R}^d}|f(x,0,0)|^p\rho^{-1}(x)dx<\infty\ {\rm
and}\
\int_{\mathbb{R}^d}\|g(x,0,0)\|_{\mathcal{L}^p_{U_0}(\mathbb{R}^{1})}^p\rho^{-1}(x)dx<\infty;$$
\item[{\rm(A.3)$'$.}] $b\in C_{l,b}^2(\mathbb{R}^{d};\mathbb{R}^{1})$, $\sigma\in
C_{l,b}^3(\mathbb{R}^{d}\times\mathbb{R}^{d};\mathbb{R}^{1})$.
Furthermore, for $p$ is given in {\rm(A.2)$'$}, if $L$ is the global
Lipschitz constant for $b$ and $\sigma$, $L$ satisfies
$K-pL-{p(p-1)\over2}L^2>0$;
\item[{\rm(A.4)$'$.}] There exists a
constant $\mu>0$ with
$2\mu-{pK}-pC-{{p(p-1)}\over2}\sum_{j=1}^{\infty}{C_j}>0$ s.t. for
any $Y_1, Y_2\in L_{\rho}^2({\mathbb{R}^{d}};{\mathbb{R}^{1}})$,
$X,Z\in L_{\rho}^2({\mathbb{R}^{d}};{\mathbb{R}^{d}})$,
measurable $U:{\mathbb{R}^{d}}\rightarrow [0,1]$,
\begin{eqnarray*}
&&\int_{\mathbb{R}^d}U(x)\big(Y_1(x)-Y_2(x)\big)\big(f(X(x),Y_1(x),Z(x))-f(X(x),Y_2(x),Z(x))\big)\rho^{-1}(x)dx\\
&\leq&-\mu\int_{\mathbb{R}^d}U(x){|Y_1(x)-Y_2(x)|}^2\rho^{-1}(x)dx.
\end{eqnarray*}
\end{description}

{\begin{rmk} We need monotone condition {\rm(A.4)}$'$ in
order to solve the infinite horizon BDSDEs. But it does not seem
obvious to replace the Lipschitz condition for $f$ in (A.1)$'$ by a weaker
condition on $f$ such as $f$ is continuous in $y$ using the
infinite horizon BSDE procedure (e.g. \cite{pa}). The difficulty
is due to the fact that we consider various conditions in the space
$L_{\rho}^2({\mathbb{R}^{d}};{\mathbb{R}^{1}})$ here rather than
pointwise ones, therefore we cannot solve the BDSDEs pointwise in $x$. However, our conjecture is that
the Lipschitz condition can be relaxed if we strengthen
the monotone condition in $L_{\rho}^2({\mathbb{R}^{d}};{\mathbb{R}^{1}})$
to a pointwise one. We will study this generality in future publications. Here due to the length of the paper, we
only consider the Lipschitz continuous function $f$ to initiate this
intrinsic method to the study of this basic problem.
\end{rmk}

We first acknowledge the two theorems below and give their proofs in
section 6.
\begin{thm}\label{qi043} Under Conditions {\rm(A.1)$'$}--{\rm(A.4)$'$}, Eq.(\ref{qi13}) has a unique solution
$(Y_{s}^{t,x},Z_{s}^{t,x})$. Moreover
$E[\sup_{s\geq0}\int_{\mathbb{R}^{d}}{\rm
e}^{{-{pK}}s}{{|{Y}_s^{t,x}|}^p}\rho^{-1}(x)dx]<\infty$.
\end{thm}
\begin{thm}\label{qi044} Under Conditions {\rm(A.1)$'$}--{\rm(A.4)$'$}, let $u(t,\cdot)\triangleq Y_{t}^{t,\cdot}$, where
$(Y_{\cdot}^{t,\cdot},Z_{\cdot}^{t,\cdot})$ is the solution of
Eq.(\ref{qi13}). Then for $t\in[0,T]$, $u(t,\cdot)$ is a weak
solution for Eq.(\ref{zhang685}). Moreover, $u(t,\cdot)$ is {\rm
a.s.} continuous w.r.t. $t$ in
$L_{\rho}^2(\mathbb{R}^d;\mathbb{R}^1)$.
\end{thm}

Then we prove the main theorem in this section.
\begin{thm}\label{qi046} Under Conditions {\rm(A.1)$'$}--{\rm(A.4)$'$}, let $u(t,\cdot)\triangleq Y_{t}^{t,\cdot}$, where
$(Y_{\cdot}^{t,\cdot},Z_{\cdot}^{t,\cdot})$ is the solution of
Eq.(\ref{qi13}). Then $u(t,\cdot)$ has an indistinguishable version
which is a "perfect" stationary solution of Eq.(\ref{zhang685}).
\end{thm}
{\em Proof}. For $Y\in
L_{\rho}^2({\mathbb{R}^{d}};{\mathbb{R}^{1}})$, $Z\in
L_{\rho}^2({\mathbb{R}^{d}};{\mathbb{R}^{d}})$, let
\begin{eqnarray*}
\hat{f}(\mathcal{T},Y,Z)=f(X_{s}^{t,\cdot},Y,Z),\ \ \
\hat{g}(\mathcal{T},Y,Z)=g(X_{s}^{t,\cdot},Y,Z).
\end{eqnarray*}
Here we take $\mathcal{T}=(s,t)$ as a dual time variable (t is
fixed). By Condition (A.1)$'$, we have
\begin{eqnarray*}
&&\|\hat{f}(\mathcal{T},Y,Z)\|^2_{L_{\rho}^2({\mathbb{R}^{d}};{\mathbb{R}^{1}})}\\
&=&\int_{\mathbb{R}^{d}}|f(X_{s}^{t,x},Y(x),Z(x))|^2\rho^{-1}(x)dx\\
&\leq&C_p\int_{\mathbb{R}^{d}}|f(X_{s}^{t,x},0,0)|^2\rho^{-1}(x)dx+C_p\int_{\mathbb{R}^{d}}|Y(x)|^2\rho^{-1}(x)dx+C_p\int_{\mathbb{R}^{d}}|Z(x)|^2\rho^{-1}(x)dx.
\end{eqnarray*}
Here and in the following, $C_p$ is a generic constant. By Lemma
\ref{qi045} and Condition (A.2)$'$,
\begin{eqnarray*}
E[\int_{0}^{\infty}\int_{\mathbb{R}^{d}}{\rm e}^{-Ks}|f(X_{s}^{t,x},0,0)|^2\rho^{-1}(x)dxds]&\leq&C_p\int_{0}^{\infty}\int_{\mathbb{R}^{d}}{\rm e}^{-Ks}|f(x,0,0)|^2\rho^{-1}(x)dxds\\
&\leq&C_p\int_{\mathbb{R}^{d}}|f(x,0,0)|^p\rho^{-1}(x)dx<\infty.
\end{eqnarray*}
We take
$\tilde{f}(\mathcal{T})=(\int_{\mathbb{R}^{d}}|f(X_{s}^{t,x},0,0)|^2\rho^{-1}(x)dx)^{1\over2}$,
then $\hat{f}(\mathcal{T},Y,Z)$ satisfies Condition (A.1). Similarly
we can also prove $\hat{g}(\mathcal{T},Y,Z)$ satisfies Condition
(A.1). On the other hand, applying $\hat{\theta}_r$ on
$\hat{f}(\mathcal{T},Y,Z)$, by (\ref{qi18}), we have for any $Y\in
L_{\rho}^2({\mathbb{R}^{d}};{\mathbb{R}^{1}})$ and $Z\in
L_{\rho}^2({\mathbb{R}^{d}};{\mathbb{R}^{d}})$,
\begin{eqnarray*}
\hat{\theta}_r\circ\hat{f}(\mathcal{T},Y,Z)=f(\hat{\theta}_r\circ
X_{s}^{t,\cdot},Y,Z)=f(X_{s+r}^{t+r,\cdot},Y,Z).
\end{eqnarray*}
Verifying $\hat{g}(\mathcal{T},Y,Z)$ in the same way, we know that
$\hat{f}(\mathcal{T},Y,Z)$ and $\hat{g}(\mathcal{T},Y,Z)$ satisfy
Condition (A.2). Since by Theorem \ref{qi043}, Eq.(\ref{qi13}) has a
unique solution $(Y_\mathcal {T},Z_\mathcal {T})$, this $(Y_\mathcal
{T},Z_\mathcal {T})$ is a stationary solution as a consequence of
Proposition \ref{qi031}. That is to say for any $t\geq0$
\begin{eqnarray*}
\hat{\theta}_r\circ Y_\mathcal {T}=\hat{\theta}_r\circ
Y_{s}^{t,\cdot}=Y_{s+r}^{t+r,\cdot}, \ \ \hat{\theta}_r\circ
Z_\mathcal {T}=\hat{\theta}_r\circ
Z_{s}^{t,\cdot}=Z_{s+r}^{t+r,\cdot}\ \ \ {\rm for}\ {\rm all}\
r\geq0,\ s\geq t\ {\rm a.s.}.
\end{eqnarray*}
In particular, for any $t\geq0$
\begin{eqnarray}\label{qi19}
\hat{\theta}_r\circ Y_{t}^{t,\cdot}=Y_{t+r}^{t+r,\cdot}\ \ \ {\rm
for}\ {\rm all}\ r\geq0\ {\rm a.s.}.
\end{eqnarray}
By Theorem \ref{qi044}, we know that $u(t,\cdot)\triangleq
Y_{t}^{t,\cdot}$ is the weak solution for Eq.(\ref{zhang685}), so we
get from (\ref{qi19}) that for any $t\geq0$
\begin{eqnarray*}
\hat{\theta}_r\circ u({t,\cdot})=u({t+r,\cdot})\ \ \ {\rm for}\ {\rm
all}\ r\geq0\ {\rm a.s.}.
\end{eqnarray*}
Until now, we know "crude" stationary property for $u({t,\cdot})$.
And by Theorem \ref{qi044}, $u(t,\cdot)$ is continuous w.r.t. $t$,
So we can get an indistinguishable version of $u({t,\cdot})$, still
denoted by $u({t,\cdot})$, s.t.
\begin{eqnarray*}
\hat{\theta}_{r}\circ u({t,\cdot})=u({t+r,\cdot})\ \ \ {\rm for}\
{\rm all}\ t,\ r\geq0\ \ \rm{a.s.}.
\end{eqnarray*}
So we proved the desired result. $\hfill\diamond$
\\

By Definition \ref{qi042}, Conditions (A.1)$'$ and (A.2)$'$, one can
calculate that $g\big(\cdot,u(s,\cdot),(\sigma^*\nabla
u)(s,\cdot)\big)\\\in\mathcal{L}^2_{U_0}(L_{\rho}^2({\mathbb{R}^{d}};{\mathbb{R}^{1}}))$
is locally square integrable in $[0,T]$. Now we consider
Eq.(\ref{zhao21}) with cylindrical Brownian motion $B$ on $U_0$. For
arbitrary $T>0$, let $Y$ be the solution of Eq.(\ref{qi13}) and
$u(t,\cdot)=Y_t^{t,\cdot}$ be the stationary solution of
Eq.(\ref{zhang685}) with $\hat{B}$ chosen as the time reversal of
$B$ from time $T$, i.e. $\hat{B}_s=B_{T-s}-B_{T}$ or
$\hat{{\beta}}_j(s)={{\beta}}_j(T-s)-{{\beta}}_j(T)$ for $s\geq0$.
By (\ref{bmrelation}) and integral transformation in
Eq.(\ref{zhang685}), we can see that $v(t,x)\triangleq u(T-t,x)$
satisfies (\ref{zhao21}) or its equivalent form
\begin{eqnarray}\label{zhang686}
&&v(t,x)=v(t,v_0)(x)=v_0(x)+\int_{0}^{t}[\mathscr{L}v(s,x)+f\big(x,v(s,x),(\sigma^*\nabla v)(s,x)\big)]ds\nonumber\\
&&\ \ \ \ \ \ \ \ \ \ \ \ \ \ \ \ \ \ \ \ \ \ \ \ \ \ \ \ \
+\sum_{j=1}^{\infty}\int_{0}^{t}g_j\big(x,v(s,x),(\sigma^*\nabla
v)(s,x)\big)d{{\beta}}_j(s),\ \ t\geq0.
\end{eqnarray}
Here $v_0(x)=v(0,x)$.

In fact, we can prove a claim that
$v(t,\cdot)(\omega)=Y_{T-t}^{T-t,\cdot}(\hat{\omega})$ does not
depend on the choice of $T$. For this, we only need to show that for
any $T'\geq T$,
$Y_{T-t}^{T-t,\cdot}(\hat{\omega})=Y_{T'-t}^{T'-t,\cdot}({\hat{\omega}}')$
when $0\leq t\leq T$, where $\hat{\omega}(s)={B}_{T-s}-{B}_{T}$ and
${\hat{\omega}}'(s)={B}_{T'-s}-{B}_{T'}$. Let $\hat{\theta}_\cdot$
and $\hat{\theta}'_\cdot$ be the shifts of $\hat{\omega}(\cdot)$ and
${\hat{\omega}}'(\cdot)$ respectively. Since by (\ref{qi19}), we
have
\begin{eqnarray*}
&&Y_{T-t}^{T-t,\cdot}(\hat{\omega})={\hat{\theta}}_{T-t}Y_{0}^{0,\cdot}(\hat{\omega})=Y_{0}^{0,\cdot}({\hat{\theta}}_{T-t}\hat{\omega}),\\
&&Y_{T'-t}^{T'-t,\cdot}({\hat{\omega}}')=\hat{\theta}'_{T'-t}Y_{0}^{0,\cdot}({\hat{\omega}}')=Y_{0}^{0,\cdot}(\hat{\theta}'_{T'-t}{\hat{\omega}}').
\end{eqnarray*}
So we just need to assert that
${\hat{\theta}}_{T-t}{\hat{\omega}}=\hat{\theta}'_{T'-t}{\hat{\omega}}'$.
Indeed we have for any $s\geq0$
\begin{eqnarray*}
({\hat{\theta}}_{T-t}{\hat{\omega}})(s)&=&{\hat{\omega}}(T-t+s)-{\hat{\omega}}(T-t)\\
&=&({B}_{T-(T-t+s)}-{B}_{T})-({B}_{T-(T-t)}-{B}_{T})\\
&=&{B}_{t-s}-{B}_{t}.
\end{eqnarray*}
Note that the right hand side of the above formula does not depend
on $T$, therefore
${\hat{\theta}}_{T-t}{\hat{\omega}}(s)={\hat{\theta}'}_{T'-t}{\hat{\omega}}'(s)={B}_{t-s}-{B}_{t}$.

On probability space $(\Omega,\mathscr{F},P)$, we define
${\theta}_{t}=(\hat{\theta}_{t})^{-1}$, $t\geq0$. Actually $\hat{B}$
is a two-sided Brownian motion, so
$(\hat{\theta}_{t})^{-1}=\hat{\theta}_{-t}$ is well defined (see
\cite{ar}). It is easy to see that ${\theta}_{t}$ is a shift w.r.t.
${B}$ satisfying
\begin{description}
\item[$(\textrm{i})$]$P\cdot({\theta}_{t})^{-1}=P$;
\item[$(\textrm{i}\textrm{i})$]${\theta}_{0}=I$;
\item[$(\textrm{i}\textrm{i}\textrm{i})$]${\theta}_{s}\circ{\theta}_{t}={\theta}_{s+t}$;
\item[$(\textrm{iv})$]${\theta}_{t}\circ{B}_s={B}_{s+t}-{B}_{t}$.
\end{description}

Since
$v(t,\cdot)(\omega)=u(T-t,\cdot)(\hat{\omega})=Y_{T-t}^{T-t,\cdot}(\hat{\omega})$
a.s., so
\begin{eqnarray*}
{\theta}_rv(t,\cdot)(\omega)=\hat{\theta}_{-r}u(T-t,\cdot)(\hat{\omega})=u(T-t-r,\cdot)(\hat{\omega})=v(t+r,\cdot)(\omega),
\end{eqnarray*}
for all $r\geq0$ and $T\geq t+r$ a.s.. In particular, let
$Y(\omega)=v_0(\omega)=Y_{T}^{T,\cdot}(\hat{\omega})$. Then above
formula implies (\ref{zhao001}):
\begin{eqnarray*}
{\theta}_tY(\omega)=Y({\theta}_t\omega)=v(t,\omega)=v(t,v_0(\omega),\omega)=v(t,Y(\omega),\omega),\
{\rm for}\ {\rm all}\ t\geq0\ \rm{a.s.}.
\end{eqnarray*}
That is to say
$v(t,\cdot)(\omega)=Y({\theta}_t\omega)(\cdot)=Y_{T-t}^{T-t,\cdot}(\hat{\omega})$
is a stationary solution of Eq.(\ref{zhao21}) w.r.t. ${\theta}$.
Therefore we proved the following theorem
\begin{thm}\label{qi047} Under Conditions {\rm(A.1)$'$}--{\rm(A.4)$'$}, for arbitrary $T$ and $t\in[0,T]$, let $v(t,\cdot)\triangleq Y_{T-t}^{T-t,\cdot}$, where $(Y_{\cdot}^{t,\cdot},Z_{\cdot}^{t,\cdot})$ is the solution of Eq.(\ref{qi13}) with $\hat{B}_s={B}_{T-s}-{B}_T$ for all $s\geq0$. Then $v(t,\cdot)$ is
a "perfect" stationary solution of Eq.(\ref{zhao21}).
\end{thm}

\section{Finite horizon BDSDEs}
\setcounter{equation}{0}

Before we study the BDSDEs on infinite horizon, we need to study the
BDSDEs on finite horizon and establish the connection with SPDEs.
For finite dimensional noise and under Lipschitz condition for a.e.
$x\in\mathbb{R}^{d}$, the problem was studied in Bally and Matoussi
\cite{ba-ma}. In this section, we consider the following BDSDE with
infinite dimensional noise on finite horizon:
\begin{eqnarray}\label{qi20}
Y_{s}^{t,x}&=&h(X_{T}^{t,x})+\int_{s}^{T}f(r,X_{r}^{t,x},Y_{r}^{t,x},Z_{r}^{t,x})dr\nonumber\\
&&-\int_{s}^{T}g(r,X_{r}^{t,x},Y_{r}^{t,x},Z_{r}^{t,x})d^\dagger{\hat{B}}_r-\int_{s}^{T}\langle
Z_{r}^{t,x},dW_r\rangle,\ \ \ 0\leq s\leq T.
\end{eqnarray}
Here $h:\Omega\times\mathbb{R}^{d}\longrightarrow{\mathbb{R}^1}$,
$f:[0,T]\times\mathbb{R}^{d}\times\mathbb{R}^1\times\mathbb{R}^{d}{\longrightarrow{\mathbb{R}^1}}$,
$g:[0,T]\times\mathbb{R}^{d}\times\mathbb{R}^1\times\mathbb{R}^{d}\longrightarrow
{\mathcal{L}^2_{U_0}(\mathbb{R}^{1})}$. Set $g_j\triangleq
g\sqrt{\lambda_j}e_j:[0,T]\times\mathbb{R}^{d}\times\mathbb{R}^1\times\mathbb{R}^{d}{\longrightarrow{\mathbb{R}^1}}$,
then Eq.(\ref{qi20}) is equivalent to
\begin{eqnarray*}\label{qi21}
Y_{s}^{t,x}&=&h(X_{T}^{t,x})+\int_{s}^{T}f(r,X_{r}^{t,x},Y_{r}^{t,x},Z_{r}^{t,x})dr\nonumber\\
&&-\sum_{j=1}^{\infty}\int_{s}^{T}g_j(r,X_{r}^{t,x},Y_{r}^{t,x},Z_{r}^{t,x})d^\dagger{\hat{\beta}}_j(r)-\int_{s}^{T}\langle
Z_{r}^{t,x},dW_r\rangle,\ \ \ 0\leq s\leq T.
\end{eqnarray*}
We assume
\begin{description}
\item[(H.1).] Function $h$ is $\mathscr{F}_{T,\infty}^{\hat{B}}\otimes\mathscr{B}_{\mathbb{R}^{d}}$ measurable and $E[\int_{\mathbb{R}^{d}}|h(x)|^2\rho^{-1}(x)dx]<\infty$;
\item[(H.2).] Functions $f$ and $g$ are $\mathscr{B}_{[0,T]}\otimes\mathscr{B}_{\mathbb{R}^{d}}\otimes\mathscr{B}_{\mathbb{R}^{1}}\otimes\mathscr{B}_{\mathbb{R}^{d}}$ measurable and there exist constants
$C,C_j,\alpha_j\geq0$ with $\sum_{j=1}^\infty C_j<\infty$ and
$\sum_{j=1}^\infty\alpha_j<{1\over2}$ s.t. for any $t\in[0,T]$,
$Y_1,Y_2\in L_{\rho}^2({\mathbb{R}^{d}};{\mathbb{R}^{1}})$,
$X,Z_1,Z_2\in L_{\rho}^2({\mathbb{R}^{d}};{\mathbb{R}^{d}})$
\begin{eqnarray*}
&&\int_{\mathbb{R}^{d}}|f(t,X(x),Y_1(x),Z_1(x))-f(t,X(x),Y_2(x),Z_2(x))|^2\rho^{-1}(x)dx\\
&\leq&C\int_{\mathbb{R}^{d}}(|Y_1(x)-Y_2(x)|^2+|Z_1(x)-Z_2(x)|^2)\rho^{-1}(x)dx,\\
&&\int_{\mathbb{R}^{d}}|g_j(t,X(x),Y_1(x),Z_1(x))-g_j(t,X(x),Y_2(x),Z_2(x))|^2\rho^{-1}(x)dx\\
&\leq&\int_{\mathbb{R}^{d}}(C_j|Y_1(x)-Y_2(x)|^2+{\alpha_j}|Z_1(x)-Z_2(x)|^2)\rho^{-1}(x)dx;
\end{eqnarray*}
\item[(H.3).] $\int_{0}^{T}\int_{\mathbb{R}^d}|f(s,x,0,0)|^2\rho^{-1}(x)dxds<\infty$ and $\int_{0}^{T}\int_{\mathbb{R}^{d}}\|g(s,x,0,0)\|^2_{\mathcal{L}^2_{U_0}(\mathbb{R}^{1})}\rho^{-1}(x)dxds<\infty$;
\item[(H.4).] $b\in C_{l,b}^2(\mathbb{R}^{d};\mathbb{R}^{d})$, $\sigma\in
C_{l,b}^3(\mathbb{R}^{d};\mathbb{R}^{d}\times\mathbb{R}^{d})$.
\end{description}

Needless to say, the conditions (H.1)-(H.4) for the existence and
uniqueness of solution of Eq.(\ref{qi20}) are weaker than what are
needed for the case of infinite horizon. We would like to point out
that for the finite horizon problem, our conditions are weaker than
those in Bally and Matoussi \cite{ba-ma}. In
(H.1), we allow the terminal function $h$ depending on
$\mathscr{F}_{t,T}$ independent sigma field
$\mathscr{F}_{T,\infty}^{\hat{B}}$. One can easily verify that it
doesn't affect the results in \cite{ba-ma}. Moreover, here we only
need Lipschitz condition in the space
$L_{\rho}^2({\mathbb{R}^{d}};{\mathbb{R}^{1}})$ instead of the
pathwise Lipschitz condition posed in \cite{ba-ma}.
\begin{defi}\label{qi051}
A pair of processes $(Y_{\cdot}^{t,\cdot},Z_{\cdot}^{t,\cdot})\in
S^{2,0}([0,T];L_{\rho}^2({\mathbb{R}^{d}};{\mathbb{R}^{1}}))\bigotimes
M^{2,0}([0,T];L_{\rho}^2({\mathbb{R}^{d}};{\mathbb{R}^{d}}))$ is
called a solution of Eq.(\ref{qi20}) if for any $\varphi\in
C_c^{0}(\mathbb{R}^d;\mathbb{R}^1)$,
\begin{eqnarray}\label{qi22}
\int_{\mathbb{R}^{d}}Y_s^{t,x}\varphi(x)dx&=&\int_{\mathbb{R}^{d}}h(X_{T}^{t,x})\varphi(x)dx+\int_{s}^{T}\int_{\mathbb{R}^{d}}f(r,X_{r}^{t,x},Y_r^{t,x},Z_r^{t,x})\varphi(x)dxdr\nonumber\\
&&-\sum_{j=1}^{\infty}\int_{s}^{T}\int_{\mathbb{R}^{d}}g_j(r,X_{r}^{t,x},Y_r^{t,x},Z_r^{t,x})\varphi(x)dxd^\dagger{\hat{\beta}}_j(r)\nonumber\\
&&-\int_{s}^{T}\langle
\int_{\mathbb{R}^{d}}Z_r^{t,x}\varphi(x)dx,dW_r\rangle\ \ \ P-{\rm
a.s.}.
\end{eqnarray}
\end{defi}

The main objective of this section is to prove
\begin{thm}\label{qi052} Under Conditions {\rm(H.1)}--{\rm(H.4)},
Eq.(\ref{qi20}) has a unique solution.
\end{thm}

This theorem is an extension of Theorem 3.1 in \cite{ba-ma}. The
idea is to start from Bally and Matoussi's results for finite
dimensional noise and then take limit to obtain the solution for the
case of infinite dimensional noise. But Bally and Matoussi's results
cannot apply immediately here as we have a weaker Lipschitz
condition and some of the key claims in the proof of Theorem 3.1
(\cite{ba-ma}) are not obvious under their conditions. Moreover, the
result $Y_{\cdot}^{t,\cdot}\in
S^{2,0}([0,T];L_{\rho}^2({\mathbb{R}^{d}};{\mathbb{R}^{1}}))$ was
not obtained in \cite{ba-ma}. We study a sequence of BDSDEs
\begin{eqnarray}\label{zhang66100}
Y_s^{t,x,n}=&&h(X_{T}^{t,x})+\int_{s}^{T}f(r,X_{r}^{t,x},Y^{t,x,n}_r,Z^{t,x,n}_r)dr\nonumber\\
&&-\sum_{j=1}^{n}\int_{s}^{T}g_j(r,X_{r}^{t,x},Y^{t,x,n}_r,Z^{t,x,n}_r)d^\dagger{\hat{\beta}}_j(r)-\int_{s}^{T}\langle
Z^{t,x,n}_r,dW_r\rangle.
\end{eqnarray}
A solution of (\ref{zhang66100}) is a pair of processes
$(Y_\cdot^{t,\cdot,n},Z_\cdot^{t,\cdot,n})\in
S^{2,0}([0,T];L_{\rho}^2({\mathbb{R}^{d}};{\mathbb{R}^{1}}))\bigotimes
M^{2,0}([0,T];\\L_{\rho}^2({\mathbb{R}^{d}};{\mathbb{R}^{d}}))$
satisfying the spatial integral form of Eq.(\ref{zhang66100}), i.e.
(\ref{qi22}) with a finite number of one dimensional backward
stochastic integrals.

First we do some preparations.
\begin{lem}\label{qi1a00} Under Conditions {\rm(H.1)}--{\rm(H.4)}, if there exists $(Y_\cdot(\cdot),Z_\cdot(\cdot))\in
M^{2,0}([t,T];L_{\rho}^2({\mathbb{R}^{d}};{\mathbb{R}^{1}}))\\\bigotimes
M^{2,0}([t,T];L_{\rho}^2({\mathbb{R}^{d}};{\mathbb{R}^{d}}))$
satisfying the spatial integral form of Eq.(\ref{zhang66100}) for
$t\leq s\leq T$, then $Y_\cdot(\cdot)\in
S^{2,0}([t,T];L_{\rho}^2({\mathbb{R}^{d}};{\mathbb{R}^{1}}))$ and
therefore $(Y_s(x),Z_s(x))$ is a solution of Eq.(\ref{zhang66100}).
\end{lem}
{\em Proof}. Let's first see ${Y}_{s}(\cdot)$ is continuous w.r.t.
$s$ in $L_{\rho}^2({\mathbb{R}^{d}};{\mathbb{R}^{1}})$. Since
$(Y_s(x),Z_s(x))$ satisfies the form of Eq.(\ref{zhang66100}) for
$t\leq s<T$, a.e. $x\in\mathbb{R}^{d}$, therefore,
\begin{eqnarray*}\label{zhangvvv}
&&\int_{\mathbb{R}^{d}}|Y_{s+\triangle
s}(x)-Y_s(x)|^2\rho^{-1}(x)dx\nonumber\\
&\leq&C_p\int_{\mathbb{R}^{d}}\int_s^{s+\triangle
s}|f(r,X_{r}^{t,x},Y_r(x),Z_r(x))|^2dr\rho^{-1}(x)dx\nonumber\\
&&+C_p\sum_{j=1}^{n}\int_{\mathbb{R}^{d}}|\int_s^{s+\triangle
s}g_j(r,X_{r}^{t,x},Y_r(x),Z_r(x))d^\dagger{\hat{\beta}}_j(r)|^2\rho^{-1}(x)dx\nonumber\\
&&+C_p\int_{\mathbb{R}^{d}}|\int_s^{s+\triangle s}\langle
Z_r(x),dW_r\rangle|^2\rho^{-1}(x)dx.
\end{eqnarray*}
For the forward stochastic integral part, it is trivial to see that
for $0\leq\triangle s\leq T-s$, $|\int_s^{s+\triangle s}\langle
Z_r(x),dW_r\rangle|^2\leq\sup_{0\leq\triangle s\leq
T-s}|\int_s^{s+\triangle s}\langle Z_r(x),dW_r\rangle|^2$ a.s.. And
we can deduce that $\int_{\mathbb{R}^{d}}\sup_{0\leq\triangle s\leq
T-s}|\int_s^{s+\triangle s}\langle
Z_r(x),dW_r\rangle|^2\rho^{-1}(x)dx<\infty$ a.s. by the B-D-G
inequality and $Z_\cdot(\cdot)\in
M^{2,0}([t,T];L_{\rho}^2({\mathbb{R}^{d}};{\mathbb{R}^{d}}))$. So by
the dominated convergence theorem, $\lim_{\triangle
s\rightarrow0^{+}}\int_{\mathbb{R}^{d}}|\int_s^{s+\triangle
s}\langle Z_r(x),\\dW_r\rangle|^2\rho^{-1}(x)dx=0$. Similarly we can
prove $\lim_{\triangle
s\rightarrow0^{-}}\int_{\mathbb{R}^{d}}|\int_{s+\triangle
s}^s\langle Z_r(x),dW_r\rangle|^2\rho^{-1}(x)dx=0$ for $t<s\leq T$.
The backward stochastic integral part tends to $0$ as $\triangle
s\rightarrow0$ can be deduced similarly. So $Y_s(\cdot)$ is
continuous w.r.t. $s$ in
$L_{\rho}^2({\mathbb{R}^{d}};{\mathbb{R}^{1}})$. From Conditions
(H.2)--(H.4) and $(Y_\cdot(\cdot),Z_\cdot(\cdot))\in
M^{2,0}([t,T];L_{\rho}^2({\mathbb{R}^{d}};{\mathbb{R}^{1}}))\bigotimes
M^{2,0}([t,T];L_{\rho}^2({\mathbb{R}^{d}};{\mathbb{R}^{d}}))$, it
follows that for a.e. $x\in{\mathbb{R}^{d}}$,
$E[\int_{t}^{T}|f(r,X_{r}^{t,x},Y_r(x),Z_r(x))|^2dr]<\infty$ and
$\sum_{j=1}^{n}E[\int_{t}^{T}|g_j(r,X_{r}^{t,x},Y_r(x),Z_r(x))|^2dr]<\infty$.
For a.e. $x\in{\mathbb{R}^{d}}$, referring to Lemma 1.4 in
\cite{pa-pe3}, we use the generalized It$\hat {\rm o}$'s formula
(c.f. Elworthy, Truman and Zhao \cite{el-tr-zh}) to
$\psi_M\big(Y_r(x)\big)$, where $\psi_M(x)=x^2I_{\{-M\leq
x<M\}}+2M(x-M)I_{\{x\geq M\}}-2M(x+M)I_{\{x<-M\}}$. Then
\begin{eqnarray}\label{zhaoaaaa}
&&\psi_M(Y_s(x))+\int_{s}^{T}I_{\{-M\leq Y_r(x)<M\}}|Z_r(x)|^2dr\nonumber\\
&=&\psi_M\big(h(X_{T}^{t,x})\big)+\int_{s}^{T}\psi_M^{'}(Y_r(x))f(r,X_{r}^{t,x},Y_r(x),Z_r(x))dr\nonumber\\
&&+\sum_{j=1}^{n}\int_{s}^{T}I_{\{-M\leq Y_r(x)<M\}}|g_j(r,X_{r}^{t,x},Y_r(x),Z_r(x))|^2dr\\
&&-\sum_{j=1}^{n}\int_{s}^{T}\psi_M^{'}(Y_r(x))g_j(r,X_{r}^{t,x},Y_r(x),Z_r(x))d^\dagger{\hat{\beta}}_j(r)-\int_{s}^{T}\langle\psi_M^{'}(Y_r(x))Z_r(x),dW_r\rangle\nonumber.
\end{eqnarray}
We can use the Fubini theorem to perfect (\ref{zhaoaaaa}) so that (\ref{zhaoaaaa})
is satisfied for a.e. $x\in \mathbb{R}^d$, on a full measure set that is independent of $x$. Taking integration in
${\mathbb{R}^{d}}$ on both sides, applying the stochastic Fubini
theorem (\cite{pr-za1}), we have
\begin{eqnarray*}
&&\int_{\mathbb{R}^{d}}\psi_M(Y_s(x))\rho^{-1}(x)dx+\int_{s}^{T}\int_{\mathbb{R}^{d}}I_{\{-M\leq Y_r(x)<M\}}|Z_r(x)|^2\rho^{-1}(x)dxdr\\
&\leq&\int_{\mathbb{R}^{d}}\psi_M\big(h(X_{T}^{t,x})\big)\rho^{-1}(x)dx+\int_{s}^{T}\int_{\mathbb{R}^{d}}\psi_M^{'}(Y_r(x))f(r,X_{r}^{t,x},0,0)\rho^{-1}(x)dxdr\nonumber\\
&&+\int_{s}^{T}\int_{\mathbb{R}^{d}}\psi_M^{'}(Y_r(x))\big(f(r,X_{r}^{t,x},Y_r(x),Z_r(x))-f(r,X_{r}^{t,x},0,0)\big)\rho^{-1}(x)dxdr\nonumber\\
&&+C_p\sum_{j=1}^{n}\int_{s}^{T}\int_{\mathbb{R}^{d}}|g_j(r,X_{r}^{t,x},Y_r(x),Z_r(x))-g_j(r,X_{r}^{t,x},0,0)|^2\rho^{-1}(x)dxdr\nonumber\\
&&+C_p\sum_{j=1}^{n}\int_{s}^{T}\int_{\mathbb{R}^{d}}|g_j(r,X_{r}^{t,x},0,0)|^2\rho^{-1}(x)dxdr-\int_{s}^{T}\langle\int_{\mathbb{R}^{d}}\psi_M^{'}(Y_r(x))Z_r(x)\rho^{-1}(x)dx,dW_r\rangle\nonumber\\
&&-\sum_{j=1}^{n}\int_{s}^{T}\int_{\mathbb{R}^{d}}\psi_M^{'}(Y_r(x))g_j(r,X_{r}^{t,x},Y_r(x),Z_r(x))\rho^{-1}(x)dxd^\dagger{\hat{\beta}}_j(r).
\end{eqnarray*}
Noting that $\psi_M\big(h(X_{T}^{t,x})\big)\leq|h(X_{T}^{t,x})|^2$
and $|\psi_M^{'}(Y_r(x))|^2\leq4|Y_r(x)|^2$, so by Lemma
\ref{qi045}, the B-D-G inequality and Cauchy-Schwartz inequality, we
have
\begin{eqnarray}\label{zhangzzz00}
&&E[\sup_{t\leq s\leq T}\int_{\mathbb{R}^{d}}\psi_M(Y_s(x))\rho^{-1}(x)dx]\nonumber\\
&\leq&C_pE[\int_{\mathbb{R}^{d}}{|h(x)|^2}\rho^{-1}(x)dx]+C_pE[\int_{t}^{T}\int_{\mathbb{R}^{d}}(|Y_r(x)|^2+|Z_r(x)|^2)\rho^{-1}(x)dxdr]\nonumber\\
&&+C_pE[\sum_{j=1}^{n}\int_{t}^{T}\int_{\mathbb{R}^{d}}(|g_j(r,x,0,0)|^2+|f(r,x,0,0)|^2)\rho^{-1}(x)dxdr]\nonumber\\
&&+C_pE[\sqrt{\int_{t}^{T}\int_{\mathbb{R}^{d}}|\psi_M^{'}(Y_s(x))|^2\rho^{-1}(x)dx\int_{\mathbb{R}^{d}}\sum_{j=1}^{n}|g_j(r,X_{r}^{t,x},Y_r(x),Z_r(x))|^2\rho^{-1}(x)dxdr}]\nonumber\\
&&+C_pE[\sqrt{\int_{t}^{T}\int_{\mathbb{R}^{d}}|\psi_M^{'}(Y_s(x))|^2\rho^{-1}(x)dx\int_{\mathbb{R}^{d}}|Z_r(x)|^2\rho^{-1}(x)dxdr}]\nonumber\\
&\leq&C_pE[\int_{\mathbb{R}^{d}}{|h(x)|^2}\rho^{-1}(x)dx]+C_pE[\int_{t}^{T}\int_{\mathbb{R}^{d}}(|Y_r(x)|^2+|Z_r(x)|^2)\rho^{-1}(x)dxdr]\nonumber\\
&&+C_pE[\sum_{j=1}^{n}\int_{t}^{T}\int_{\mathbb{R}^{d}}(|g_j(r,x,0,0)|^2+|f(r,x,0,0)|^2)\rho^{-1}(x)dxdr]\nonumber\\
&&+{1\over5}E[\sup_{t\leq s\leq
T}\int_{\mathbb{R}^{d}}|\psi_M^{'}(Y_s(x))|^2\rho^{-1}(x)dx].
\end{eqnarray}
Since $(Y_\cdot(\cdot),Z_\cdot(\cdot))\in
M^{2,0}([t,T];L_{\rho}^2({\mathbb{R}^{d}};{\mathbb{R}^{1}}))\bigotimes
M^{2,0}([t,T];L_{\rho}^2({\mathbb{R}^{d}};{\mathbb{R}^{d}}))$,
taking the limit as $M\rightarrow\infty$ and applying the monotone
convergence theorem, we have $E[\sup_{t\leq s\leq
T}\int_{\mathbb{R}^{d}}|Y_s(x)|^2\rho^{-1}(x)dx]<\infty$. So
$Y_\cdot(\cdot)\in
S^{2,0}([t,T];L_{\rho}^2({\mathbb{R}^{d}};{\mathbb{R}^{1}}))$
follows. That is to say $(Y_s(x),Z_s(x))$ is a solution of
Eq.(\ref{zhang66100}). $\hfill\diamond$\\

For the rest of our paper, we will leave out the similar
localization argument as in the proof of Lemma \ref{qi1a00} when
applying It$\hat {\rm o}$'s formula to save the space of this paper.
\begin{prop}\label{000000}
Under Conditions {\rm(H.1)}--{\rm(H.4)}, assume
Eq.(\ref{zhang66100}) have a unique solution
$(Y_r^{t,x,n},Z_r^{t,x,n})$, then for any $t\leq s\leq T$,
$Y_r^{s,X_s^{t,x},n}=Y_r^{t,x,n}$ and
$Z_r^{s,X_s^{t,x},n}=Z_r^{t,x,n}$ for any $r\in[s,T]$ and a.e.
$x\in{\mathbb{R}^{d}}$ a.s..
\end{prop}
{\em Proof}. For $t\leq s\leq r\leq T$, note that
$(Y_r^{s,\cdot,n},Z_r^{s,\cdot,n})$ is
$\mathscr{F}^{\hat{B}}_{r,\infty}\otimes\mathscr{F}^W_{s,r}$
measurable, so is independent of $\mathscr{F}^W_{t,s}$. Thus by
Lemma \ref{qi045}, we have
\begin{eqnarray*}
&&E[\int_{s}^{T}\int_{\mathbb{R}^{d}}(|Y_r^{s,X^{t,x}_s,n}|^2+|Z_r^{s,X^{t,x}_s,n}|^2)\rho^{-1}(x)dxdr]\nonumber\\
&\leq&C_pE[\int_{s}^{T}\int_{\mathbb{R}^{d}}(|Y_r^{s,x,n}|^2+|Z_r^{s,x,n}|^2)\rho^{-1}(x)dxdr]<\infty.
\end{eqnarray*}
Moreover, it is easy to see that $X_r^{s,X^{t,x}_s}=X^{t,x}_r$ and
$(Y_r^{s,X^{t,x}_s,n},Z_r^{s,X^{t,x}_s,n})$ is
$\mathscr{F}_{r,\infty}^{\hat{B}}\otimes\mathscr{F}_{t,r}^W$
measurable, so
$(Y_\cdot^{s,X^{t,\cdot}_s,n},Z_\cdot^{s,X^{t,\cdot}_s,n})\in
M^{2,0}([s,T];L_{\rho}^2({\mathbb{R}^{d}};{\mathbb{R}^{1}}))
\bigotimes
M^{2,0}([s,T];L_{\rho}^2({\mathbb{R}^{d}};{\mathbb{R}^{d}}))$ and
$(Y_r^{s,X^{t,x}_s,n},Z_r^{s,X^{t,x}_s,n})$ satisfies the spatial
integral form of Eq.(\ref{zhang66100}) for $s\leq r\leq T$. Define
$Y_r^{s,X^{t,x}_s,n}=Y_r^{t,x,n}$, $Z_r^{s,X^{t,x}_s,n}=Z_r^{t,x,n}$
when $t\leq r<s$. Then $(Y_r^{s,X^{t,x}_s,n},Z_r^{s,X^{t,x}_s,n})$
satisfies the spatial integral form of Eq.(\ref{zhang66100}) for
$t\leq r\leq T$ and
$(Y_\cdot^{s,X^{t,\cdot}_s,n},Z_\cdot^{s,X^{t,\cdot}_s,n})\in
M^{2,0}([t,T];L_{\rho}^2({\mathbb{R}^{d}};{\mathbb{R}^{1}}))
\bigotimes
M^{2,0}([t,T];L_{\rho}^2({\mathbb{R}^{d}};{\mathbb{R}^{d}}))$.
Therefore, by Lemma \ref{qi1a00},
$(Y_r^{s,X^{t,x}_s,n},Z_r^{s,X^{t,x}_s,n})$ is the solution of
Eq.(\ref{zhang66100}). By the uniqueness of the solution of
Eq.(\ref{zhang66100}), we have for any $s\in[t,T]$,
$(Y_r^{s,X^{t,x}_s,n},Z_r^{s,X^{t,x}_s,n})=(Y_r^{t,x,n},Z_r^{t,x,n})$
for any $r\in[s,T]$ and a.e. $x\in{\mathbb{R}^{d}}$ a.s.. $\hfill\diamond$\\

\begin{thm}\label{jia} Under Conditions {\rm(H.1)}--{\rm(H.4)}, Eq.(\ref{zhang66100}) has a unique solution, i.e. there
exists a unique $(Y^{t,\cdot,n}_{\cdot},Z^{t,\cdot,n}_{\cdot})\in
S^{2,0}([t,T];
L_{\rho}^2({\mathbb{R}^{d}};{\mathbb{R}^{1}}))\bigotimes
M^{2,0}([t,T];L_{\rho}^2({\mathbb{R}^{d}};{\mathbb{R}^{d}}))$ such
that for an arbitrary $\varphi\in
C_c^{0}(\mathbb{R}^d;\mathbb{R}^1)$
\begin{eqnarray}\label{zhang662}
\int_{\mathbb{R}^{d}}Y_s^{t,x,n}\varphi(x)dx&=&\int_{\mathbb{R}^{d}}h(X_{T}^{t,x})\varphi(x)dx+\int_{s}^{T}\int_{\mathbb{R}^{d}}f(r,X_{r}^{t,x},Y_r^{t,x,n},Z_r^{t,x,n})\varphi(x)dxdr\nonumber\\
&&-\sum_{j=1}^{n}\int_{s}^{T}\int_{\mathbb{R}^{d}}g_j(r,X_{r}^{t,x},Y_r^{t,x,n},Z_r^{t,x,n})\varphi(x)dxd^\dagger{\hat{\beta}}_j(r)\nonumber\\
&&-\int_{s}^{T}\langle\int_{\mathbb{R}^{d}}Z_r^{t,x,n}\varphi(x)dx,dW_r\rangle\
\ \ P-{\rm a.s.}.
\end{eqnarray}
\end{thm}
{\em Proof}. \underline{Uniqueness}. Assume there exists another
$(\hat{Y}_{s}^{t,x,n}$, $\hat{Z}_{s}^{t,x,n})\in
S^{2,0}([t,T];L_{\rho}^2({\mathbb{R}^{d}};{\mathbb{R}^{1}}))\bigotimes
M^{2,0}\\([t,T];L_{\rho}^2({\mathbb{R}^{d}};{\mathbb{R}^{d}}))$
satisfying (\ref{zhang662}). Define
$\bar{Y}_s^{t,x,n}=Y_{s}^{t,x,n}-\hat{Y}_{s}^{t,x,n}$ and
$\bar{Z}_s^{t,x,n}=Z_{s}^{t,x,n}-\hat{Z}_{s}^{t,x,n}$, $t\leq s\leq
T$. From Conditions (H.2)--(H.4) and
$(Y^{t,\cdot,n}_{\cdot},Z^{t,\cdot,n}_{\cdot})\in S^{2,0}([t,T];
L_{\rho}^2({\mathbb{R}^{d}};{\mathbb{R}^{1}}))\bigotimes
M^{2,0}\\([t,T];L_{\rho}^2({\mathbb{R}^{d}};{\mathbb{R}^{d}}))$, it
follows that for a.e. $x\in{\mathbb{R}^{d}}$,
$E[\int_{t}^{T}|f(r,X_{r}^{t,x},{Y}_{r}^{t,x,n},{Z}_{r}^{t,x,n})-f(r,X_{r}^{t,x},{\hat{Y}}_{r}^{t,x,n},\\{\hat{Z}}_{r}^{t,x,n})|^2dr]<\infty$
and
$\sum_{j=1}^{n}E[\int_{t}^{T}|g_j(r,X_{r}^{t,x},{Y}_{r}^{t,x,n},{Z}_{r}^{t,x,n})-g_j(r,X_{r}^{t,x},{\hat{Y}}_{r}^{t,x,n},{\hat{Z}}_{r}^{t,x,n})|^2dr]<\infty$.
For a.e. $x\in{\mathbb{R}^{d}}$, similar as in (\ref{zhaoaaaa}), we
use generalized It$\hat {\rm o}$'s formula to ${\rm
e}^{Kr}\psi_M(\bar{Y}_r^{t,x,n})$ where $K\in\mathbb{R}^1$, then
take integration in ${\mathbb{R}^{d}}\times\Omega$ on both sides and
apply the stochastic Fubini theorem. Note that the stochastic
integrals are martingales, so taking the limit as
$M\rightarrow\infty$, we have
\begin{eqnarray}\label{zhao002}
&&E[{\rm e}^{Ks}\int_{\mathbb{R}^{d}}{{|\bar{Y}_{s}^{t,x,n}|}^2}\rho^{-1}(x)dx]+(K-2C-\sum_{j=1}^{\infty}C_j-{1\over2})E[\int_{s}^{T}\int_{\mathbb{R}^{d}}{\rm e}^{Kr}|{\bar{Y}}_{r}^{t,x,n}|^2\rho^{-1}(x)dxdr]\nonumber\\
&&+({1\over2}-\sum_{j=1}^{\infty}\alpha_j)E[\int_{s}^{T}\int_{\mathbb{R}^{d}}{\rm
e}^{Kr}|{\bar{Z}}_{r}^{t,x,n}|^2\rho^{-1}(x)dxdr]\leq0.
\end{eqnarray}
All the terms on the left hand side of (\ref{zhao002}) are positive
when taking $K$ sufficiently large, so it is easy to see that for
each $s\in [t,T]$, $\bar{Y}_s^{t,x}=0$ a.e. $x\in\mathbb{R}^{d}$
a.s.. By a "standard" argument taking $s$ in the rational number
space and noting $\int_{\mathbb{R}^{d}}{\rm
e}^{Ks}|\bar{Y}_s^{t,x,n}|^2\rho^{-1}(x)dx$ is continuous w.r.t.
$s$, we have $\bar{Y}_s^{t,x,n}=0$ for all $s\in[t,T]$, a.e.
$x\in\mathbb{R}^{d}$ a.s.. Also by (\ref{zhao002}), for a.e.
$s\in[t,T]$, $\bar{Z}_s^{t,x,n}=0$ a.e. $x\in\mathbb{R}^{d}$, a.s..
We can modify the values of $Z$ at the measure zero exceptional set
of $s$ such that $\bar{Z}_s^{t,x,n}=0$ for all $s\in[t,T]$, a.e.
$x\in\mathbb{R}^{d}$ a.s..

\underline{Existence}. Step 1: We prove for the following equation:
\begin{eqnarray}\label{zhang66111}
\tilde{Y}_s^{t,x,n}=h(X_{T}^{t,x})+\int_{s}^{T}\tilde{f}(r,X_{r}^{t,x})dr-\sum_{j=1}^{n}\int_{s}^{T}\tilde{g}_j(r,X_{r}^{t,x})d^\dagger{\hat{\beta}}_j(r)-\int_{s}^{T}\langle
\tilde{Z}^{t,x,n}_r,dW_r\rangle,
\end{eqnarray}
if (H.1) and (H.4) are satisfied, and
$\tilde{f}(\cdot,X_\cdot^{t,\cdot})$,
$\tilde{g}_j(\cdot,X_\cdot^{t,\cdot})\in
M^{2,0}([t,T];L_{\rho}^2({\mathbb{R}^{d}};{\mathbb{R}^{1}}))$, then
there exists a unique solution. For this, we can first use a similar
method as in the proof of Theorem 2.1 in \cite{ba-ma} to prove there
exists
$(\tilde{Y}_{\cdot}^{t,\cdot,n},\tilde{Z}_{\cdot}^{t,\cdot,n})\in
M^{2,0}([t,T];L_{\rho}^2({\mathbb{R}^{d}};{\mathbb{R}^{1}}))\bigotimes
M^{2,0}([t,T];L_{\rho}^2({\mathbb{R}^{d}};{\mathbb{R}^{d}}))$ such
that for an arbitrary $\varphi\in
C_c^{0}(\mathbb{R}^d;\mathbb{R}^1)$
\begin{eqnarray*}\label{zhang66200}
&&\int_{\mathbb{R}^{d}}\tilde{Y}_s^{t,x,n}\varphi(x)dx=\int_{\mathbb{R}^{d}}h(X_{T}^{t,x})\varphi(x)dx+\int_{s}^{T}\int_{\mathbb{R}^{d}}\tilde{f}(r,X_{r}^{t,x})\varphi(x)dxdr\nonumber\\
&&-\sum_{j=1}^{n}\int_{s}^{T}\int_{\mathbb{R}^{d}}\tilde{g}_j(r,X_{r}^{t,x})\varphi(x)dxd^\dagger{\hat{\beta}}_j(r)-\int_{s}^{T}\langle\int_{\mathbb{R}^{d}}\tilde{Z}_r^{t,x,n}\varphi(x)dx,dW_r\rangle\
\ \ P-{\rm a.s.}.
\end{eqnarray*}
By Lemma \ref{qi1a00}, $\tilde{Y}_{\cdot}^{t,\cdot,n}\in
S^{2,0}([t,T];L_{\rho}^2({\mathbb{R}^{d}};{\mathbb{R}^{1}}))$. Then
Step 1 follows.

Step 2: Given $(Y_{s}^{t,x,n,N-1},Z_{s}^{t,x,n,N-1})\in
S^{2,0}([t,T];L_{\rho}^2({\mathbb{R}^{d}};{\mathbb{R}^{1}}))\bigotimes
M^{2,0}([t,T];L_{\rho}^2({\mathbb{R}^{d}};{\mathbb{R}^{d}}))$,
define $(Y_{s}^{t,x,n,N},Z_{s}^{t,x,n,N})$ as follows:
\begin{eqnarray}\label{jia3}
Y_{s}^{t,x,n,N}&=&h(X_{T}^{t,x})+\int_{s}^{T}f(r,X_{r}^{t,x},Y_{r}^{t,x,n,N-1},Z_{r}^{t,x,n,N-1})dr\nonumber\\
&&-\sum_{j=1}^{n}\int_{s}^{T}g_j(r,X_{r}^{t,x},Y_{r}^{t,x,n,N-1},Z_{r}^{t,x,n,N-1})d^\dagger{\hat{\beta}}_j(r)-\int_{s}^{T}\langle
Z_{r}^{t,x,n,N},dW_r\rangle.
\end{eqnarray}
Let $(Y_{r}^{t,x,n,0},Z_{r}^{t,x,n,0})=(0,0)$. By Conditions (H.1),
(H.3), (H.4) and Lemma \ref{qi045}, we know $h$, $f(r,X_r^{t,x},0,0)$ and
$g_j(r,X_r^{t,x},0,0)$ satisfy the conditions in Step 1, so
Eq.(\ref{zhang66111}) has a unique solution
$(Y_{\cdot}^{t,\cdot,n,1},Z_{\cdot}^{t,\cdot,n,1})\in
M^{2,0}([t,T];L_{\rho}^2({\mathbb{R}^{d}};{\mathbb{R}^{1}}))\bigotimes
M^{2,0}([t,T];L_{\rho}^2({\mathbb{R}^{d}};{\mathbb{R}^{d}}))$ when
$\tilde{f}(r,X_r^{t,x})=f(r,X_r^{t,x},0,0)$ and
$\tilde{g}(r,X_r^{t,x})=g(r,X_r^{t,x},0,0)$. From Proposition
\ref{000000} and the Fubini theorem, we have
$Y_{r}^{t,x,n,1}=Y_r^{r,X_r^{t,x},n,1}$ and
$Z_{r}^{t,x,n,1}=Z_r^{r,X_r^{t,x},n,1}$ for a.e. $r\in[t,T]$,
$x\in{\mathbb{R}^{d}}$ a.s.. Thus by Conditions (H.1)--(H.4) and
Lemma \ref{qi045}, we have $h$,
$f(r,X_r^{t,x},Y_{r}^{t,x,n,1},Z_{r}^{t,x,n,1})=f(r,X_r^{t,x},Y_r^{r,X_r^{t,x},n,1},Z_r^{r,X_r^{t,x},n,1})$
and
$g_j(r,X_r^{t,x},Y_{r}^{t,x,n,1},Z_{r}^{t,x,n,1})=g_j(r,X_r^{t,x},Y_r^{r,X_r^{t,x},n,1},Z_r^{r,X_r^{t,x},n,1})$
satisfy the conditions in Step 1. Following the same procedure, we
obtain $(Y_{\cdot}^{t,\cdot,n,2},Z_{\cdot}^{t,\cdot,n,2})\in
M^{2,0}([t,T];L_{\rho}^2({\mathbb{R}^{d}};{\mathbb{R}^{1}}))\bigotimes
M^{2,0}([t,T];L_{\rho}^2({\mathbb{R}^{d}};{\mathbb{R}^{d}}))$. In
general, we see (\ref{jia3}) is an iterated mapping from
$S^{2,0}([t,T];L_{\rho}^2({\mathbb{R}^{d}};{\mathbb{R}^{1}}))\bigotimes
M^{2,0}([t,T];L_{\rho}^2({\mathbb{R}^{d}};{\mathbb{R}^{d}}))$ to
itself and obtain a sequence
$\{({Y}_r^{t,x,n,i},{Z}_r^{t,x,n,i})\}_{i=0,1,2\cdot\cdot\cdot}$. We
will prove that (\ref{jia3}) is a contraction mapping. For this,
define
\begin{eqnarray*}
&&\bar{Y}_s^{t,x,n,i}={Y}_s^{t,x,n,i}-{Y}_s^{t,x,n,i-1}, \ \
\bar{Z}_s^{t,x,n,i}={Z}_s^{t,x,n,i}-{Z}_s^{t,x,n,i-1},\\
&&\bar{f}^{i}(s,x)=f(s,X_{s}^{t,x},{Y}_s^{t,x,n,i},{Z}_s^{t,x,n,i})-f(s,X_{s}^{t,x},{Y}_s^{t,x,n,i-1},{Z}_s^{t,x,n,i-1}),\\
&&\bar{g}^{i}_j(s,x)=g_j(s,X_{s}^{t,x},{Y}_s^{t,x,n,i},{Z}_s^{t,x,n,i})-g_j(s,X_{s}^{t,x},{Y}_s^{t,x,n,i-1},{Z}_s^{t,x,n,i-1}),\
i=1,2,\cdot\cdot\cdot,\ t\leq s\leq T.
\end{eqnarray*}
Then, for a.e. $x\in{\mathbb{R}^{d}}$,
$(\bar{Y}_s^{t,x,n,N},\bar{Z}_s^{t,x,n,N})$ satisfies
\begin{eqnarray*}
\bar{Y}_{s}^{t,x,n,N}&=&\int_{s}^{T}\bar{f}^{N-1}(r,x)dr-\sum_{j=1}^{n}\int_{s}^{T}\bar{g}^{N-1}_j(r,x)d^\dagger{\hat{\beta}}_j(r)-\int_{s}^{T}\langle\bar{Z}_{r}^{t,x,n,N},dW_r\rangle.
\end{eqnarray*}
Applying generalized It$\hat {\rm o}$'s formula to ${\rm
e}^{Kr}|\bar{Y}_{r}^{t,x,n,N}|^2$ for a.e. $x\in\mathbb{R}^{d}$, by
the Young inequality and Condition (H.2), we can deduce that
\begin{eqnarray}\label{zhangsss}
&&\int_{\mathbb{R}^{d}}{\rm e}^{Ks}|\bar{Y}_{s}^{t,x,n,N}|^2\rho^{-1}(x)dx+K\int_{s}^{T}\int_{\mathbb{R}^{d}}{\rm e}^{Kr}|\bar{Y}_{r}^{t,x,n,N}|^2\rho^{-1}(x)dxdr\nonumber\\
&&+\int_{s}^{T}\int_{\mathbb{R}^{d}}{\rm e}^{Kr}|\bar{Z}_{r}^{t,x,n,N}|^2\rho^{-1}(x)dxdr\nonumber\\
&\leq&\int_{s}^{T}\int_{\mathbb{R}^{d}}{\rm e}^{Kr}\big(2C|\bar{Y}_{r}^{t,x,n,N}|^2+{1\over2}|\bar{Y}_{r}^{t,x,n,N-1}|^2+{1\over2}|\bar{Z}_{r}^{t,x,n,N-1}|^2\big)\rho^{-1}(x)dxdr\nonumber\\
&&+\int_{s}^{T}\int_{\mathbb{R}^{d}}{\rm e}^{Kr}(\sum_{j=1}^{\infty}C_j|\bar{Y}_{r}^{t,x,n,N-1}|^2+\sum_{j=1}^{\infty}\alpha_j|\bar{Z}_{r}^{t,x,n,N-1}|^2\big)\rho^{-1}(x)dxdr\nonumber\\
&&-\sum_{j=1}^{n}\int_{s}^{T}\int_{\mathbb{R}^{d}}{\rm e}^{Kr}2\bar{Y}_{r}^{t,x,n,N}\bar{g}^{N-1}_j(r,x)\rho^{-1}(x)dxd^\dagger{\hat{\beta}}_j(r)\nonumber\\
&&-\int_{s}^{T}\langle\int_{\mathbb{R}^{d}}{\rm
e}^{Kr}2\bar{Y}_{r}^{t,x,n,N}\bar{Z}_{r}^{t,x,n,N}\rho^{-1}(x)dx,dW_r\rangle.
\end{eqnarray}
Then we have
\begin{eqnarray*}
&&(K-2C)E[\int_{s}^{T}\int_{\mathbb{R}^{d}}{\rm e}^{Kr}|{\bar{Y}}_{r}^{t,x,n,N}|^2\rho^{-1}(x)dxdr]+E[\int_{s}^{T}\int_{\mathbb{R}^{d}}{\rm e}^{Kr}|{\bar{Z}}_{r}^{t,x,n,N}|^2\rho^{-1}(x)dxdr]\\
&\leq&({1\over2}+\sum_{j=1}^{\infty}\alpha_j)E[\int_{s}^{T}\int_{\mathbb{R}^{d}}{\rm
e}^{Kr}\big((1+2\sum_{j=1}^{\infty}C_j)|{\bar{Y}}_{r}^{t,x,n,N-1}|^2+|{\bar{Z}}_{r}^{t,x,n,N-1}|^2\big)\rho^{-1}(x)dxdr].
\end{eqnarray*}
Letting $K=1+2C+2\sum_{j=1}^{\infty}C_j$, we have
\begin{eqnarray}\label{zhangttt}
&&E[\int_{s}^{T}\int_{\mathbb{R}^{d}}{\rm e}^{Kr}\big((1+2\sum_{j=1}^{\infty}C_j)|{\bar{Y}}_{r}^{t,x,n,N}|^2+|{\bar{Z}}_{r}^{t,x,n,N}|^2\big)\rho^{-1}(x)dxdr]\\
&\leq&({1\over2}+\sum_{j=1}^{\infty}\alpha_j)E[\int_{s}^{T}\int_{\mathbb{R}^{d}}{\rm
e}^{Kr}\big((1+2\sum_{j=1}^{\infty}C_j)|{\bar{Y}}_{r}^{t,x,n,N-1}|^2+|{\bar{Z}}_{r}^{t,x,n,N-1}|^2\big)\rho^{-1}(x)dxdr].\nonumber
\end{eqnarray}
Note that $E[\int_{t}^{T}\int_{\mathbb{R}^{d}}{\rm
e}^{Kr}\big((1+2\sum_{j=1}^{\infty}C_j)|\cdot|^2+|\cdot|^2\big)\rho^{-1}(x)dxdr]$
is equivalent to
$E[\int_{t}^{T}\int_{\mathbb{R}^{d}}\big(|\cdot|^2+|\cdot|^2\big)\rho^{-1}(x)dxdr]$.
From the contraction principle, the mapping (\ref{jia3}) has a pair
of fixed point $(Y_{\cdot}^{t,\cdot,n},Z_{\cdot}^{t,\cdot,n})$ that
is the limit of the Cauchy sequence
${\{(Y_{\cdot}^{t,\cdot,n,N},Z_{\cdot}^{t,\cdot,n,N})\}}_{N=1}^{\infty}$
in
$M^{2,0}([t,T];L_{\rho}^2({\mathbb{R}^{d}};{\mathbb{R}^{1}}))\bigotimes
M^{2,0}([t,T];L_{\rho}^2({\mathbb{R}^{d}};{\mathbb{R}^{d}}))$. We
then prove $Y_{\cdot}^{t,\cdot,n}$ is also the limit of
$Y_{\cdot}^{t,\cdot,n,N}$ in
$S^{2,0}([t,T];L_{\rho}^2({\mathbb{R}^{d}};{\mathbb{R}^{1}}))$ as
$N\rightarrow\infty$. For this, we only need to prove
${\{Y_{\cdot}^{t,\cdot,n,N}\}}_{N=1}^{\infty}$ is a Cauchy sequence
in $S^{2,0}([t,T];L_{\rho}^2({\mathbb{R}^{d}};{\mathbb{R}^{1}}))$.
Similar as in (\ref{zhangzzz00}), by the B-D-G inequality and
Cauchy-Schwartz inequality, from (\ref{zhangsss}), we have
\begin{eqnarray}\label{zhanguuu}
&&E[\sup_{t\leq s\leq T}\int_{\mathbb{R}^{d}}{\rm e}^{Ks}|\bar{Y}_{s}^{t,x,n,N}|^2\rho^{-1}(x)dx]\\
&\leq&M_3E[\int_{s}^{T}\int_{\mathbb{R}^{d}}{\rm
e}^{Kr}\big(|\bar{Y}_{r}^{t,x,n,N-1}|^2+|\bar{Z}_{r}^{t,x,n,N-1}|^2+|\bar{Y}_{r}^{t,x,n,N}|^2+|\bar{Z}_{r}^{t,x,n,N}|^2\big)\rho^{-1}(x)dxdr],\nonumber
\end{eqnarray}
where $M_3>0$ is independent of $n$ and $N$. Without losing any
generality, assume that $M\geq N$. We can deduce from
(\ref{zhangttt}) and (\ref{zhanguuu}) that
\begin{eqnarray*}
&&\big(E[\sup_{t\leq s\leq T}\int_{\mathbb{R}^{d}}{{|{Y}_{s}^{t,x,n,M}-{Y}_{s}^{t,x,n,N}|}^2}\rho^{-1}(x)dx]\big)^{1\over2}\nonumber\\
&\leq&\sum_{i=N+1}^{M}\big(E[\sup_{t\leq s\leq T}\int_{\mathbb{R}^{d}}{{|\bar{Y}_{s}^{t,x,n,i}|}^2}\rho^{-1}(x)dx]\big)^{1\over2}\nonumber\\
&\leq&\sum_{i=N+1}^{M}\big(M_3E[\int_{t}^{T}\int_{\mathbb{R}^{d}}{\rm e}^{Kr}\big(|\bar{Y}_{r}^{t,x,n,i-1}|^2+|\bar{Z}_{r}^{t,x,n,i-1}|^2+|\bar{Y}_{r}^{t,x,n,i}|^2+|\bar{Z}_{r}^{t,x,n,i}|^2\big)\rho^{-1}(x)dxdr]\big)^{1\over2}\nonumber\\
&\leq&\sum_{i=N+1}^{M}\big(2M_3E[\int_{t}^{T}\int_{\mathbb{R}^{d}}{\rm e}^{Kr}\big((1+2\sum_{j=1}^{\infty}C_j)|\bar{Y}_{r}^{t,x,n,i-1}|^2+|\bar{Z}_{r}^{t,x,n,i-1}|^2\big)\rho^{-1}(x)dxdr]\big)^{1\over2}\nonumber\\
&\leq&\sum_{i=N+1}^{\infty}({{1\over2}+\sum_{j=1}^{\infty}\alpha_j})^{{i-2}\over2}\big(2M_3E[\int_{t}^{T}\int_{\mathbb{R}^{d}}{\rm
e}^{Kr}\big((1+2\sum_{j=1}^{\infty}C_j)|{Y}_{r}^{t,x,n,1}|^2+|{Z}_{r}^{t,x,n,1}|^2\big)\rho^{-1}(x)dxdr]\big)^{1\over2}\nonumber\\
&\longrightarrow&0\ {\rm as}\ M,\ N\longrightarrow\infty.
\end{eqnarray*}
The lemma is proved. $\hfill\diamond$\\

Following a similar procedure as in the proof of Lemma \ref{qi1a00},
and using It$\hat {\rm o}$'s formula to ${\rm
e}^{Kr}{|{Y}_{r}^{t,x,n}|}^2$, by the B-D-G inequality, we have the
following estimation for the solution of Eq.(\ref{zhang66100}):
\begin{prop}\label{qi053}
Under the conditions of Theorem \ref{qi052},
$({Y}_\cdot^{t,\cdot,n},Z^{t,\cdot,n}_{\cdot})$ satifies
\begin{eqnarray*}
\sup_nE[\sup_{t\leq s\leq
T}\int_{\mathbb{R}^d}|Y_s^{t,x,n}|^2\rho^{-1}(x)dx]+\sup_nE[\int_{t}^{T}\int_{\mathbb{R}^{d}}|Z_r^{t,x,n}|^2\rho^{-1}(x)dxdr]<\infty.
\end{eqnarray*}
\end{prop}

\begin{rmk}\label{qi315} For $s\in[0,t]$, Eq.(\ref{zhang66100}) is equivalent to the
following BDSDE
\begin{eqnarray}\label{zhang663}
Y_s^{x,n}=&&Y_t^{t,x,n}+\int_{s}^{t}f(r,x,Y^{x,n}_r,Z^{x,n}_r)dr\nonumber\\
&&-\sum_{j=1}^{n}\int_{s}^{t}g_j(r,x,Y^{x,n}_r,Z^{x,n}_r)d^\dagger{\hat{\beta}}_j(r)-\int_{s}^{t}\langle
Z^{x,n}_r,dW_r\rangle.
\end{eqnarray}
Note that $Y_t^{t,x,n}$ satisfies Condition (H.1). By a similar
method as in the proof of Theorem \ref{jia} and Proposition
\ref{qi053}, we can obtain a
$({Y}_\cdot^{\cdot,n},{Z}_\cdot^{\cdot,n})\in
S^{2,0}([0,t];L_{\rho}^2({\mathbb{R}^{d}};{\mathbb{R}^{1}}))\bigotimes
M^{2,0}([0,t];L_{\rho}^2({\mathbb{R}^{d}};{\mathbb{R}^{d}}))$, is
the unique solution of Eq.(\ref{zhang663}). Moreover,
\begin{eqnarray*}
\sup_nE[\sup_{0\leq s\leq
t}\int_{\mathbb{R}^d}|Y_s^{x,n}|^2\rho^{-1}(x)dx]+\sup_nE[\int_{0}^{t}\int_{\mathbb{R}^{d}}|Z_r^{x,n}|^2\rho^{-1}(x)dxdr]<\infty.
\end{eqnarray*}
To unify the notation, we define
$({Y}_s^{t,x,n},{Z}_s^{t,x,n})=({Y}_s^{x,n},{Z}_s^{x,n})$ when
$s\in[0,t)$. Then $({Y}_\cdot^{t,\cdot,n},{Z}_\cdot^{t,\cdot,n})\in
S^{2,0}([0,T];L_{\rho}^2({\mathbb{R}^{d}};{\mathbb{R}^{1}}))\bigotimes
M^{2,0}([0,T];L_{\rho}^2({\mathbb{R}^{d}};{\mathbb{R}^{d}}))$.
Furthermore, we have
\begin{eqnarray}\label{zhang668}
\sup_nE[\sup_{0\leq s\leq
T}\int_{\mathbb{R}^d}|Y_s^{t,x,n}|^2\rho^{-1}(x)dx]+\sup_nE[\int_{0}^{T}\int_{\mathbb{R}^{d}}|Z_r^{t,x,n}|^2\rho^{-1}(x)dxdr]<\infty.
\end{eqnarray}
\end{rmk}
{\em Proof of Theorem \ref{qi052}}. The proof of the uniqueness is
rather similar to the uniqueness proof in Theorem \ref{jia}.

\underline{Existence}. By Theorem \ref{jia} and Remark \ref{qi315},
for each n, there exists a unique solution
$({Y}_\cdot^{t,\cdot,n},{Z}_\cdot^{t,\cdot,n})\in
S^{2,0}([0,T];L_{\rho}^2({\mathbb{R}^{d}};{\mathbb{R}^{1}}))\bigotimes
M^{2,0}([0,T];L_{\rho}^2({\mathbb{R}^{d}};{\mathbb{R}^{d}}))$ to
Eq.(\ref{zhang66100}). We will prove
$({Y}_\cdot^{t,\cdot,n},{Z}_\cdot^{t,\cdot,n})$ is a Cauchy sequence
in
$S^{2,0}([0,T];L_{\rho}^2({\mathbb{R}^{d}};{\mathbb{R}^{1}}))\bigotimes
M^{2,0}([0,T];L_{\rho}^2({\mathbb{R}^{d}};{\mathbb{R}^{d}}))$.
Without losing any generality, assume that $m\geq n$, and define
\begin{eqnarray*}
&&\bar{Y}_{s}^{t,x,m,n}={Y}_{s}^{t,x,m}-{Y}_{s}^{t,x,n}, \ \
\bar{Z}_{s}^{t,x,m,n}={Z}_{s}^{t,x,m}-{Z}_{s}^{t,x,n},\\
&&\bar{f}^{m,n}(s,x)=f(s,X_{s}^{t,x},{Y}_{s}^{t,x,m},{Z}_{s}^{t,x,m})-f(s,X_{s}^{t,x},Y_{s}^{t,x,n},Z_{s}^{t,x,n}),\\
&&\bar{g}^{m,n}_j(s,x)=g_j(s,X_{s}^{t,x},{Y}_{s}^{t,x,m},{Z}_{s}^{t,x,m})-g_j(s,X_{s}^{t,x},Y_{s}^{t,x,n},Z_{s}^{t,x,n}),\
\ \ \ \ \ 0\leq s\leq T.
\end{eqnarray*}
Then for $0\leq s\leq T$ and a.e. $x\in\mathbb{R}^d$,
\begin{eqnarray*}
\left\{\begin{array}{l}
d\bar{Y}_{s}^{t,x,m,n}=-\bar{f}^{m,n}(s,x)ds+\sum_{j=1}^{n}\bar{g}^{m,n}_j(s,x)d^\dagger{\hat{\beta}}_j(s)\\
\ \ \ \ \ \ \ \ \ \ \ \ \ \ \ \ +\sum_{j=n+1}^{m}{g}_j(s,X_{s}^{t,x},{Y}_{s}^{t,x,m},{Z}_{s}^{t,x,m})d^\dagger{\hat{\beta}}_j(s)+\langle\bar{Z}_{s}^{t,x,m,n},dW_s\rangle\\
\bar{Y}_{T}^{t,x,m,n}=0\ \ \ a.s..
\end{array}\right.
\end{eqnarray*}
Applying It$\hat {\rm o}$'s formula to ${\rm
e}^{Kr}{{|\bar{Y}_r^{t,x,m,n}|}^2}$ for a.e. $x\in\mathbb{R}^{d}$,
we have
\begin{eqnarray}\label{zhang667}
&&\int_{\mathbb{R}^{d}}{\rm e}^{Ks}{|\bar{Y}_s^{t,x,m,n}|}^2\rho^{-1}(x)dx+(K-2C-\sum_{j=1}^{\infty}{C_j}-{1\over2})\int_{s}^{T}\int_{\mathbb{R}^{d}}{\rm e}^{Kr}{{|\bar{Y}_r^{t,x,m,n}|}^2}\rho^{-1}(x)dxdr\nonumber\\
&&+({1\over2}-\sum_{j=1}^{\infty}\alpha_j)\int_{s}^{T}\int_{\mathbb{R}^{d}}{\rm e}^{Kr}|\bar{Z}_r^{t,x,m,n}|^2\rho^{-1}(x)dxdr\nonumber\\
&\leq&C_p\sum_{j=n+1}^{m}\{(C_j+\alpha_j)\big(\int_{s}^{T}\int_{\mathbb{R}^{d}}(|{Y}_{r}^{t,x,m}|^2+|{Z}_{r}^{t,x,m}|^2)\rho^{-1}(x)dxdr\nonumber\\
&&+\int_{s}^{T}\int_{\mathbb{R}^{d}}|{g}_j(r,X_{r}^{t,x},0,0)|^2\rho^{-1}(x)dxdr\big)\}-\sum_{j=1}^{n}\int_{s}^{T}\int_{\mathbb{R}^{d}}2{\rm e}^{Kr}\bar{Y}_r^{t,x,m,n}\bar{g}^{m,n}_j(r,x)\rho^{-1}(x)dxd^\dagger{\hat{\beta}}_j(r)\nonumber\\
&&-\sum_{j=n+1}^{m}\int_{s}^{T}\int_{\mathbb{R}^{d}}2{\rm e}^{Kr}\bar{Y}_r^{t,x,m,n}{g}_j(r,X_{r}^{t,x},{Y}_{r}^{t,x,m},{Z}_{r}^{t,x,m})\rho^{-1}(x)dxd^\dagger{\hat{\beta}}_j(r)\nonumber\\
&&-\int_{s}^{T}\langle\int_{\mathbb{R}^{d}}2{\rm
e}^{Kr}\bar{Y}_r^{t,x,m,n}\bar{Z}_r^{t,x,m,n}\rho^{-1}(x)dx,dW_r\rangle.
\end{eqnarray}
All the terms on the left hand side of (\ref{zhang667}) are positive
when taking $K$ sufficiently large. Take expectation on both sides
of (\ref{zhang667}), then by Lemma \ref{qi045} and (\ref{zhang668}),
we have
\begin{eqnarray}\label{zhang669}
&&E[\int_{0}^{T}\int_{\mathbb{R}^{d}}{\rm e}^{Kr}{|\bar{Y}_r^{t,x,m,n}|}^2\rho^{-1}(x)dxdr]+E[\int_{0}^{T}\int_{\mathbb{R}^{d}}{\rm e}^{Kr}|\bar{Z}_r^{t,x,m,n}|^2\rho^{-1}(x)dxdr]\nonumber\\
&\leq&C_p\sum_{j=n+1}^{m}\{(C_j+\alpha_j)\big(\sup_nE[\int_{0}^{T}\int_{\mathbb{R}^{d}}(|{Y}_{r}^{t,x,n}|^2+|{Z}_{r}^{t,x,n}|^2)\rho^{-1}(x)dxdr]\nonumber\\
&&+\int_{0}^{T}\int_{\mathbb{R}^{d}}|{g}_j(r,x,0,0)|^2\rho^{-1}(x)dxdr\big)\}\longrightarrow0,\
\ \ \ \ {\rm as}\ n,\ m\longrightarrow\infty.
\end{eqnarray}
Also by the B-D-G inequality, from (\ref{zhang667}) we have
\begin{eqnarray*}
&&E[\sup_{0\leq s\leq T}\int_{\mathbb{R}^{d}}{\rm e}^{Ks}{|\bar{Y}_s^{t,x,m,n}|}^2\rho^{-1}(x)dx]\\
&\leq&C_pE[\int_{0}^{T}\int_{\mathbb{R}^{d}}{\rm e}^{Kr}(|\bar{Y}_r^{t,x,m,n}|^2+|\bar{Z}_r^{t,x,m,n}|^2)\rho^{-1}(x)dxdr]\\
&&+C_p\sum_{j=n+1}^{m}(C_j+\alpha_j)\big(\sup_nE[\int_{0}^{T}\int_{\mathbb{R}^{d}}(|{Y}_{r}^{t,x,n}|^2+|{Z}_{r}^{t,x,n}|^2)\rho^{-1}(x)dxdr]\big)\\
&&+C_p\sum_{j=n+1}^{m}\int_{0}^{T}\int_{\mathbb{R}^{d}}|{g}_j(r,x,0,0)|^2\rho^{-1}(x)dxdr.
\end{eqnarray*}
So by (\ref{zhang668}), (\ref{zhang669}) and Condition (H.3), we
have
\begin{eqnarray*}
E[\sup_{0\leq s\leq T}\int_{\mathbb{R}^{d}}{\rm
e}^{Ks}{{|\bar{Y}_s^{t,x,m,n}|}^2}\rho^{-1}(x)dx]\longrightarrow0,\
\ \ \ \ {\rm as}\ n,\ m\longrightarrow\infty.
\end{eqnarray*}
Therefore $({Y}_\cdot^{t,\cdot,n},{Z}_\cdot^{t,\cdot,n})$ is a
Cauchy sequence in
$S^{2,0}([0,T];L_{\rho}^2({\mathbb{R}^{d}};{\mathbb{R}^{1}}))\bigotimes
M^{2,0}([0,T];L_{\rho}^2({\mathbb{R}^{d}};{\mathbb{R}^{d}}))$ with
its limit denoted by $({Y}_s^{t,x},{Z}_s^{t,x})$. We will show that
$({Y}_\cdot^{t,\cdot},{Z}_\cdot^{t,\cdot})$ is the solution of
Eq.(\ref{qi20}), i.e. $({Y}_\cdot^{t,\cdot},{Z}_\cdot^{t,\cdot})$
satisfies (\ref{qi22}) for an arbitrary $\varphi\in
C_c^{0}(\mathbb{R}^d;\mathbb{R}^1)$. For this, we will prove that
Eq.(\ref{zhang662}) converges to Eq.(\ref{qi22}) in $L^2(\Omega)$
term by term as $n\longrightarrow\infty$. Here we only show the
convergence of the third term,
\begin{eqnarray*}
&&E[\ |\sum_{j=1}^{n}\int_{s}^{T}\int_{\mathbb{R}^{d}}g_j(r,X_{r}^{t,x},Y_r^{t,x,n},Z_r^{t,x,n})\varphi(x)dxd^\dagger{\hat{\beta}}_j(r)\\
&&\ \ \ \ \ \ -\sum_{j=1}^{\infty}\int_{s}^{T}\int_{\mathbb{R}^{d}}g_j(r,X_{r}^{t,x},Y_r^{t,x},Z_r^{t,x})\varphi(x)dxd^\dagger{\hat{\beta}}_j(r)|^2]\\
&\leq&2E[\ |\sum_{j=1}^{n}\int_{s}^{T}\int_{\mathbb{R}^{d}}\big(g_j(r,X_{r}^{t,x},Y_r^{t,x,n},Z_r^{t,x,n})-g_j(r,X_{r}^{t,x},Y_r^{t,x},Z_r^{t,x})\big)\varphi(x)dxd^\dagger{\hat{\beta}}_j(r)|^2]\\
&&+2E[\ |\sum_{j=n+1}^{\infty}\int_{s}^{T}\int_{\mathbb{R}^{d}}g_j(r,X_{r}^{t,x},Y_r^{t,x},Z_r^{t,x})\varphi(x)dxd^\dagger{\hat{\beta}}_j(r)|^2]\\
&\leq&C_p\sum_{j=1}^{\infty}({C_j+\alpha_j})E[\int_{s}^{T}\int_{\mathbb{R}^{d}}(|Y_r^{t,x,n}-Y_r^{t,x}|^2+|Z_r^{t,x,n}-Z_r^{t,x}|^2)\rho^{-1}(x)dxdr]\\
&&+C_pE[\ |\sum_{j=n+1}^{\infty}\int_{s}^{T}\int_{\mathbb{R}^{d}}\big(g_j(r,X_{r}^{t,x},Y_r^{t,x},Z_r^{t,x})-g_j(r,X_{r}^{t,x},0,0)\big)\varphi(x)dxd^\dagger{\hat{\beta}}_j(r)|^2]\\
&&+C_pE[\
|\sum_{j=n+1}^{\infty}\int_{s}^{T}\int_{\mathbb{R}^{d}}g_j(r,X_{r}^{t,x},0,0)\varphi(x)dxd^\dagger{\hat{\beta}}_j(r)|^2].
\end{eqnarray*}
Note
\begin{eqnarray}\label{qii6}
&&E[\ |\sum_{j=n+1}^{\infty}\int_{s}^{T}\int_{\mathbb{R}^{d}}\big(g_j(r,X_{r}^{t,x},Y_r^{t,x},Z_r^{t,x})-g_j(r,X_{r}^{t,x},0,0)\big)\varphi(x)dxd^\dagger{\hat{\beta}}_j(r)|^2]\nonumber\\
&=&E[\int_{s}^{T}\|\int_{\mathbb{R}^{d}}\big(g(r,X_{r}^{t,x},Y_r^{t,x},Z_r^{t,x})-g(r,X_{r}^{t,x},0,0)\big)\varphi(x)dx(\sum_{j=n+1}^{\infty}{\lambda_j}e_j\otimes e_j)^{1\over2}\|_{L_U}^2dr]\nonumber\\
&=&E[\int_{s}^{T}\sum_{i=1}^{\infty}|\int_{\mathbb{R}^{d}}\big(g(r,X_{r}^{t,x},Y_r^{t,x},Z_r^{t,x})-g(r,X_{r}^{t,x},0,0)\big)\varphi(x)dx\sum_{j=n+1}^{\infty}\sqrt{{\lambda_j}}e_j\langle e_j,e_i\rangle|^2dr]\nonumber\\
&=&E[\sum_{j=n+1}^{\infty}\int_{s}^{T}|\int_{\mathbb{R}^{d}}\big(g_j(r,X_{r}^{t,x},Y_r^{t,x},Z_r^{t,x})-g_j(r,X_{r}^{t,x},0,0)\big)\varphi(x)dx|^2dr]\nonumber\\
&\leq&C_pE[\sum_{j=n+1}^{\infty}\int_{s}^{T}\int_{\mathbb{R}^{d}}|g_j(r,X_{r}^{t,x},Y_r^{t,x},Z_r^{t,x})-g_j(r,X_{r}^{t,x},0,0)|^2\rho^{-1}(x)dxdr]\nonumber\\
&\leq&C_p\sum_{j=n+1}^{\infty}({C_j+\alpha_j})E[\int_{s}^{T}\int_{\mathbb{R}^{d}}(|Y_r^{t,x}|^2+|Z_r^{t,x}|^2)\rho^{-1}(x)dxdr]\longrightarrow0.
\end{eqnarray}
Here we used $(\sum_{j=n+1}^{\infty}{\lambda_j}e_j\otimes
e_j)^{1\over2}=\sum_{j=n+1}^{\infty}\sqrt{{\lambda_j}}e_j\otimes
e_j$. This can be verified as follows: for an arbitrary $u\in U$, by
definition of tensor operator,
\begin{eqnarray*}
(\sum_{j=n+1}^{\infty}\sqrt{{\lambda_j}}e_j\otimes e_j)(\sum_{i=n+1}^{\infty}\sqrt{{\lambda_i}}e_i\otimes e_i)u&=&\sum_{j=n+1}^{\infty}\sqrt{{\lambda_j}}e_j\langle e_j,\sum_{i=n+1}^{\infty}\sqrt{{\lambda_i}}e_i\langle e_i,u\rangle\rangle\\
&=&\sum_{j=n+1}^{\infty}\sqrt{{\lambda_j}}e_j\langle\sqrt{{\lambda_j}}e_j,e_j\rangle\langle e_j,u\rangle\\
&=&(\sum_{j=n+1}^{\infty}{\lambda_j}e_j\otimes e_j)u.
\end{eqnarray*}
Similarly we have
\begin{eqnarray}\label{qib}
&&C_pE[\ |\sum_{j=n+1}^{\infty}\int_{s}^{T}\int_{\mathbb{R}^{d}}g_j(r,X_{r}^{t,x},0,0)\varphi(x)dxd^\dagger{\hat{\beta}}_j(r)|^2]\nonumber\\
&\leq&C_p\int_{s}^{T}\int_{\mathbb{R}^{d}}\sum_{j=n+1}^{\infty}|g_j(r,x,0,0)|^2\rho^{-1}(x)dxdr\longrightarrow0.
\end{eqnarray}
That is to say ${(Y_s^{t,x}, Z_s^{t,x})}_{0\leq s\leq T}$ satisfies
Eq.(\ref{qi22}). The proof of Theorem \ref{qi052} is completed.
$\hfill\diamond$

\section{Weak solutions of the corresponding SPDEs}
\setcounter{equation}{0}

In section 3, we proved the existence and uniqueness of the weak
solution of BDSDE (\ref{qi20}). We obtained the solution
$({Y}_s^{t,x},{Z}_s^{t,x})$ by taking the limit of
$({Y}_s^{t,x,n},{Z}_s^{t,x,n})$ of the solutions of
Eq.(\ref{zhang66100}) in the space
$S^{2,0}([0,T];L_{\rho}^2({\mathbb{R}^{d}};{\mathbb{R}^{1}}))\bigotimes
M^{2,0}([0,T];L_{\rho}^2({\mathbb{R}^{d}};{\mathbb{R}^{d}}))$. We
still start from Eq.(\ref{zhang66100}) in this section. A direct
application of Proposition \ref{000000} and Fubini theorem
immediately leads to
\begin{prop}\label{0000} Under Conditions {\rm(H.1)}--{\rm(H.4)}, if we define $u^n(t,x)=Y_t^{t,x,n}$, $v^n(t,x)=Z_t^{t,x,n}$, then $u^n(s,X^{t,x}_s)=Y_s^{t,x,n}$,
$v^n(s,X^{t,x}_s)=Z_s^{t,x,n}$ for a.e. $s\in[t,T]$,
$x\in{\mathbb{R}^{d}}$ a.s..
\end{prop}

We first use the idea of Bally and Matoussi \cite{ba-ma} to give the
correspondence between the weak solutions of SPDEs and BDSDEs with finite
dimensional noise. Consider the BDSDEs (\ref{zhang66111}). Define
the mollifier $K^m(x)=mc\exp\{{1\over{(mx-1)^2-1}}\}$, if
$0<x<{2\over m}$; $K^m(x)=0$ otherwise, where $c$ is chosen such
that $\int_{-\infty}^{+\infty}K^m(x)dx=1$. Define
$h^m(x)=\int_{\mathbb{R}^{d}}h(y)K^m(x-y)dy$,
$\tilde{f}^m(r,x)=\int_{\mathbb{R}^{d}}\tilde{f}(r,y)K^m(x-y)dy$ and
$\tilde{g}_j^m(r,x)=\int_{\mathbb{R}^{d}}\tilde{g}_j(r,y)K^m(x-y)dy$.
It is easy to see from standard results in analysis that
$h^m(\cdot)\rightarrow h(\cdot)$,
$\tilde{f}^m(r,\cdot)\rightarrow\tilde{f}(r,\cdot)$ and
$\tilde{g}_j^m(r,\cdot)\rightarrow\tilde{g}_j(r,\cdot)$ in
$L_{\rho}^2({\mathbb{R}^{d}};{\mathbb{R}^{1}})$ respectively. Denote
by $(\tilde{Y}_{s,m}^{t,x,n},\tilde{Z}_{s,m}^{t,x,n})$ the solution
of the following BDSDEs:
\begin{eqnarray*}\label{zhang66300}
\tilde{Y}_{s,m}^{t,x,n}=\tilde{h}^m(X_T^{t,x})+\int_{s}^{T}\tilde{f}^m(r,X_r^{t,x})dr-\sum_{j=1}^{n}\int_{s}^{T}\tilde{g}_j^m(r,X_r^{t,x})d^\dagger{\hat{\beta}}_j(r)-\int_{s}^{T}\langle
\tilde{Z}_{r,m}^{t,x,n},dW_r\rangle.
\end{eqnarray*}
Let $u^m(t,x)=Y^{t,x,n}_{t,m}$. Then following classical results of
Pardoux and Peng \cite{pa-pe3}, we have
$\tilde{Z}_{t,m}^{t,x,n}=\sigma^*\nabla \tilde{u}^n_m(t,x)$, and
$\tilde{Y}_{s,m}^{t,x,n}=\tilde{u}^n_m(s,X_s^{t,x})=\tilde{Y}_{s,m}^{s,X_s^{t,x},n}$,
$\tilde{Z}_{s,m}^{t,x,n}=\sigma^*\nabla
\tilde{u}^n_m(s,X_s^{t,x})=\tilde{Z}_{s,m}^{s,X_s^{t,x},n}$.
Moreover $\tilde{u}^n_m(t,x)$ satisfies the smootherized SPDE. In
particular, for any smooth test function $\Psi\in
C_c^{1,\infty}([0,T]\times\mathbb{R}^d;\mathbb{R}^1)$, we still have
\begin{eqnarray}\label{zhang66301}
&&\int_{t}^{T}\int_{\mathbb{R}^{d}}\tilde{u}^n_m(s,x)\partial_s\Psi(s,x)dxds+\int_{\mathbb{R}^{d}}\tilde{u}^n_m(t,x)\Psi(t,x)dx-\int_{\mathbb{R}^{d}}\tilde{h}^m(x)\Psi(T,x)dx\nonumber\\
&&-{1\over2}\int_{t}^{T}\int_{\mathbb{R}^{d}}(\sigma^*\nabla \tilde{u}^n_m)(s,x)(\sigma^*\nabla\Psi)(s,x)dxds-\int_{t}^{T}\int_{\mathbb{R}^{d}}\tilde{u}^n_m(s,x)\nabla\big((b-\tilde{A})\Psi\big)(s,x)dxds\nonumber\\
&=&\int_{t}^{T}\int_{\mathbb{R}^{d}}\tilde{f}^m(s,x)\Psi(s,x)dxds-\sum_{j=1}^{n}\int_{t}^{T}\int_{\mathbb{R}^{d}}\tilde{g}_j^m(s,x)\Psi(s,x)dxd^\dagger{\hat{\beta}}_j(s)\
\ \ P-{\rm a.s.}.
\end{eqnarray}
But by standard estimates
\begin{eqnarray*}\label{zhang66302}
E[\int_{t}^{T}\int_{\mathbb{R}^{d}}(|\tilde{Y}_{s,m}^{t,x,n}-\tilde{Y}_s^{t,x,n}|^2+|\tilde{Z}_{s,m}^{t,x,n}-\tilde{Z}_s^{t,x,n}|^2)\rho^{-1}(x)dxds]\longrightarrow0\
\ \ {\rm as}\ m\rightarrow\infty.
\end{eqnarray*}
And as $m_1$, $m_2\rightarrow\infty$
\begin{eqnarray}\label{zhang66303}
&&E[\int_{t}^{T}\int_{\mathbb{R}^{d}}(|\tilde{u}^n_{m_1}(s,X_s^{t,x})-\tilde{u}^n_{m_2}(s,X_s^{t,x})|^2+|\sigma^*\nabla
\tilde{u}^n_{m_1}(s,X_s^{t,x})-\sigma^*\nabla
\tilde{u}^n_{m_2}(s,X_s^{t,x})|^2)\rho^{-1}(x)dxds]\nonumber\\
&=&E[\int_{t}^{T}\int_{\mathbb{R}^{d}}(|\tilde{Y}_{s,m_1}^{t,x,n}-\tilde{Y}_{s,m_2}^{t,x,n}|^2+|\tilde{Z}_{s,m_1}^{t,x,n}-\tilde{Z}_{s,m_2}^{t,x,n}|^2)\rho^{-1}(x)dxds]\longrightarrow0.
\end{eqnarray}
We define $\mathcal{H}$ to be the set of random
fields $\{w(s,x);\ s\in[0,T],\ x\in\mathbb{R}^{d}\}$ such that
$(w,\sigma^*\nabla w)\in
M^{2,0}([0,T];L_{\rho}^2({\mathbb{R}^{d}};{\mathbb{R}^{1}}))\bigotimes
M^{2,0}([0,T];L_{\rho}^2({\mathbb{R}^{d}};{\mathbb{R}^{d}}))$ with
the norm
$(E[\int_{0}^{T}\int_{\mathbb{R}^{d}}(|w(s,x)|^2
\linebreak +|(\sigma^*\nabla)w(s,x)|^2)\rho^{-1}(x)dxds)^{1\over2}$.
Following a standard argument as in the proof of the completeness of the Sobolev
spaces, we can prove $\mathcal{H}$ is complete. Now by the generalized
equivalence of norm principle
and (\ref{zhang66303}),
we can see that $\tilde{u}^n_m$ is a Cauchy sequence in
$\mathcal{H}$. So there exists $\tilde{u}^n\in\mathcal{H}$ such that
$(\tilde{u}^n_m,\sigma^*\nabla
\tilde{u}^n_m)\rightarrow(\tilde{u}^n,\sigma^*\nabla \tilde{u}^n)$
in
$M^{2,0}([0,T];L_{\rho}^2({\mathbb{R}^{d}};{\mathbb{R}^{1}}))\bigotimes
M^{2,0}([0,T];L_{\rho}^2({\mathbb{R}^{d}};{\mathbb{R}^{d}}))$.}
Moreover $\tilde{Y}_s^{t,x,n}=u^n(s,X_s^{t,x})$,
$\tilde{Z}_s^{t,x,n}=\sigma^*\nabla u^n(s,X_s^{t,x})$ for a.e.
$s\in[t,T]$, $x\in\mathbb{R}^{d}$ a.s.. Now it is easy to pass the
limit as $m\rightarrow\infty$ in (\ref{zhang66301}) to conclude that
$\tilde{u}^n$ is a weak solution of the corresponding SPDEs. For the
nonlinear case, we can regard
$\tilde{f}(r,x)=f(r,x,\tilde{u}^n(r,x),\sigma^*\nabla
\tilde{u}^n(r,x))$,
$\tilde{g}_j(r,x)=g_j(r,x,\tilde{u}^n(r,x),\sigma^*\nabla
\tilde{u}^n(r,x))$, and $\tilde{f}$, $\tilde{g}_j$ satisfy the
conditions in the above argument. Using a similar proof as in the
proof of Theorem 3.1 in \cite{ba-ma} together with Theorem \ref{jia}
and Proposition \ref{0000}, we have, under Conditions (H.1)--(H.4),
$v^n(t,x)=(\sigma^*\nabla u^n)(t,x)$. Moreover, $(u^n,\sigma^*\nabla
u^n)\in
M^{2,0}([0,T];L_{\rho}^2({\mathbb{R}^{d}};{\mathbb{R}^{1}}))\bigotimes
M^{2,0}([0,T];L_{\rho}^2({\mathbb{R}^{d}};{\mathbb{R}^{d}}))$,
$u^n(t,x)$ is the weak solution of the following SPDE:
\begin{eqnarray*}
u^n(t,x)&=&h(x)+\int_{t}^{T}[\mathscr{L}u^n(s,x)+f\big(s,x,u^n(s,x),(\sigma^*\nabla u^n)(s,x)\big)]ds\\
&&-\sum_{j=1}^{n}\int_{t}^{T}g_j\big(s,x,u^n(s,x),(\sigma^*\nabla
u^n)(s,x)\big)d^\dagger{\hat{\beta}}_j(s),\ \ \ \ \ 0\leq t\leq
s\leq T.
\end{eqnarray*}
That is to say, for any $\Psi\in
C_c^{1,\infty}([0,T]\times\mathbb{R}^d;\mathbb{R}^1)$, we have
\begin{eqnarray}\label{qi25}
&&\int_{t}^{T}\int_{\mathbb{R}^{d}}u^n(s,x)\partial_s\Psi(s,x)dxds+\int_{\mathbb{R}^{d}}u^n(t,x)\Psi(t,x)dx-\int_{\mathbb{R}^{d}}h(x)\Psi(T,x)dx\nonumber\\
&&-{1\over2}\int_{t}^{T}\int_{\mathbb{R}^{d}}(\sigma^*\nabla u^n)(s,x)(\sigma^*\nabla\Psi)(s,x)dxds-\int_{t}^{T}\int_{\mathbb{R}^{d}}u^n(s,x)\nabla\big((b-\tilde{A})\Psi\big)(s,x)dxds\nonumber\\
&=&\int_{t}^{T}\int_{\mathbb{R}^{d}}f\big(s,x,u^n(s,x),(\sigma^*\nabla u^n)(s,x)\big)\Psi(s,x)dxds\nonumber\\
&&-\sum_{j=1}^{n}\int_{t}^{T}\int_{\mathbb{R}^{d}}g_j\big(s,x,u^n(s,x),(\sigma^*\nabla
u^n)(s,x)\big)\Psi(s,x)dxd^\dagger{\hat{\beta}}_j(s)\ \ \ P-{\rm
a.s.}.
\end{eqnarray}
By intuition if we define $u(t,x)={Y}_t^{t,x}$, it should be a "weak
solution" of the Eq.(\ref{zhang685}) with $u(T,x)=h(x)$. We will
prove this result.

First we need some necessary preparations.
\begin{prop}\label{qi061}
Under Conditions {\rm(H.1)}--{\rm(H.4)}, let $(Y_s^{t,x},Z_s^{t,x})$
be the solution of Eq.(\ref{qi20}). If we define $u(t,x)=Y_t^{t,x}$,
then $\sigma^*\nabla u(t,x)$ exists for a.e. $t\in[0,T]$,
$x\in\mathbb{R}^{d}$ a.s., and $u(s,X_s^{t,x})=Y_s^{t,x}$,
$(\sigma^*\nabla u)(s,X^{t,x}_s)=Z_s^{t,x}$ for a.e. $s\in[t,T]$,
$x\in\mathbb{R}^{d}$ a.s..
\end{prop}
{\em Proof}. First we prove $u^n$ is a Cauchy sequence in $\mathcal{H}$. For this, by Lemma \ref{qi045} and Proposition \ref{0000}, as $m$,
$n\rightarrow\infty$, we have
\begin{eqnarray*}
&&E[\int_{0}^{T}\int_{\mathbb{R}^{d}}(|u^m(s,x)-u^n(s,x)|^2+|(\sigma^*\nabla u^m)(s,x)-(\sigma^*\nabla u^n)(s,x)|^2)\rho^{-1}(x)dxds]\\
&\leq&C_pE[\int_{0}^{T}\int_{\mathbb{R}^{d}}(|u^m(s,X_s^{0,x})-u^n(s,X_s^{0,x})|^2+|(\sigma^*\nabla u^m)(s,X_s^{0,x})-(\sigma^*\nabla u^n)(s,X_s^{0,x})|^2)\rho^{-1}(x)dxds]\\
&=&C_pE[\int_{0}^{T}\int_{\mathbb{R}^{d}}(|Y_s^{0,x,m}-Y_s^{0,x,n}|^2+|Z_s^{0,x,m}-Z_s^{0,x,n}|^2)\rho^{-1}(x)dxds]\longrightarrow0.
\end{eqnarray*}
So there exists $\tilde{u}\in\mathcal{H}$ as the limit of $u^n$ such
that $\nabla\tilde{u}(s,x)$ exists for a.e. $s\in[0,T]$,
$x\in{\mathbb{R}^{d}}$ a.s. and
$E[\int_{0}^{T}\int_{\mathbb{R}^{d}}(|u^n(s,x)-\tilde{u}(s,x)|^2+|(\sigma^*\nabla
u^n)(s,x)-(\sigma^*\nabla\tilde{u})(s,x)|^2)\rho^{-1}(x)dxds]\longrightarrow0$.
We define $u(t,x)=Y_t^{t,x}$, then similar to the proof as in
Proposition \ref{0000}, by the uniqueness of solution of
Eq.(\ref{qi20}), we have $u(s,X_s^{t,x})=Y_s^{t,x}$ for a.e.
$s\in[t,T]$, $x\in{\mathbb{R}^{d}}$ a.s.. Since
\begin{eqnarray*}
&&E[\int_{0}^{T}\int_{\mathbb{R}^{d}}|u(s,x)-\tilde{u}(s,x)|^2\rho^{-1}(x)dxds]\\
&\leq&2E[\int_{0}^{T}\int_{\mathbb{R}^{d}}(|u(s,x)-u^n(s,x)|^2+|u^n(s,x)-\tilde{u}(s,x)|^2)\rho^{-1}(x)dxds]\\
&\leq&C_pE[\int_{0}^{T}\int_{\mathbb{R}^{d}}(|Y_s^{0,x}-Y_s^{0,x,n}|^2+|u^n(s,x)-\tilde{u}(s,x)|^2)\rho^{-1}(x)dxds]\longrightarrow0,
\end{eqnarray*}
$u(t,x)=\tilde{u}(t,x)$ for a.e. $t\in[0,T]$, $x\in\mathbb{R}^{d}$,
a.s.. So $\sigma^*\nabla u(t,x)$ exists for a.e. $t\in[0,T]$,
$x\in\mathbb{R}^{d}$, a.s.. Using Lemma \ref{qi045} again, we have
\begin{eqnarray*}
&&E[\int_{t}^{T}\int_{\mathbb{R}^{d}}(|u(s,X_s^{t,x})-Y_s^{t,x}|^2+|(\sigma^*\nabla u)(s,X_s^{t,x})-Z_s^{t,x}|^2)\rho^{-1}(x)dxds]\\
&\leq&2E[\int_{t}^{T}\int_{\mathbb{R}^{d}}(|u(s,X_s^{t,x})-u^n(s,X_s^{t,x})|^2+|u^n(s,X_s^{t,x})-Y_s^{t,x}|^2)\rho^{-1}(x)dxds]\\
&&+2E[\int_{t}^{T}\int_{\mathbb{R}^{d}}(|(\sigma^*\nabla u)(s,X_s^{t,x})-(\sigma^*\nabla u^n)(s,X_s^{t,x})|^2+|(\sigma^*\nabla u^n)(s,X_s^{t,x})-Z_s^{t,x}|^2)\rho^{-1}(x)dxds]\\
&\leq&C_pE[\int_{t}^{T}\int_{\mathbb{R}^{d}}(|u(s,x)-\tilde{u}(s,x)|^2+|\tilde{u}(s,x)-u^n(s,x)|^2+|Y_s^{t,x,n}-Y_s^{t,x}|^2)\rho^{-1}(x)dxds]\\
&&+C_pE[\int_{t}^{T}\int_{\mathbb{R}^{d}}(|(\sigma^*\nabla
u)(s,x)-(\sigma^*\nabla\tilde{u})(s,x)|^2+|(\sigma^*\nabla\tilde{u})(s,x)-(\sigma^*\nabla
u^n)(s,x)|^2\\
&&\ \ \ \ \ \ \ \ \ \ \ \ \ \ \ \ \ \ \ \
+|Z_s^{t,x,n}-Z_s^{t,x}|^2)\rho^{-1}(x)dxds]\longrightarrow0.
\end{eqnarray*}
So $u(s,X_s^{t,x})=Y_s^{t,x}$, $(\sigma^*\nabla
u)(s,X^{t,x}_s)=Z_s^{t,x}$ for a.e. $s\in[t,T]$,
$x\in{\mathbb{R}^d}$ a.s.. $\hfill\diamond$
\\

From Proposition \ref{qi061} and Lemma \ref{qi045}, it is easy to
know that
\begin{eqnarray*}
&&E[\int_{t}^{T}\int_{\mathbb{R}^d}|u^n(s,x)-u(s,x)|^2\rho^{-1}(x)dxds]\\
&&+E[\int_{t}^{T}\int_{\mathbb{R}^d}|(\sigma^*\nabla u^n)(s,x)-(\sigma^*\nabla u)(s,x)|^2\rho^{-1}(x)dxds]\\
&\leq&C_pE[\int_{t}^{T}\int_{\mathbb{R}^d}|u^n(s,X^{t,x}_s)-u(s,X^{t,x}_s)|^2\rho^{-1}(x)dxds]\\
&&+C_pE[\int_{t}^{T}\int_{\mathbb{R}^d}|(\sigma^*\nabla u^n)(s,X^{t,x}_s)-(\sigma^*\nabla u)(s,X^{t,x}_s)|^2\rho^{-1}(x)dxds]\\
&=&C_pE[\int_{t}^{T}\int_{\mathbb{R}^d}|{Y}_s^{t,x,n}-{Y}_s^{t,x}|^2\rho^{-1}(x)dxds]+C_pE[\int_{t}^{T}\int_{\mathbb{R}^d}|{Z}_s^{t,x,n}-{Z}_s^{t,x}|^2\rho^{-1}(x)dxds]\longrightarrow0,
\end{eqnarray*}
as $n\to \infty$. This will be used in the following theorem.
\begin{thm}\label{qi062} Under Conditions {\rm(H.1)}--{\rm(H.4)},
if we define $u(t,x)=Y_t^{t,x}$, where $(Y_s^{t,x},Z_s^{t,x})$ is
the solution of Eq.(\ref{qi20}), then $u(t,x)$ is the unique weak
solution of Eq.(\ref{zhang685}) with $u(T,x)=h(x)$. Moreover,
$u(s,X_s^{t,x})=Y_s^{t,x}$, $(\sigma^*\nabla
u)(s,X^{t,x}_s)=Z_s^{t,x}$ for a.e. $s\in[t,T]$,
$x\in\mathbb{R}^{d}$ a.s..
\end{thm}
{\em Proof}. From Proposition \ref{qi061}, we only need to verify
that this $u$ is the unique weak solution of Eq.(\ref{zhang685})
with $u(T,x)=h(x)$. By Lemma \ref{qi045}, it is easy to see that
$(\sigma^*\nabla u)(t,x)=Z_t^{t,x}$ for a.e. $t\in[0,T]$,
$x\in\mathbb{R}^d$, a.s.. Furthermore, by the generalized
equivalence of norm principle again we have
\begin{eqnarray*}
&&E[\int_{0}^{T}\int_{\mathbb{R}^d}(|u(s,x)|^2+|(\sigma^*\nabla u)(s,x)|^2)\rho^{-1}(x)dxds]\\
&\leq&C_pE[\int_{0}^{T}\int_{\mathbb{R}^d}(|u(s,X^{0,x}_s)|^2+|(\sigma^*\nabla u)(s,X^{0,x}_s)|^2)\rho^{-1}(x)dxds]\\
&=&C_pE[\int_{0}^{T}\int_{\mathbb{R}^d}(|Y_s^{0,x}|^2+|Z_s^{0,x}|^2)\rho^{-1}(x)dxds]<\infty.
\end{eqnarray*}
Now we verify that $u(t,x)$ satisfies (\ref{qi16}) with
$u(T,x)=h(x)$ by passing the limit in $L^2(\Omega)$ to (\ref{qi25}).
We only show the convergence of the last term. The last term
includes infinite dimensional integral, but
\begin{eqnarray*}
&&E[\ |\sum_{j=1}^{n}\int_{t}^{T}\int_{\mathbb{R}^{d}}g_j\big(s,x,u^n(s,x),(\sigma^*\nabla u^n)(s,x)\big)\Psi(s,x)dxd^\dagger{\hat{\beta}}_j(s)\\
&&\ \ \ \ \ \ -\sum_{j=1}^{\infty}\int_{t}^{T}\int_{\mathbb{R}^{d}}g_j\big(s,x,u(s,x),(\sigma^*\nabla u)(s,x)\big)\Psi(s,x)dxd^\dagger{\hat{\beta}}_j(s)|^2]\\
&\leq&2E[\ |\sum_{j=1}^{n}\int_{t}^{T}\int_{\mathbb{R}^{d}}\big(g_j\big(s,x,u^n(s,x),(\sigma^*\nabla u^n)(s,x)\big)-g_j\big(s,x,u(s,x),(\sigma^*\nabla u)(s,x)\big)\big)\Psi(s,x)dxd^\dagger{\hat{\beta}}_j(s)|^2]\\
&&+2E[\ |\sum_{j=n+1}^{\infty}\int_{t}^{T}\int_{\mathbb{R}^{d}}g_j\big(s,x,u(s,x),(\sigma^*\nabla u)(s,x)\big)\Psi(s,x)dxd^\dagger{\hat{\beta}}_j(s)|^2]\\
&\leq&C_pE[\sum_{j=1}^{\infty}(C_j+\alpha_j)\int_{t}^{T}\int_{\mathbb{R}^{d}}(|u^n(t,x)-u(t,x)|^2+|(\sigma^*\nabla u^n)(s,x)-(\sigma^*\nabla u)(s,x)|^2)\rho^{-1}(x)dxds]\\
&&+C_pE[\ |\sum_{j=n+1}^{\infty}\int_{t}^{T}\int_{\mathbb{R}^{d}}\big(g_j\big(s,x,u(s,x),(\sigma^*\nabla u)(s,x)\big)-g_j(s,x,0,0)\big)\Psi(s,x)dxd^\dagger{\hat{\beta}}_j(s)|^2]\\
&&+C_pE[\
|\sum_{j=n+1}^{\infty}\int_{t}^{T}\int_{\mathbb{R}^{d}}g_j(s,x,0,0)\Psi(s,x)dxd^\dagger{\hat{\beta}}_j(s)|^2].
\end{eqnarray*}
It is obvious that the first term tends to zero as
$n\rightarrow\infty$. The last two terms can be treated using a
similar method as (\ref{qii6}) and (\ref{qib}).

Therefore $u(t,x)$ satisfies (\ref{qi16}), so is a weak solution of
Eq.(\ref{zhang685}) with $u(T,x)=h(x)$. The uniqueness can be proved
following a similar argument of Theorem 3.1 in Bally and Matoussi
\cite{ba-ma}. $\hfill\diamond$

\section{Infinite horizon BDSDEs}
\setcounter{equation}{0}

We consider the following BDSDE with infinite dimensional noise on
infinite horizon,
\begin{eqnarray}\label{qi30}
{\rm e}^{-Ks}Y_{s}^{t,x}&=&\int_{s}^{\infty}{\rm e}^{-Kr}f(r,X_{r}^{t,x},Y_{r}^{t,x},Z_{r}^{t,x})dr+\int_{s}^{\infty}K{\rm e}^{-Kr}Y_{r}^{t,x}dr\nonumber\\
&&-\int_{s}^{\infty}{\rm e}^{-Kr}
g(r,X_{r}^{t,x},Y_{r}^{t,x},Z_{r}^{t,x})d^\dagger{\hat{B}}_r-\int_{s}^{\infty}{\rm
e}^{-Kr}\langle Z_{r}^{t,x},dW_r\rangle.
\end{eqnarray}
Here
$f:[0,\infty)\times\mathbb{R}^{d}\times\mathbb{R}^1\times\mathbb{R}^{d}{\longrightarrow{\mathbb{R}^1}}$,
$g:[0,\infty)\times\mathbb{R}^{d}\times\mathbb{R}^1\times\mathbb{R}^{d}\longrightarrow
{\mathcal{L}^2_{U_0}(\mathbb{R}^1)}$. Eq.(\ref{qi30}) is equivalent
to
\begin{eqnarray*}\label{qi31}
{\rm e}^{-Ks}Y_{s}^{t,x}&=&\int_{s}^{\infty}{\rm e}^{-Kr}f(r,X_{r}^{t,x},Y_{r}^{t,x},Z_{r}^{t,x})dr+\int_{s}^{\infty}K{\rm e}^{-Kr}Y_{r}^{t,x}dr\nonumber\\
&&-\sum_{j=1}^{\infty}\int_{s}^{\infty}{\rm
e}^{-Kr}g_j(r,X_{r}^{t,x},Y_{r}^{t,x},Z_{r}^{t,x})d^\dagger{\hat{\beta}}_j(r)-\int_{s}^{\infty}{\rm
e}^{-Kr}\langle Z_{r}^{t,x},dW_r\rangle.
\end{eqnarray*}

We assume
\begin{description}
\item[(H.5).] Change $\mathscr{B}_{[0,T]}$ to
$\mathscr{B}_{\mathbb{R}^{+}}$ and $t\in[0,T]$ to $t\geq0$ in (H.2);
\item[(H.6).] Change $\int_{0}^T$ to $\int_{0}^\infty {\rm
e}^{-Ks}$ in (H.3);
\item[(H.7).] There exists a constant $\mu>0$ with $2\mu-K-2C-\sum_{j=1}^{\infty}{C_j}>0$ s.t. for any
$t\geq0$, $Y_1,Y_2\in
L_{\rho}^2({\mathbb{R}^{d}};{\mathbb{R}^{1}})$, $X,Z\in
L_{\rho}^2({\mathbb{R}^{d}};{\mathbb{R}^{d}})$,
\begin{eqnarray*}
&&\int_{\mathbb{R}^d}(Y_1(x)-Y_2(x))\big(f(t,X(x),Y_1(x),Z(x))-f(t,X(x),Y_2(x),Z(x))\big)\rho^{-1}(x)dx\\
&\leq&-\mu\int_{\mathbb{R}^d}{|Y_1(x)-Y_2(x)|}^2\rho^{-1}(x)dx.
\end{eqnarray*}
\end{description}

The main objective of this section is to prove
\begin{thm}\label{qi071} Under Conditions {\rm(H.4)}--{\rm(H.7)},
Eq.(\ref{qi30}) has a unique solution.
\end{thm}
{\em Proof}. \underline{Uniqueness}. Let $(Y_{s}^{t,x},Z_{s}^{t,x})$
and $(\hat{Y}_{s}^{t,x},\hat{Z}_{s}^{t,x})$ be two solutions of
Eq.(\ref{qi30}). Define
\begin{eqnarray*}
&&\bar{Y}_{s}^{t,x}=\hat{Y}_{s}^{t,x}-Y_{s}^{t,x},\ \ \bar{Z}_{s}^{t,x}=\hat{Z}_{s}^{t,x}-Z_{s}^{t,x},\\
&&\bar{f}(s,x)=f(s,X_s^{t,x},\hat{Y}_s^{t,x},\hat{Z}_s^{t,x})-f(s,X_s^{t,x},Y_s^{t,x},Z_s^{t,x}),\\
&&\bar{g}(s,x)=g(s,X_s^{t,x},\hat{Y}_s^{t,x},\hat{Z}_s^{t,x})-g(s,X_s^{t,x},Y_s^{t,x},Z_s^{t,x}),\
\ \ \ \ \ s\geq0.
\end{eqnarray*}
Then for $s\geq0$ and a.e. $x\in\mathbb{R}^d$,
$(Y_{s}^{t,x},Z_{s}^{t,x})$ and
$(\hat{Y}_{s}^{t,x},\hat{Z}_{s}^{t,x})$ satisfy
\begin{eqnarray*}
\left\{\begin{array}{l}
d\bar{Y}_{s}^{t,x}=-\bar{f}(s,x)ds+\sum_{j=1}^{\infty}\bar{g}_j(s,x)d^\dagger{\hat{\beta}}_j(s)+\langle\bar{Z}_{s}^{t,x},dW_s\rangle\\
\lim_{T\longrightarrow\infty}{\rm e}^{-KT}\bar{Y}_{T}^{t,x}=0\ \ \
a.s..
\end{array}\right.
\end{eqnarray*}
For a.e. $x\in\mathbb{R}^d$, applying It$\hat {\rm o}$'s formula for
infinite dimensional noise to ${\rm
e}^{-Ks}{{|\bar{Y}_s^{t,x}|}^2}$, and by Young inequality and
Conditions (H.5), (H.7), we obtain
\begin{eqnarray}\label{qi33}
&&E[\int_{\mathbb{R}^{d}}{\rm e}^{-Ks}{|\bar{Y}_s^{t,x}|}^2\rho^{-1}(x)dx]+(2\mu-K-2C-\sum_{j=1}^{\infty}{C_j})E[\int_{s}^{T}\int_{\mathbb{R}^{d}}{\rm e}^{-Kr}{{|\bar{Y}_r^{t,x}|}^2}\rho^{-1}(x)dxdr]\nonumber\\
&&+({1\over2}-\sum_{j=1}^{\infty}\alpha_j)E[\int_{s}^{T}\int_{\mathbb{R}^{d}}{\rm
e}^{-Kr}|\bar{Z}_r^{t,x}|^2\rho^{-1}(x)dxdr]\nonumber\\
&\leq&E[\int_{\mathbb{R}^{d}}{\rm
e}^{-KT}{|\bar{Y}_T^{t,x}|}^2\rho^{-1}(x)dx].
\end{eqnarray}
Taking $K'>K$ s.t. $2\mu-K'-2C-\sum_{j=1}^{\infty}{C_j}>0$ as well,
we can see that (\ref{qi33}) remains true with $K$ replaced by $K'$.
In particular,
\begin{eqnarray*}
E[\int_{\mathbb{R}^{d}}{\rm
e}^{-K's}{|\bar{Y}_s^{t,x}|}^2\rho^{-1}(x)dx]\leq
E[\int_{\mathbb{R}^{d}}{\rm
e}^{-K'T}{|\bar{Y}_T^{t,x}|}^2\rho^{-1}(x)dx].
\end{eqnarray*}
Therefore, we have
\begin{eqnarray}\label{qi34}
E[\int_{\mathbb{R}^{d}}{\rm
e}^{-K's}{|\bar{Y}_s^{t,x}|}^2\rho^{-1}(x)dx]\leq {\rm
e}^{-(K'-K)T}E[\int_{\mathbb{R}^{d}}{\rm
e}^{-KT}{|\bar{Y}_T^{t,x}|}^2\rho^{-1}(x)dx].
\end{eqnarray}
Since $\hat{Y}_{s}^{t,x}, Y_{s}^{t,x}\in S^{2,-K}\bigcap
M^{2,-K}([0,\infty);L_{\rho}^2({\mathbb{R}^{d}};{\mathbb{R}^{1}}))$,
so
\begin{eqnarray*}
\sup_{T\geq0}E[\int_{\mathbb{R}^{d}}{\rm
e}^{-KT}{|\bar{Y}_T^{t,x}|}^2\rho^{-1}(x)dx]\leq
E[\sup_{T\geq0}\int_{\mathbb{R}^{d}}{\rm
e}^{-KT}(2{|\hat{Y}_T^{t,x}|}^2+2{|{Y}_T^{t,x}|}^2)\rho^{-1}(x)dx]<\infty.
\end{eqnarray*}
Therefore, taking the limit as $T\to \infty$ in (\ref{qi34}), we have
\begin{eqnarray*}
E[\int_{\mathbb{R}^{d}}{\rm
e}^{-K's}{|\bar{Y}_s^{t,x}|}^2\rho^{-1}(x)dx]=0.
\end{eqnarray*}
Then the uniqueness is proved.

\underline{Existence}. For each $n\in\mathbb{N}$, we define a
sequence of BDSDEs (\ref{qi20}) with $h=0$ and $T=n$ and denote it
by Eq.($\ref{qi20}_n$). It is easy to verify that for each n, these
BDSDEs satisfy conditions of Theorem \ref{qi052}. Therefore, for
each n, there exists a $({Y}_s^{t,x,n},{Z}_s^{t,x,n})\in
S^{2,0}([0,n];L_{\rho}^2({\mathbb{R}^{d}};{\mathbb{R}^{1}}))\bigotimes
M^{2,0}([0,n];L_{\rho}^2({\mathbb{R}^{d}};{\mathbb{R}^{d}}))$ which
is equivalent to the space
$S^{2,-K}([0,n];\\L_{\rho}^2({\mathbb{R}^{d}};{\mathbb{R}^{1}}))\bigotimes
M^{2,-K}([0,n];L_{\rho}^2({\mathbb{R}^{d}};{\mathbb{R}^{d}}))$ and
$({Y}_s^{t,x,n},{Z}_s^{t,x,n})$ is the unique solution of
Eq.($\ref{qi20}_n$). That is to say, for an arbitrary $\varphi\in
C_c^{0}(\mathbb{R}^d;\mathbb{R}^1)$, $({Y}_s^{t,x,n},{Z}_s^{t,x,n})$
satisfies
\begin{eqnarray}\label{zhang682}
&&\int_{\mathbb{R}^{d}}{\rm e}^{-Ks}Y_s^{t,x,n}\varphi(x)dx=\int_{s}^{n}\int_{\mathbb{R}^{d}}{\rm e}^{-Kr}f(r,X_{r}^{t,x},Y_{r}^{t,x,n},Z_{r}^{t,x,n})\varphi(x)dxdr\nonumber\\
&&+\int_{s}^{n}\int_{\mathbb{R}^{d}}K{\rm e}^{-Kr}Y_{r}^{t,x,n}\varphi(x)dxdr-\sum_{j=1}^{\infty}\int_{s}^{n}\int_{\mathbb{R}^{d}}{\rm e}^{-Kr}g_j(r,X_{r}^{t,x},Y_{r}^{t,x,n},Z_{r}^{t,x,n})\varphi(x)dxd^\dagger{\hat{\beta}}_j(r)\nonumber\\
&&-\int_{s}^{n}\langle \int_{\mathbb{R}^{d}}{\rm
e}^{-Kr}Z_{r}^{t,x,n}\varphi(x)dx,dW_r\rangle\ \ \ P-{\rm a.s.}.
\end{eqnarray}
Let $(Y_{t}^{n}, Z_{t}^{n})_{t>n}=(0, 0)$, then
$({Y}_s^{t,x,n},{Z}_s^{t,x,n})\in S^{2,-K}\bigcap
M^{2,-K}([0,\infty);L_{\rho}^2({\mathbb{R}^{d}};{\mathbb{R}^{1}}))\bigotimes
M^{2,-K}\\([0,\infty);L_{\rho}^2({\mathbb{R}^{d}};{\mathbb{R}^{d}}))$.
We will prove $({Y}_s^{t,x,n},{Z}_s^{t,x,n})$ is a Cauchy sequence.
For this, let $(Y_{s}^{t,x,m},Z_{s}^{t,x,m})$ and
$({Y}_{s}^{t,x,n},{Z}_{s}^{t,x,n})$ be the solutions of
Eq.($\ref{qi20}_m$) and Eq.($\ref{qi20}_n$) respectively. Without
losing any generality, assume that $m\geq n$, and define
\begin{eqnarray*}
&&\bar{Y}_{s}^{t,x,m,n}={Y}_{s}^{t,x,m}-{Y}_{s}^{t,x,n}, \ \
\bar{Z}_{s}^{t,x,m,n}={Z}_{s}^{t,x,m}-{Z}_{s}^{t,x,n},\\
&&\bar{f}^{m,n}(s,x)=f(s,X_{s}^{t,x},{Y}_{s}^{t,x,m},{Z}_{s}^{t,x,m})-f(s,X_{s}^{t,x},Y_{s}^{t,x,n},Z_{s}^{t,x,n}),\\
&&\bar{g}_j^{m,n}(s,x)=g_j(s,X_{s}^{t,x},{Y}_{s}^{t,x,m},{Z}_{s}^{t,x,m})-g_j(s,X_{s}^{t,x},Y_{s}^{t,x,n},Z_{s}^{t,x,n}),\
\ \ \ \ \ s\geq0.
\end{eqnarray*}
Consider two cases:\\
($\textrm{i}$) When $n\leq s\leq m$,
$\bar{Y}_{s}^{t,x,m,n}={Y}_{s}^{t,x,m}$. Since
$({Y}_s^{t,x,m},{Z}_s^{t,x,m})$ is the solution of
Eq.($\ref{qi20}_m$), we have for any $m\in\mathbb{N}$,
\begin{eqnarray*}
\left\{\begin{array}{l}
d{Y}_{s}^{t,x,m}=-{f}(s,X_{s}^{t,x},{Y}_{s}^{t,x,m},{Z}_{s}^{t,x,m})ds+\sum_{j=1}^{\infty}{g}_j(s,X_{s}^{t,x},{Y}_{s}^{t,x,m},{Z}_{s}^{t,x,m})d^\dagger{\hat{\beta}}_j(s)+\langle{Z}_{s}^{t,x,m},dW_s\rangle\\
{Y}_{m}^{t,x,m}=0\ \ \ {\rm for}\ s\in[0,m],\ {\rm a.e.}\
x\in{\mathbb{R}^{d}},\ {\rm a.s.}.
\end{array}\right.
\end{eqnarray*}
Noting that
$E[\int_{0}^{m}\|g(r,X_{r}^{t,x},Y_{r}^{t,x,m},Z_{r}^{t,x,m})\|^2_{\mathcal{L}^2_{U_0}(\mathbb{R}^{1})}dr]<\infty$
for a.e. $x\in\mathbb{R}^{d}$, we can apply It$\hat {\rm o}$'s
formula to ${\rm e}^{-Kr}|Y_{r}^{t,x,m}|^2$ for a.e.
$x\in\mathbb{R}^d$, then taking integration in
$\mathbb{R}^{d}\times\Omega$, we have
\begin{eqnarray}\label{zhang671}
&&\int_{\mathbb{R}^{d}}{\rm e}^{-Ks}{|{Y}_s^{t,x,m}|}^2\rho^{-1}(x)dx\nonumber\\
&&+\big(2\mu-K-2C-\sum_{j=1}^{\infty}{C_j}-(1+\sum_{j=1}^{\infty}{C_j})\varepsilon\big)\int_{s}^{m}\int_{\mathbb{R}^{d}}{\rm e}^{-Kr}{{|{Y}_r^{t,x,m}|}^2}\rho^{-1}(x)dxdr\nonumber\\
&&+({1\over2}-\sum_{j=1}^{\infty}\alpha_j-\sum_{j=1}^{\infty}\alpha_j\varepsilon)\int_{s}^{m}\int_{\mathbb{R}^{d}}{\rm e}^{-Kr}|{Z}_r^{t,x,m}|^2\rho^{-1}(x)dxdr\nonumber\\
&\leq&C_p\int_{s}^{m}\int_{\mathbb{R}^{d}}{\rm e}^{-Kr}|{f}(r,X_{r}^{t,x},0,0)|^2\rho^{-1}(x)dxdr\nonumber\\
&&+C_p\int_{s}^{m}\int_{\mathbb{R}^{d}}{\rm e}^{-Kr}\sum_{j=1}^{\infty}|{g}_j(r,X_{r}^{t,x},0,0)|^2\rho^{-1}(x)dxdr\nonumber\\
&&-\sum_{j=1}^{\infty}\int_{s}^{m}\int_{\mathbb{R}^{d}}2{\rm e}^{-Kr}{Y}_r^{t,x,m}{g}_j(r,X_{r}^{t,x},{Y}_{r}^{t,x,m},{Z}_{r}^{t,x,m})\rho^{-1}(x)dxd^\dagger{\hat{\beta}}_j(r)\nonumber\\
&&-\int_{s}^{m}\langle\int_{\mathbb{R}^{d}}2{\rm
e}^{-Kr}{Y}_r^{t,x,m}{Z}_{r}^{t,x,m}\rho^{-1}(x)dx,dW_r\rangle.
\end{eqnarray}
Note that the constant $\varepsilon$ can be chosen to be
sufficiently small s.t. all the terms on the left hand side of
(\ref{zhang671}) are positive. By (\ref{zhang671}), as $n$,
$m\longrightarrow\infty$ we have
\begin{eqnarray}\label{zhang674}
&&E[\int_{n}^{m}\int_{\mathbb{R}^{d}}{\rm e}^{-Kr}|Y_{r}^{t,x,m}|^2\rho^{-1}(x)dxdr]+E[\int_{n}^{m}\int_{\mathbb{R}^{d}}{\rm e}^{-Kr}|Z_{r}^{t,x,m}|^2\rho^{-1}(x)dxdr]\nonumber \\
&\leq&C_pE[\int_{n}^{m}\int_{\mathbb{R}^{d}}{\rm
e}^{-Kr}(|{f}(r,X_{r}^{t,x},0,0)|^2+\sum_{j=1}^{\infty}|{g}_j(r,X_{r}^{t,x},0,0)|^2)\rho^{-1}(x)dxdr]\longrightarrow0.
\end{eqnarray}
Note that the right hand side of (\ref{zhang674}) converges to 0
follows from the generalized equivalence of norm principle. Also
using the B-D-G inequality to deal with (\ref{zhang671}) in the
interval $[n,m]$, by (\ref{zhang674}), as $n$,
$m\longrightarrow\infty$ we have
\begin{eqnarray}\label{zhang676}
&&E[\sup_{n\leq s\leq m}\int_{\mathbb{R}^d}{\rm e}^{-Ks}|{Y}_{s}^{t,x,m}|^2\rho^{-1}dx]\nonumber\\
&\leq&C_pE[\int_{n}^{m}\int_{\mathbb{R}^{d}}{\rm e}^{-Kr}(|{f}(r,X_{r}^{t,x},0,0)|^2+\sum_{j=1}^{\infty}|{g}_j(r,X_{r}^{t,x},0,0)|^2)\rho^{-1}(x)dxdr]\nonumber\\
&&+C_pE[\int_{n}^{m}\int_{\mathbb{R}^{d}}{\rm
e}^{-Kr}(|Y_{r}^{t,x,m}|^2+|Z_{r}^{t,x,m}|^2)\rho^{-1}(x)dxdr]\longrightarrow0.
\end{eqnarray}
$(\textrm{i}$$\textrm{i})$ When $0\leq s\leq n$,
\begin{eqnarray*}
\bar{Y}_{s}^{t,x,m,n}={Y}_{n}^{t,x,m}+\int_{s}^{n}\bar{f}^{m,n}(r,x)dr-\sum_{j=1}^{\infty}\int_{s}^{n}\bar{g}^{m,n}_j(r,x)d^\dagger{\hat{\beta}}_j(r)-\int_{s}^{n}\langle\bar{Z}_{r}^{t,x,m,n},dW_r\rangle.
\end{eqnarray*}
Apply It$\hat {\rm o}$'s formula to ${\rm
e}^{-Kr}{{|\bar{Y}_r^{t,x,m,n}|}^2}$ for a.e. $x\in\mathbb{R}^d$,
then
\begin{eqnarray}\label{zhang677}
&&\int_{\mathbb{R}^{d}}{\rm e}^{-Ks}{|\bar{Y}_s^{t,x,m,n}|}^2\rho^{-1}(x)dx+(2\mu-K-2C-\sum_{j=1}^{\infty}C_j)\int_{s}^{n}\int_{\mathbb{R}^{d}}{\rm e}^{-Kr}{{|\bar{Y}_r^{t,x,m,n}|}^2}\rho^{-1}(x)dxdr\nonumber\\
&&+({1\over2}-\sum_{j=1}^{\infty}\alpha_j)\int_{s}^{n}\int_{\mathbb{R}^{d}}{\rm e}^{-Kr}|\bar{Z}_r^{t,x,m,n}|^2\rho^{-1}(x)dxdr\nonumber\\
&\leq&\int_{\mathbb{R}^{d}}{\rm e}^{-Kn}{|{Y}_n^{t,x,m}|}^2\rho^{-1}(x)dx-\sum_{j=1}^{\infty}\int_{s}^{n}\int_{\mathbb{R}^{d}}2{\rm e}^{-Kr}\bar{Y}_r^{t,x,m,n}\bar{g}^{m,n}_j(r,x)\rho^{-1}(x)dxd^\dagger{\hat{\beta}}_j(r)\nonumber\\
&&-\int_{s}^{n}\langle\int_{\mathbb{R}^{d}}2{\rm
e}^{-Kr}\bar{Y}_r^{t,x,m,n}\bar{Z}_{r}^{t,x,m,n}\rho^{-1}(x)dx,dW_r\rangle.
\end{eqnarray}
Taking expectation on both sides of (\ref{zhang677}), as $n$,
$m\longrightarrow\infty$, using (\ref{zhang676}), we have
\begin{eqnarray}\label{zhang678}
&&E[\int_{s}^{n}\int_{\mathbb{R}^{d}}{\rm e}^{-Kr}{{|\bar{Y}_r^{t,x,m,n}|}^2}\rho^{-1}(x)dxdr]+E[\int_{s}^{n}\int_{\mathbb{R}^{d}}{\rm e}^{-Kr}|\bar{Z}_r^{t,x,m,n}|^2\rho^{-1}(x)dxdr]\nonumber\\
&\leq&C_pE[\sup_{n\leq s\leq m}\int_{\mathbb{R}^{d}}{\rm
e}^{-Ks}{|{Y}_s^{t,x,m}|}^2\rho^{-1}(x)dx]\longrightarrow0.
\end{eqnarray}
Also by the B-D-G inequality, (\ref{zhang676}), (\ref{zhang677}) and
(\ref{zhang678}), as $n$, $m\longrightarrow\infty$, we have
\begin{eqnarray*}\label{zhang680}
E[\sup_{0\leq s\leq n}\int_{\mathbb{R}^{d}}{\rm
e}^{-Ks}{{|\bar{Y}_s^{t,x,m,n}|}^2}\rho^{-1}(x)dx]\leq
C_pE[\sup_{n\leq s\leq m}\int_{\mathbb{R}^{d}}{\rm
e}^{-Ks}{|{Y}_s^{t,x,m}|}^2\rho^{-1}(x)dx]\longrightarrow0.
\end{eqnarray*}
Therefore taking a combination of cases ($\textrm{i}$) and
$(\textrm{i}$$\textrm{i})$, as $n$, $m\longrightarrow\infty$, we
have
\begin{eqnarray*}
&&E[\sup_{s\geq0}\int_{\mathbb{R}^d}{\rm e}^{-Ks}|\bar{Y}_{s}^{t,x,m,n}|^2\rho^{-1}(x)dx]+E[\int_{0}^{\infty}\int_{\mathbb{R}^{d}}{\rm e}^{-Kr}|\bar{Y}_{r}^{t,x,m,n}|^2\rho^{-1}(x)dxdr]\\
&&+E[\int_{0}^{\infty}\int_{\mathbb{R}^{d}}{\rm
e}^{-Kr}|\bar{Z}_{r}^{t,x,m,n}|^2\rho^{-1}(x)dxdr]\longrightarrow0.
\end{eqnarray*}
That is to say $({Y}_s^{t,x,n},{Z}_s^{t,x,n})$ is a Cauchy sequence.
Take $({Y}_s^{t,x},{Z}_s^{t,x})$ as the limit of
$({Y}_s^{t,x,n},{Z}_s^{t,x,n})$ in the space $S^{2,-K}\bigcap
M^{2,-K}([0,\infty);L_{\rho}^2({\mathbb{R}^{d}};{\mathbb{R}^{1}}))\bigotimes
M^{2,-K}([0,\infty);L_{\rho}^2({\mathbb{R}^{d}};{\mathbb{R}^{d}}))$
and we will show that $({Y}_s^{t,x},{Z}_s^{t,x})$ is the solution of
Eq.(\ref{qi30}). We only need to verify that for arbitrary
$\varphi\in C_c^{0}(\mathbb{R}^d;\mathbb{R}^1)$,
$({Y}_s^{t,x},{Z}_s^{t,x})$ satisfies ($\ref{qi15}_r$), where
($\ref{qi15}_r$) means a more general form of (\ref{qi15}) with $f$
and $g_j$ also depending on $r\in[0,\infty)$. Since
$({Y}_s^{t,x,n},{Z}_s^{t,x,n})$ satisfies Eq.(\ref{zhang682}), so we
verify that Eq.(\ref{zhang682}) converges to Eq.($\ref{qi15}_r$) in
$L^2(\Omega)$ term by term as $n\longrightarrow\infty$. We only show
the infinite dimensional stochastic integral term:
\begin{eqnarray*}
&&E[\ |\sum_{j=1}^{\infty}\int_{s}^{n}\int_{\mathbb{R}^{d}}{\rm e}^{-Kr}g_j(r,X_{r}^{t,x},Y_r^{t,x,n},Z_r^{t,x,n})\varphi(x)dxd^\dagger{\hat{\beta}}_j(r)\\
&&\ \ \ \ \ \ -\sum_{j=1}^{\infty}\int_{s}^{\infty}\int_{\mathbb{R}^{d}}{\rm e}^{-Kr}g_j(r,X_{r}^{t,x},Y_r^{t,x},Z_r^{t,x})\varphi(x)dxd^\dagger{\hat{\beta}}_j(r)|^2]\\
&\leq&2E[\ |\sum_{j=1}^{\infty}\int_{s}^{n}\int_{\mathbb{R}^{d}}{\rm e}^{-Kr}\big(g_j(r,X_{r}^{t,x},Y_r^{t,x,n},Z_r^{t,x,n})-g_j(r,X_{r}^{t,x},Y_r^{t,x},Z_r^{t,x})\big)\varphi(x)dxd^\dagger{\hat{\beta}}_j(r)|^2]\\
&&+2E[\
|\sum_{j=1}^{\infty}\int_{n}^{\infty}\int_{\mathbb{R}^{d}}{\rm
e}^{-Kr}g_j(r,X_{r}^{t,x},Y_r^{t,x},Z_r^{t,x})\varphi(x)dxd^\dagger{\hat{\beta}}_j(r)|^2].
\end{eqnarray*}
We see that each term in the above formula tends to zero as
$n\rightarrow\infty$ since
\begin{eqnarray*}
&&E[\ |\sum_{j=1}^{\infty}\int_{s}^{n}\int_{\mathbb{R}^{d}}{\rm e}^{-Kr}\big(g_j(r,X_{r}^{t,x},Y_r^{t,x,n},Z_r^{t,x,n})-g_j(r,X_{r}^{t,x},Y_r^{t,x},Z_r^{t,x})\big)\varphi(x)dxd^\dagger{\hat{\beta}}_j(r)|^2]\\
&\leq&C_pE[\int_{0}^{\infty}\int_{\mathbb{R}^{d}}{\rm
e}^{-Kr}(|Y_r^{t,x,n}-Y_r^{t,x}|^2+|Z_r^{t,x,n}-Z_r^{t,x}|^2)\rho^{-1}(x)dxdr]\longrightarrow0,\
\ {\rm as}\ n\rightarrow\infty,
\end{eqnarray*}
and
\begin{eqnarray*}
&&E[\ |\sum_{j=1}^{\infty}\int_{n}^{\infty}\int_{\mathbb{R}^{d}}{\rm e}^{-Kr}g_j(r,X_{r}^{t,x},Y_r^{t,x},Z_r^{t,x})\varphi(x)dxd^\dagger{\hat{\beta}}_j(r)|^2]\\
&\leq&C_pE[\int_{n}^{\infty}\int_{\mathbb{R}^{d}}{\rm e}^{-Kr}(|Y_r^{t,x}|^2+|Z_r^{t,x}|^2)\rho^{-1}(x)dxdr]\\
&&+C_p\int_{n}^{\infty}\int_{\mathbb{R}^{d}}\sum_{j=1}^{\infty}{\rm
e}^{-Kr}|g_j(r,x,0,0)|^2\rho^{-1}(x)dxdr\longrightarrow0,\ \ {\rm
as}\ n\rightarrow\infty.
\end{eqnarray*}
That is to say ${(Y_s^{t,x}, Z_s^{t,x})}_{s\geq0}$ satisfies
Eq.($\ref{qi15}_r$). The proof of Theorem \ref{qi071} is completed.
$\hfill\diamond$\\

By similar method as in the proof of existence part case (i) in
Theorem \ref{qi071}, we have the following estimation:
\begin{prop}\label{qi072}
Let $({Y}_s^{t,x,n},{Z}_s^{t,x,n})$ be the solution of
Eq.{\rm($\ref{qi20}_n$)}, then under the conditions of Theorem
\ref{qi071},
\begin{eqnarray*}
&&\sup_nE[\sup_{s\geq0}\int_{\mathbb{R}^d}{\rm
e}^{-Ks}|Y_s^{t,x,n}(x)|^2\rho^{-1}(x)dx]+\sup_nE[\int_{0}^{\infty}\int_{\mathbb{R}^{d}}{\rm
e}^{-Kr}|Y_r^{t,x,n}(x)|^2\rho^{-1}(x)dxdr]\nonumber\\
&&+\sup_nE[\int_{0}^{\infty}\int_{\mathbb{R}^{d}}{\rm
e}^{-Kr}|Z_r^{t,x,n}(x)|^2\rho^{-1}(x)dxdr]<\infty.
\end{eqnarray*}
\end{prop}

\section{The continuity of the solution of the infinite horizon BDSDEs as the solution of the corresponding SPDEs}
\setcounter{equation}{0}

Now we study BDSDE (\ref{qi13}), a simpler form of Eq.(\ref{qi30}).

{\em Proof of Theorem \ref{qi043}}. Since conditions here are
stronger than those in Theorem \ref{qi071}, so there exists a unique
solution $(Y_{s}^{t,x},Z_{s}^{t,x})$. We only need to prove
$E[\sup_{s\geq0}\int_{\mathbb{R}^{d}}{\rm
e}^{{-{pK}}s}{{|{Y}_s^{t,x}|}^p}\rho^{-1}(x)dx]<\infty$. Let
$\varphi_{N,p}(x)=x^{p\over2}I_{\{0\leq
x<N\}}+{p\over2}N^{{p-2}\over2}(x-N)I_{\{x\geq N\}}$. We apply
generalized It$\hat {\rm o}$'s formula to ${\rm
e}^{{-{pK}}r}\varphi_{N,p}\big(\psi_M(Y_{r}^{t,x})\big)$ for a.e.
$x\in{\mathbb{R}^{d}}$ to have the following estimation
\begin{eqnarray}\label{zhang690}
&&{\rm e}^{{-{pK}}s}\varphi_{N,p}\big(\psi_M(Y_{s}^{t,x})\big)-{pK}\int_{s}^{T}{\rm e}^{{-{pK}}r}\varphi_{N,p}\big(\psi_M(Y_{r}^{t,x})\big)dr\nonumber\\
&&+{1\over2}\int_{s}^{T}{\rm e}^{{-{pK}}r}\varphi^{''}_{N,p}\big(\psi_M(Y_{r}^{t,x})\big)|\psi_M^{'}(Y_{r}^{t,x})|^2|Z_{r}^{t,x}|^2dr\nonumber\\
&&+\int_{s}^{T}{\rm e}^{{-{pK}}r}\varphi^{'}_{N,p}\big(\psi_M(Y_{r}^{t,x})\big)I_{\{-M\leq{Y}_{r}^{t,x}<M\}}|{Z}_{r}^{t,x}|^2dr\nonumber\\
&\leq&{\rm e}^{{-{pK}}T}\varphi_{N,p}\big(\psi_M(Y_{T}^{t,x})\big)+\int_{s}^{T}{\rm e}^{{-{pK}}r}\varphi^{'}_{N,p}\big(\psi_M(Y_{r}^{t,x})\big)\psi_M^{'}(Y_{r}^{t,x}){f}(X_{r}^{t,x},Y_{r}^{t,x},Z_{r}^{t,x})dr\nonumber\\
&&+\int_{s}^{T}{\rm e}^{{-{pK}}r}\varphi^{'}_{N,p}\big(\psi_M(Y_{r}^{t,x})\big)I_{\{-M\leq{Y}_{r}^{t,x}<M\}}\sum_{j=1}^{\infty}|{g}_j(X_{r}^{t,x},Y_{r}^{t,x},Z_{r}^{t,x})|^2dr\nonumber\\
&&+{1\over2}\int_{s}^{T}{\rm e}^{{-{pK}}r}\varphi^{''}_{N,p}\big(\psi_M(Y_{r}^{t,x})\big)|\psi_M^{'}(Y_{r}^{t,x})|^2\sum_{j=1}^{\infty}|{g}_j(X_{r}^{t,x},Y_{r}^{t,x},Z_{r}^{t,x})|^2dr\nonumber\\
&&-\sum_{j=1}^{\infty}\int_{s}^{T}{\rm e}^{{-{pK}}r}\varphi^{'}_{N,p}\big(\psi_M(Y_{r}^{t,x})\big)\psi_M^{'}(Y_{r}^{t,x}){g}_j(X_{r}^{t,x},{Y}_{r}^{t,x},{Z}_{r}^{t,x})d^\dagger{\hat{\beta}}_j(r)\nonumber\\
&&-\int_{s}^{T}\langle{\rm
e}^{{-{pK}}r}\varphi^{'}_{N,p}\big(\psi_M(Y_{r}^{t,x})\big)\psi_M^{'}(Y_{r}^{t,x}){Z}_{r}^{t,x},dW_r\rangle.
\end{eqnarray}
Note that $\lim_{T\rightarrow\infty}{\rm
e}^{{-{pK}}T}\varphi_{N,p}\big(\psi_M(Y_{T}^{t,x})\big)=0$, so after
taking limit as $T\to \infty$, we take the integration on
$\Omega\times\mathbb{R}^{d}$. As
$(Y_{\cdot}^{t,\cdot},Z_{\cdot}^{t,\cdot})\in S^{2,-K}\bigcap
M^{2,-K}([0,\infty);L_{\rho}^2({\mathbb{R}^{d}};{\mathbb{R}^{1}}))\bigotimes
M^{2,-K}([0,\infty);L_{\rho}^2({\mathbb{R}^{d}};{\mathbb{R}^{d}}))$
and
$\varphi^{'}_{N,p}\big(\psi_M(Y_{r}^{t,x})\big)\psi_M^{'}(Y_{r}^{t,x})$
is bounded, we can use the stochastic Fubini theorem and all the
stochastic integrals have zero expectation. Using Conditions
(A.1)$'$--(A.4)$'$, and taking the limit as $M\to \infty$ first,
then the limit as $N\to \infty$, by the monotone convergence theorem,
we have
\begin{eqnarray}\label{zhang691}
&&\big(p\mu-{pK}-pC-{{p(p-1)}\over2}\sum_{j=1}^{\infty}C_j-(3+{{p(p-1)}\over2}\sum_{j=1}^{\infty}{C_j})\varepsilon\big)E[\int_{s}^{\infty}\int_{\mathbb{R}^{d}}{\rm e}^{{-{pK}}r}{|{Y}_r^{t,x}|}^p\rho^{-1}(x)dxdr]\nonumber\\
&&+{p\over4}\big(2p-3-(2p-2)\sum_{j=1}^{\infty}\alpha_j-(2p-2)\sum_{j=1}^{\infty}\alpha_j\varepsilon\big)E[\int_{s}^{\infty}\int_{\mathbb{R}^{d}}{\rm e}^{{-{pK}}r}{{|{Y}_r^{t,x}|}^{p-2}}|{Z}_{r}^{t,x}|^2\rho^{-1}(x)dxdr]\nonumber\\
&\leq&C_p\int_{\mathbb{R}^{d}}|{f}(x,0,0)|^p\rho^{-1}(x)dx+C_p\int_{\mathbb{R}^{d}}\sum_{j=1}^{\infty}|{g}_j(x,0,0)|^p\rho^{-1}(x)dx]<\infty.
\end{eqnarray}
Note that the constant $\varepsilon$ can be chosen to be
sufficiently small s.t. all the terms on the left hand side of
(\ref{zhang691}) are positive. Also by the B-D-G inequality,
Cauchy-Schwartz inequality and Young inequality, from
(\ref{zhang690}) we have
\begin{eqnarray*}\label{zhang692}
&&E[\sup_{s\geq0}\int_{\mathbb{R}^{d}}{\rm e}^{{-{pK}}s}{{|{Y}_s^{t,x}|}^p}\rho^{-1}(x)dx]\nonumber\\
&\leq&C_p\int_{\mathbb{R}^{d}}|{f}(x,0,0)|^p\rho^{-1}(x)dx+C_p\int_{\mathbb{R}^{d}}\sum_{j=1}^{\infty}|{g}_j(x,0,0)|^p\rho^{-1}(x)dx\\
&&+C_pE[\int_{0}^{\infty}\int_{\mathbb{R}^{d}}{\rm
e}^{{-{pK}}r}{{|{Y}_r^{t,x}|}^{p-2}}|{Z}_{r}^{t,x}|^2\rho^{-1}(x)dxdr]
+C_pE[\int_{0}^{\infty}\int_{\mathbb{R}^{d}}{\rm
e}^{{-{pK}}r}|Y_{r}^{t,x}|^p\rho^{-1}(x)dxdr].\nonumber
\end{eqnarray*}
So by (\ref{zhang691}), Theorem \ref{qi043} is proved.
$\hfill\diamond$
\\

We need to prove two lemmas before giving a proof of Theorem
\ref{qi044}.
\begin{lem}\label{qi081} Under Condition {\rm(A.3)$'$}, for arbitrary $T>0$, $t,t'\in[0,T]$,
\begin{eqnarray*}
E[\int_{0}^{\infty}\int_{\mathbb{R}^{d}}{\rm
e}^{-Kr}|X_r^{t',x}-X_r^{t,x}|^p\rho^{-1}(x)dxdr]&\leq&C_p|t'-t|^{p\over2}\
\ \ \rm{a.s.}.
\end{eqnarray*}
\end{lem}
{\em Proof}. It is not difficult to deduce from Lemma 4.5.6 in
\cite{ku2}, so we omit the proof. $\hfill\diamond$

\begin{lem}\label{qi082} Under Conditions {\rm(A.1)$'$}--{\rm(A.4)$'$}, for arbitrary $T>0$, $t,t'\in[0,T]$, let $(Y_{s}^{t',x})_{s\geq0}$, $(Y_{s}^{t,x})_{s\geq0}$ be the solutions of Eq.(\ref{qi30}), then
\begin{eqnarray*}
E[\sup_{s\geq0}\int_{\mathbb{R}^{d}}{\rm
e}^{-{pK}s}|Y_{s}^{t',x}-Y_{s}^{t,x}|^p\rho^{-1}(x)dx]\leq
C_p|t'-t|^{p\over2}.
\end{eqnarray*}
\end{lem}
{\em Proof}. Let
\begin{eqnarray*}
&&\bar{Y}_s=Y_{s}^{t',x}-Y_{s}^{t,x},\ \
\ \bar{Z}_s=Z_{s}^{t',x}-Z_{s}^{t,x},\\
&&\bar{f}(s)=f(X_{s}^{t',x},Y_{s}^{t',x},Z_{s}^{t',x})-f(X_{s}^{t,x},Y_{s}^{t,x},Z_{s}^{t,x}),\\
&&\bar{g}_j(s)=g_j(X_{s}^{t',x},Y_{s}^{t',x},Z_{s}^{t',x})-g_j(X_{s}^{t,x},Y_{s}^{t,x},Z_{s}^{t,x}),\
\ \ \ \ \ s\geq0.
\end{eqnarray*}
Then
\begin{eqnarray*}
\left\{\begin{array}{l}
d\bar{Y}_s=-\bar{f}(s)ds+\sum_{j=1}^{\infty}\bar{g}_j(s)d^\dagger{\hat{\beta}}_j(s)+\langle\bar{Z}_s,dW_s\rangle\\
\lim_{T\rightarrow\infty}{\rm e}^{-KT} \bar{Y}_T=0\ \ \ {\rm for}\
{\rm a.e.}\ x\in\mathbb{R}^d\ {\rm a.s.}.
\end{array}\right.
\end{eqnarray*}
First note that from Theorem \ref{qi043}, we know
$E[\sup_{s\geq0}\int_{\mathbb{R}^{d}}{\rm
e}^{-{pK}s}|\bar{Y}_s|^p\rho^{-1}(x)dx]<\infty$. Applying It$\hat
{\rm o}$'s formula to ${\rm e}^{-{pK}r}|\bar{Y}_r|^p$ for a.e.
$x\in\mathbb{R}^d$ (we leave out procedure of localization as in
(\ref{zhang690}) for simplicity) and taking integration on
$\mathbb{R}^{d}$, we have
\begin{eqnarray}\label{zhang695}
&&\int_{\mathbb{R}^{d}}{\rm e}^{-{pK}s}|\bar{Y}_s|^p\rho^{-1}(x)dx\nonumber\\
&&+\big(p\mu-{pK}-pC-{{p(p-1)}\over2}\sum_{j=1}^{\infty}C_j-3\varepsilon\big)\int_{s}^{\infty}\int_{\mathbb{R}^{d}}{\rm e}^{-{pK}r}|\bar{Y}_{r}|^p\rho^{-1}(x)dxdr\nonumber\\
&&+{p\over4}\big(2p-3-(2p-2)\sum_{j=1}^{\infty}\alpha_j\big)\int_{s}^{\infty}\int_{\mathbb{R}^{d}}{\rm e}^{-{pK}r}|\bar{Y}_{r}|^{p-2}|\bar{Z}_r|^2\rho^{-1}(x)dxdr\nonumber\\
&\leq&C_p\int_{s}^{\infty}\int_{\mathbb{R}^{d}}{\rm e}^{-{pK}r}|\bar{X}_{r}|^p\rho^{-1}(x)dxdr-p\sum_{j=1}^{\infty}\int_{s}^{\infty}\int_{\mathbb{R}^{d}}{\rm e}^{-{pK}r}|\bar{Y}_{r}|^{p-2}\bar{Y}_{r}\bar{g}_j(r)\rho^{-1}(x)dxd^\dagger{\hat{\beta}}_j(r)\nonumber\\
&&-p\int_{s}^{\infty}\langle\int_{\mathbb{R}^{d}}{\rm
e}^{-{pK}r}|\bar{Y}_{r}|^{p-2}\bar{Y}_{r}\bar{Z}_r\rho^{-1}(x)dx,dW_r\rangle.
\end{eqnarray}
Note that the constant $\varepsilon$ can be chosen to be
sufficiently small s.t. all the terms on the left hand side of
(\ref{zhang695}) are positive. Taking integration on $\Omega$ on
both sides of (\ref{zhang695}), by Lemma \ref{qi081} we have
\begin{eqnarray}\label{zhang696}
&&E[\int_{s}^{\infty}\int_{\mathbb{R}^{d}}{\rm e}^{-{pK}r}|\bar{Y}_{r}|^p\rho^{-1}(x)dxdr]+E[\int_{s}^{\infty}\int_{\mathbb{R}^{d}}{\rm e}^{-{pK}r}|\bar{Y}_{r}|^{p-2}|\bar{Z}_r|^2\rho^{-1}(x)dxdr]\nonumber\\
&\leq&C_pE[\int_{s}^{\infty}\int_{\mathbb{R}^{d}}{\rm e}^{-{pK}r}|\bar{X}_{r}|^p\rho^{-1}(x)dxdr]\nonumber\\
&\leq&C_p|t'-t|^{p\over2}.
\end{eqnarray}
Also by the B-D-G inequality, from (\ref{zhang695}) and
(\ref{zhang696}), we have
\begin{eqnarray*}
&&E[\sup_{s\geq0}\int_{\mathbb{R}^{d}}{\rm e}^{-{pK}s}|\bar{Y}_s|^p\rho^{-1}(x)dx]\\
&\leq&C_pE[\int_{0}^{\infty}\int_{\mathbb{R}^{d}}{\rm e}^{-{pK}r}|\bar{X}_{r}|^p\rho^{-1}(x)dxdr]+C_pE[\int_{0}^{\infty}\int_{\mathbb{R}^{d}}{\rm e}^{-{pK}r}|\bar{Y}_r|^p\rho^{-1}(x)dxdr]\nonumber\\
&&+C_pE[\int_{0}^{\infty}\int_{\mathbb{R}^{d}}{\rm e}^{-{pK}r}{{|\bar{Y}_r|}^{p-2}}|\bar{Z}_r|^2\rho^{-1}(x)dxdr]\\
&\leq&C_p|t'-t|^{p\over2}.
\end{eqnarray*}
$\hfill\diamond$

{\em Proof of Theorem \ref{qi044}}. By Lemma \ref{qi082}, we have
\begin{eqnarray*}
&&E([\sup_{s\geq0}\int_{\mathbb{R}^{d}}{\rm e}^{-{2K}s}|Y_{s}^{t',x}-Y_{s}^{t,x}|^2\rho^{-1}(x)dx])^{p\over2}\nonumber\\
&\leq&C_pE[\sup_{s\geq0}\int_{\mathbb{R}^{d}}{\rm e}^{-{pK}r}|Y_{s}^{t',x}-Y_{s}^{t,x}|^p\rho^{-1}(x)dx]\big(\int_{\mathbb{R}^{d}}\rho^{-1}(x)dx\big)^{{p-2}\over2}\nonumber\\
&\leq&C_p|t'-t|^{p\over2}.
\end{eqnarray*}
Noting $p>2$, by the Kolmogorov continuity theorem (see \cite{ku2}),
we have $t\longrightarrow Y_{s}^{t,x}$ is a.s. continuous for
$t\in[0,T]$ under the norm $(\sup_{s\geq0}\int_{\mathbb{R}^{d}}{\rm
e}^{-{2K}s}|\cdot|^2\rho^{-1}(x)dx)^{1\over2}$. Without losing any
generality, assume that $t'\geq t$. Then we can see
\begin{eqnarray*}
\lim_{t'\rightarrow t}(\int_{\mathbb{R}^{d}}{\rm
e}^{-{2K}t'}|Y_{t'}^{t',x}-Y_{t'}^{t,x}|^2\rho^{-1}(x)dx)^{1\over2}\leq\lim_{t'\rightarrow
t}(\sup_{s\geq0}\int_{\mathbb{R}^{d}}{\rm
e}^{-{2K}s}|Y_{s}^{t',x}-Y_{s}^{t,x}|^2\rho^{-1}(x)dx)^{1\over2}=0\
\ \ \ \rm{a.s.}.
\end{eqnarray*}
Notice $t'\in[0,T]$, so
\begin{eqnarray}\label{zhang706}
\lim_{t'\rightarrow
t}(\int_{\mathbb{R}^{d}}|Y_{t'}^{t',x}-Y_{t'}^{t,x}|^2\rho^{-1}(x)dx)^{1\over2}=0\
\ \ \ \rm{a.s.}.
\end{eqnarray}
Since $Y_{\cdot}^{t,\cdot}\in
S^{2,-K}([0,\infty);L_{\rho}^2({\mathbb{R}^{d}};{\mathbb{R}^{1}}))$,
$Y_{t'}^{t,\cdot}$ is continuous w.r.t. $t'$ in
$L_{\rho}^2(\mathbb{R}^d;\mathbb{R}^1)$. That is to say for each
$t$,
\begin{eqnarray}\label{zhang705}
\lim_{t'\rightarrow
t}(\int_{\mathbb{R}^{d}}|{Y}_{t'}^{t,x}-{Y}_{t}^{t,x}|^2\rho^{-1}(x)dx)^{1\over2}=0\
\ \ \ \rm{a.s.}.
\end{eqnarray}
Now by (\ref{zhang706}) and (\ref{zhang705})
\begin{eqnarray*}
&&\lim_{t'\rightarrow t}(\int_{\mathbb{R}^{d}}|{Y}_{t'}^{t',x}-{Y}_{t}^{t,x}|^2\rho^{-1}(x)dx)^{1\over2}\\
&\leq&\lim_{t'\rightarrow t}(\int_{\mathbb{R}^{d}}|{Y}_{t'}^{t',x}-{Y}_{t'}^{t,x}|^2\rho^{-1}(x)dx)^{1\over2}+\lim_{t'\rightarrow t}(\int_{\mathbb{R}^{d}}|{Y}_{t'}^{t,x}-{Y}_{t}^{t,x}|^2\rho^{-1}(x)dx)^{1\over2}\\
&=&0\ \ \ \ \rm{a.s.}.
\end{eqnarray*}
For arbitrary $T>0$, $0\leq t\leq T$, define
$u(t,\cdot)=Y_t^{t,\cdot}$, then $u(t,\cdot)$ is a.s. continuous
w.r.t. $t$ in $L_{\rho}^2(\mathbb{R}^d;\mathbb{R}^1)$. Since
$Y_{\cdot}^{t,\cdot}\in
S^{2,-K}([0,\infty);L_{\rho}^2({\mathbb{R}^{d}};{\mathbb{R}^{1}}))$,
$Y_T^{T,x}$ is
$\mathscr{F}_{T,\infty}^{\hat{B}}\otimes\mathscr{B}_{\mathbb{R}^{d}}$
measurable and
$E[\int_{\mathbb{R}^{d}}|Y_T^{T,x}|^2\rho^{-1}(x)dx]<\infty$. It
follows that Condition (H.1) is satisfied. Moreover, Conditions
{\rm(A.1)$'$}--{\rm(A.3)$'$} are stronger than Conditions
(H.2)--(H.4), so by Theorem \ref{qi062}, $u(t,x)$ is a weak solution
of Eq.(\ref{zhang685}). Theorem \ref{qi044} is proved.
$\hfill\diamond$
\\
\\
{\bf Acknowledgements}. It is our great pleasure to thank Z. Brzezniak, Z. Ma, E.
Pardoux, S. Peng and T. Zhang for useful conversations.

\end{document}